\documentclass[final,leqno,oneeqnum,onefignum,onetabnum]{siamltex1213}

\usepackage{moreverb}

\newif\ifpdf
\pdffalse

\usepackage{tikz,environ}
\usetikzlibrary{arrows, positioning, calc, shapes, backgrounds, decorations.markings}

\usepackage{amsmath,amssymb,amsbsy}
\usepackage{bm,multirow}

\ifpdf
\else
\usepackage{psfrag}
\fi

\usepackage{booktabs}
\usepackage{mathrsfs}
\usepackage{algorithm,algorithmic}
\usepackage{setspace}

\ifpdf
\usepackage{epstopdf}
\fi

\graphicspath{{../ps/}}

\makeatletter
\newsavebox{\measure@tikzpicture}
\NewEnviron{scaletikzpicturetowidth}[1]{%
  \def\tikz@width{#1}%
  \def\tikzscale{1}\begin{lrbox}{\measure@tikzpicture}%
  \BODY
  \end{lrbox}%
  \pgfmathparse{#1/\wd\measure@tikzpicture}%
  \edef\tikzscale{\pgfmathresult}%
  \BODY
}
\makeatother

\tikzset{->-/.style={decoration={markings,mark=at position #1 with {\arrow{>}}},postaction={decorate}}}
\tikzset{-<-/.style={decoration={markings,mark=at position #1 with {\arrow{<}}},postaction={decorate}}}

\def\equationautorefname~#1\null{%
  Eq.~(#1)\null
}
\def\algorithmautorefname~#1\null{%
  Algorithm~#1\null
}
\def\sectionautorefname~#1\null{%
  Section~#1\null
}
\def\subsectionautorefname~#1\null{%
  Subsection~#1\null
}

\title{The inverse fast multipole method: using a fast approximate direct solver as a preconditioner for dense linear systems}

\author{Pieter Coulier\footnotemark[2]\ \footnotemark[3] \and Hadi Pouransari\footnotemark[3] \and Eric Darve\footnotemark[3]}

\newcommand\mO{\mathcal{O}}
\newcommand\mH{\mathcal{H}}

\begin{document}
\maketitle
\slugger{sisc}{xxxx}{xx}{x}{x--x}

\renewcommand{\thefootnote}{\fnsymbol{footnote}}

\footnotetext[2]{KU~Leuven, Department of Civil Engineering, Kasteelpark Arenberg 40, 3001 Leuven, Belgium (\email{pieter.coulier@bwk.kuleuven.be}).}
\footnotetext[3]{Stanford University, Department of Mechanical Engineering, 496 Lomita Mall, 94305 Stanford, CA, USA (\email{[pcoulier,hadip,darve]@stanford.edu}).}

\renewcommand{\thefootnote}{\arabic{footnote}}

\begin{abstract}
Although some preconditioners are available for solving dense linear systems, there are still many matrices for which preconditioners are lacking, in particular in cases where the size of the matrix $N$ becomes very large. Examples of preconditioners include ILU preconditioners that sparsify the matrix based on some threshold, algebraic multigrid, and specialized preconditioners, e.g., Calder\'on and other analytical approximation methods when available. Despite these methods, there remains a great need to develop general purpose preconditioners whose cost scales well with the matrix size $N$. In this paper, we propose a preconditioner, with broad applicability and with cost $\mO(N)$, for dense matrices, when the matrix is given by a smooth kernel. Extending the method using the same framework to general $\mH^2$-matrices (i.e., algebraic instead of defined in terms of an analytical kernel) is relatively straightforward, but this won't be discussed here. These preconditioners have a controlled accuracy (e.g., machine accuracy can be achieved if needed) and scale linearly with $N$. They are based on an approximate direct solve of the system. The linear scaling of the algorithm is achieved by means of two key ideas. First, the $\mH^2$-structure of the dense matrix is exploited to obtain an extended sparse system of equations. Second, fill-ins arising when performing the elimination of the latter are compressed as low-rank matrices if they correspond to well-separated interactions. This ensures that the sparsity pattern of the extended sparse matrix is preserved throughout the elimination, hence resulting in a very efficient algorithm with $\mO(N \log(1/\varepsilon)^2 )$ computational cost and $\mO(N \log 1/\varepsilon )$ memory requirement, for an error tolerance $0 < \varepsilon < 1$. The solver is inexact, although the error can be controlled and made as small as needed. These solvers are related to ILU in the sense that the fill-in is controlled. However, in ILU, most of the fill-in (i.e., below a certain tolerance) is simply discarded, whereas here it is approximated using low-rank blocks, with a prescribed tolerance. Numerical examples are discussed to demonstrate the linear scaling of the method and to illustrate its effectiveness as a preconditioner.
\end{abstract}



\pagestyle{myheadings}
\thispagestyle{plain}
\markboth{PIETER COULIER, HADI POURANSARI, AND ERIC DARVE}{THE INVERSE FAST MULTIPOLE METHOD}

\section{Introduction}
Large linear systems of equations $\mathbf{A}\mathbf{x} = \mathbf{b}$ involving dense matrices $\mathbf{A}$ arise in a variety of applications in science and engineering (e.g., in boundary element methods~\cite{bonn95a,hsia06a}, Kalman filtering~\cite{kalmij63a}, radial basis function interpolation~\cite{buhm03a,gume07a}, etc.). Storing a dense matrix entrywise requires a quadratic amount of memory $\mO(N^2)$ (where $N$ indicates the matrix size), while solving the corresponding system of equations with a direct solver involves $\mO(N^3)$ flops. In order to alleviate these stringent requirements, fast methods have been developed in the past decades, including the fast multipole method (FMM)~\cite{darv00a,gree87a,nish02a,rokh85a,ying04a}, the panel clustering technique~\cite{hack89a,saut98a}, and methods based on hierarchical matrices~\cite{bebe08a,borm03a,hack00a,hack02a,hack00b}. These methods rely on the observation that many dense matrices can be represented in a hierarchical format, i.e., the original matrix can be subdivided into a hierarchy of smaller block matrices and most of these blocks are replaceable by low-rank approximations. Several categories of hierarchical matrix representations can be distinguished. In the simplest case, all off-diagonal blocks are assumed to be low-rank, which is referred to as a weak admissibility criterion~\cite{hack04a}. If the matrix entries can be interpreted as interactions between points in a physical domain, this criterion implies that all the interactions of a certain cluster of points with any other cluster (except itself) are low-rank. This assumption leads to hierarchically off-diagonal low-rank (HODLR) matrices~\cite{amin14a} and hierarchically semi-separable (HSS) matrices~\cite{chan06a,shen07b}. In the latter, a nested approach is used (i.e., the low-rank basis at a certain hierarchical level is constructed using the low-rank basis at the child level), while this is not the case in the former. If a more general strong admissibility criterion is adopted, $\mH$- and $\mH^2$-matrices are obtained for the non-nested and nested cases, respectively~\cite{bebe08a,borm03a}\footnote{There is a confusing point of terminology associated with strong and weak, which cannot really be removed. If we use a strong admissibility criterion, we are considering that all well-separated blocks are low-rank. This therefore must correspond to weak hierarchical matrices or w$\mH$ (e.g., FMM matrices). Similarly, if we use a weak admissibility criterion, we get strong hierarchical matrices or s$\mH$ (e.g., HODLR and HSS). This terminology is imposed by the fact that we must have s$\mH$ $\subset$ w$\mH$.}. In these matrices, only interactions between non-neighboring clusters of points are assumed to be low-rank. The FMM can be considered as a subcategory of $\mH^2$-matrices~\cite{bebe08a}. The reader is referred to~\cite{bebe08a,borm03a} for a detailed overview of the various hierarchical structures. See also \autoref{classification} for additional discussion of this classification.

The use of the aforementioned hierarchical matrix representations facilitates matrix-vector multiplications, hence paving the way for iterative solvers such as Krylov subspace methods for solving the linear system $\mathbf{A}\mathbf{x} = \mathbf{b}$~\cite{saad81a}. The precarious convergence behavior --- and thus the number of iterations --- is a major disadvantage of iterative solvers. This necessitates the incorporation of a preconditioning strategy to improve the latter, which is often a difficult and case-dependent problem. Typical preconditioners for dense matrices include the incomplete LU factorization (with thresholding)~\cite{lee03a}, methods based on the sparse approximate inverse~\cite{carp04a}, geometric and algebraic multigrid methods~\cite{lang05a,trot00a}, as well as analytical preconditioning techniques such as Calder\'on preconditioners~\cite{chri02a,darb13a}. Although these techniques are usually well suited for the specific problem they have been designed for, their performance as robust and efficient black-box preconditioners for general problems remains rather limited. Another disadvantage of iterative solvers is that they do not lead to a factorization of $\mathbf{A}$ that can easily be re-used, implying that the algorithm needs to be restarted if the right hand side $\mathbf{b}$ of the system is modified. Direct solvers, on the other hand, are known to be more robust and are able to cope with multiple right hand sides, but at a high computational cost.

Significant progress has been made in recent years regarding the development of fast direct solvers for hierarchical matrices, especially for the aforementioned category of HODLR and HSS matrices. For example, direct $\mO(N)$ algorithms for HSS matrices have been presented by, among others, Martinsson and Rokhlin~\cite{mart05b}, Chandrasekaran et al.~\cite{chan06a}, Xia et al.~\cite{xia10a,xia10b}, and Gillman et al.~\cite{gill14a}, while fast solvers for HODLR matrices can be found in~\cite{amin14a,kong11a}. These methods can be applied as high accuracy direct solvers but are often more effective if employed as preconditioners in an iterative solver. The framework used for HODLR and HSS matrices is appropriate for one-dimensional (1D) applications. For higher dimensional problems (2D or 3D), the assumption that all off-diagonal blocks are low-rank does not hold in that case, that is the rank typically grows with $N$. For 2D problems (points distributed on a 2D manifold or 2D plane), the rank grows like $\mO(\sqrt{N})$~\cite{ambi13b}, hence jeopardizing the computational efficiency of the algorithms. It is worth noting that the approach presented in~\cite{mart05b} has recently been extended by Corona et al.~\cite{coro15a} to obtain a direct solver for HSS matrices that scales even for 2D problems as $\mO(N)$. This improvement in complexity is achieved by introducing an additional level of compression, which also results in an improved memory efficiency of the algorithm. Extensions of the method towards surfaces in 3D are currently underway~\cite{coro15a}. Another recent development towards fast direct solvers is the hierarchical interpolative factorization recently introduced by Ho and Ying~\cite{ho15b,ho15a}, which combines ideas from \cite{coro15a} and \cite{ho12a}. In order to obtain a direct solver with a linear complexity $\mO(N)$ for general 2D and 3D problems, the most general framework of $\mH$- and $\mH^2$-matrices seems to be very well suited as well.

In this paper, the inverse fast multipole method (IFMM) is presented as a fast direct solver for dense $\mH^2$-matrices, where ``direct'' refers to the use of an approximate LU factorization~\cite{ambi14a}. We focus on matrices given by a smooth kernel. The solver is inexact, although the error can be controlled and made as small as needed. A low accuracy solver is hence applicable as a robust and efficient preconditioner for iterative schemes. The IFMM relies on two key ideas. First, the dense $\mH^2$-matrix is converted to an extended sparse matrix through the introduction of auxiliary variables. Second, fill-ins\footnote{The term {\em fill-in} is used here to denote the non-zero matrix entries.} arising during the elimination\footnote{The term {\em elimination} refers to removing variables from a system of equations and updating the coefficients of the remaining variables appropriately.} of the latter are compressed as low-rank matrices if they correspond to well-separated (i.e., non-neighboring) interactions. As a result, the sparsity of the extended sparse matrix is maintained throughout the elimination. Specifically, there is an upper bound in $\mO(N)$ on the number of non-zero entries, at all steps during the elimination. This results in a very efficient algorithm with linear complexity. To the authors' knowledge, this is the first direct solver that achieves this complexity for the most general category of hierarchical matrices~\cite{ambi14a}. This paper only considers dense matrices, but similar ideas can be employed to develop a fast hierarchical solver with linear complexity for sparse matrices. Such a solver is described in a companion paper~\cite{pour15b}.

The proposed strategy shows similarities to multigrid methods as it uses a hierarchical decomposition of the domain, as well as to the incomplete LU (ILU) factorization and its variants as the sparsity pattern is preserved. Contrary to ILU, the sparsity is maintained using low-rank compression instead of thresholding.

The outline of this paper is as follows. \autoref{sec:IFMM} introduces the IFMM and focuses on the key ideas required to achieve linear complexity; special attention is paid to a graph representation of the method. Numerical examples are subsequently discussed in \autoref{sec:NumericalExamples} to illustrate the computational efficiency of the methodology, both as a direct solver and as a preconditioner in an iterative scheme. The applicability of the method to a Stokes flow problem is finally demonstrated in \autoref{sec:Application}, while summarizing conclusions are drawn in \autoref{sec:Conclusions}. In \autoref{classification}, we included a discussion of theoretical connections between the method discussed in this paper and the algorithms developed by Hackbusch et al.~\cite{bebe08a,borm03a,hack00a,hack02a,hack00b}.

\section{The inverse fast multipole method: key ideas and algorithm}\label{sec:IFMM}
\subsection{Extended sparsification}\label{subsec:extended_sparsification}
It is the aim to obtain a fast solver for the dense linear system $\mathbf{A} \mathbf{x} = \mathbf{b}$, where $\mathbf{A}$ is an $\mH^{2}$-matrix, i.e., a matrix constructed based on strong admissibility criteria (well-separated clusters are assumed to be low-rank only), and with a nested basis for the low-rank representations. 

The first step in obtaining such a solver is to exploit the $\mH^{2}$-structure of the matrix in order to convert the original dense system into an extended sparse system, as solving the latter is more efficient from a computational point of view. The sparsification procedure is based on recursively decomposing a dense matrix into a sparse matrix (representing the ``near field'' interactions) and a low-rank approximation (representing the ``far field'' interactions), combined with the introduction of auxiliary variables. This is done in a multilevel fashion. Similar ideas have been considered for HSS matrices in~\cite{chan07a,ho12a}.

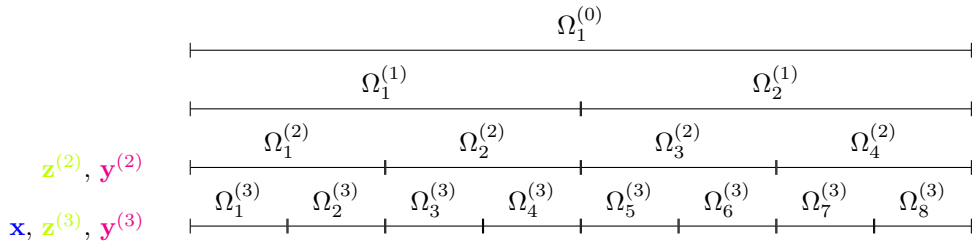
\begin{figure}[hbtp]
	\begin{center}
		\begin{scaletikzpicturetowidth}{\textwidth}
\begin{tikzpicture}[scale=\tikzscale]

		\node (0) at (0,1.75) {};
		\draw [|-|] (0,1.75) -- (10,1.75);
		\node [above] at (5,1.75) (0) {$\Omega_{1}^{(0)}$};
		\node (1) at (0,1.00) {};
		\draw [|-|] (0,1.00) -- (5,1.00);
		\draw [|-|] (5,1.00) -- (10,1.00);
		\node [above] at (2.5,1.00) {$\Omega_{1}^{(1)}$};
		\node [above] at (7.5,1.00) {$\Omega_{2}^{(1)}$};
	
		\node (2) at (0,0.25) {};
		\draw [|-|] (0,0.25) -- (2.5,0.25);
		\draw [|-|] (2.5,0.25) -- (5,0.25);
		\draw [|-|] (5,0.25) -- (7.5,0.25);
		\draw [|-|] (7.5,0.25) -- (10,0.25);
		\node [above] at (1.25,0.25) {$\Omega_{1}^{(2)}$};
		\node [above] at (3.75,0.25) {$\Omega_{2}^{(2)}$};
		\node [above] at (6.25,0.25) {$\Omega_{3}^{(2)}$};
		\node [above] at (8.75,0.25) {$\Omega_{4}^{(2)}$};

		\node[left=3.5mm of 2] () {\textcolor{lime}{$\mathbf{z}^{(2)}$}, \textcolor{magenta}{$\mathbf{y}^{(2)}$}};

		\node (3) at (0,-0.50) {};
		\draw [|-|] (0,-0.50) -- (1.25,-0.50);
		\draw [|-|] (1.25,-0.50) -- (2.50,-0.50);
		\draw [|-|] (2.50,-0.50) -- (3.75,-0.50);
		\draw [|-|] (3.75,-0.50) -- (5,-0.50);
		\draw [|-|] (5,-0.50) -- (6.25,-0.50);
		\draw [|-|] (6.25,-0.50) -- (7.5,-0.50);
		\draw [|-|] (7.5,-0.50) -- (8.75,-0.50);
		\draw [|-|] (8.75,-0.50) -- (10,-0.50);
		\node [above] at (0.625,-0.50) {$\Omega_{1}^{(3)}$};
		\node [above] at (1.875,-0.50) {$\Omega_{2}^{(3)}$};
		\node [above] at (3.125,-0.50) {$\Omega_{3}^{(3)}$};
		\node [above] at (4.375,-0.50) {$\Omega_{4}^{(3)}$};
		\node [above] at (5.625,-0.50) {$\Omega_{5}^{(3)}$};
		\node [above] at (6.875,-0.50) {$\Omega_{6}^{(3)}$};
		\node [above] at (8.125,-0.50) {$\Omega_{7}^{(3)}$};
		\node [above] at (9.375,-0.50) {$\Omega_{8}^{(3)}$};
		
		\node[left=3.5mm of 3] () {\textcolor{blue}{$\mathbf{x}$}, \textcolor{lime}{$\mathbf{z}^{(3)}$}, \textcolor{magenta}{$\mathbf{y}^{(3)}$}};
		
	\end{tikzpicture}
\end{scaletikzpicturetowidth}%
	\end{center}
	\caption{One-dimensional domain $\Omega$ hierarchically subdivided into clusters $\Omega_{i}^{(l)}$.}
	\label{fig:FMM_1D_new}
\end{figure}

\begin{figure}[hbtp]
	\begin{center}
		\begin{scaletikzpicturetowidth}{\textwidth}
\begin{tikzpicture}[scale=\tikzscale]

	\node(0) at (-0.50,-0.50) {(a)};
	\fill[white,draw=black] (0.00,0.00) rectangle (4.80,4.80);
	
	\fill[red,draw=black,very thin] (0.00,4.20) rectangle (0.60,4.80);
	\fill[red,draw=black,very thin] (0.60,4.20) rectangle (1.20,4.80);
	\fill[red,draw=black,very thin] (1.20,4.20) rectangle (1.80,4.80);
	\fill[red,draw=black,very thin] (1.80,4.20) rectangle (2.40,4.80);
	\fill[red,draw=black,very thin] (2.40,4.20) rectangle (3.00,4.80);
	\fill[red,draw=black,very thin] (3.00,4.20) rectangle (3.60,4.80);
	\fill[red,draw=black,very thin] (3.60,4.20) rectangle (4.20,4.80);
	\fill[red,draw=black,very thin] (4.20,4.20) rectangle (4.80,4.80);

	\fill[red,draw=black,very thin] (0.00,3.60) rectangle (0.60,4.20);
	\fill[red,draw=black,very thin] (0.60,3.60) rectangle (1.20,4.20);
	\fill[red,draw=black,very thin] (1.20,3.60) rectangle (1.80,4.20);
	\fill[red,draw=black,very thin] (1.80,3.60) rectangle (2.40,4.20);
	\fill[red,draw=black,very thin] (2.40,3.60) rectangle (3.00,4.20);
	\fill[red,draw=black,very thin] (3.00,3.60) rectangle (3.60,4.20);
	\fill[red,draw=black,very thin] (3.60,3.60) rectangle (4.20,4.20);
	\fill[red,draw=black,very thin] (4.20,3.60) rectangle (4.80,4.20);

	\fill[red,draw=black,very thin] (0.00,3.00) rectangle (0.60,3.60);
	\fill[red,draw=black,very thin] (0.60,3.00) rectangle (1.20,3.60);
	\fill[red,draw=black,very thin] (1.20,3.00) rectangle (1.80,3.60);
	\fill[red,draw=black,very thin] (1.80,3.00) rectangle (2.40,3.60);
	\fill[red,draw=black,very thin] (2.40,3.00) rectangle (3.00,3.60);
	\fill[red,draw=black,very thin] (3.00,3.00) rectangle (3.60,3.60);
	\fill[red,draw=black,very thin] (3.60,3.00) rectangle (4.20,3.60);
	\fill[red,draw=black,very thin] (4.20,3.00) rectangle (4.80,3.60);
	
	\fill[red,draw=black,very thin] (0.00,2.40) rectangle (0.60,3.00);
	\fill[red,draw=black,very thin] (0.60,2.40) rectangle (1.20,3.00);
	\fill[red,draw=black,very thin] (1.20,2.40) rectangle (1.80,3.00);
	\fill[red,draw=black,very thin] (1.80,2.40) rectangle (2.40,3.00);
	\fill[red,draw=black,very thin] (2.40,2.40) rectangle (3.00,3.00);
	\fill[red,draw=black,very thin] (3.00,2.40) rectangle (3.60,3.00);
	\fill[red,draw=black,very thin] (3.60,2.40) rectangle (4.20,3.00);
	\fill[red,draw=black,very thin] (4.20,2.40) rectangle (4.80,3.00);
	
	\fill[red,draw=black,very thin] (0.00,1.80) rectangle (0.60,2.40);
	\fill[red,draw=black,very thin] (0.60,1.80) rectangle (1.20,2.40);
	\fill[red,draw=black,very thin] (1.20,1.80) rectangle (1.80,2.40);
	\fill[red,draw=black,very thin] (1.80,1.80) rectangle (2.40,2.40);
	\fill[red,draw=black,very thin] (2.40,1.80) rectangle (3.00,2.40);
	\fill[red,draw=black,very thin] (3.00,1.80) rectangle (3.60,2.40);
	\fill[red,draw=black,very thin] (3.60,1.80) rectangle (4.20,2.40);
	\fill[red,draw=black,very thin] (4.20,1.80) rectangle (4.80,2.40);

	\fill[red,draw=black,very thin] (0.00,1.20) rectangle (0.60,1.80);
	\fill[red,draw=black,very thin] (0.60,1.20) rectangle (1.20,1.80);
	\fill[red,draw=black,very thin] (1.20,1.20) rectangle (1.80,1.80);
	\fill[red,draw=black,very thin] (1.80,1.20) rectangle (2.40,1.80);
	\fill[red,draw=black,very thin] (2.40,1.20) rectangle (3.00,1.80);
	\fill[red,draw=black,very thin] (3.00,1.20) rectangle (3.60,1.80);
	\fill[red,draw=black,very thin] (3.60,1.20) rectangle (4.20,1.80);
	\fill[red,draw=black,very thin] (4.20,1.20) rectangle (4.80,1.80);
	
	\fill[red,draw=black,very thin] (0.00,0.60) rectangle (0.60,1.20);
	\fill[red,draw=black,very thin] (0.60,0.60) rectangle (1.20,1.20);
	\fill[red,draw=black,very thin] (1.20,0.60) rectangle (1.80,1.20);
	\fill[red,draw=black,very thin] (1.80,0.60) rectangle (2.40,1.20);
	\fill[red,draw=black,very thin] (2.40,0.60) rectangle (3.00,1.20);
	\fill[red,draw=black,very thin] (3.00,0.60) rectangle (3.60,1.20);
	\fill[red,draw=black,very thin] (3.60,0.60) rectangle (4.20,1.20);
	\fill[red,draw=black,very thin] (4.20,0.60) rectangle (4.80,1.20);
	
	\fill[red,draw=black,very thin] (0.00,0.00) rectangle (0.60,0.60);
	\fill[red,draw=black,very thin] (0.60,0.00) rectangle (1.20,0.60);
	\fill[red,draw=black,very thin] (1.20,0.00) rectangle (1.80,0.60);
	\fill[red,draw=black,very thin] (1.80,0.00) rectangle (2.40,0.60);
	\fill[red,draw=black,very thin] (2.40,0.00) rectangle (3.00,0.60);
	\fill[red,draw=black,very thin] (3.00,0.00) rectangle (3.60,0.60);
	\fill[red,draw=black,very thin] (3.60,0.00) rectangle (4.20,0.60);
	\fill[red,draw=black,very thin] (4.20,0.00) rectangle (4.80,0.60);

	\node (1) at (5.30,2.40) {$\simeq$};

	\fill[white,draw=black] (0.00+5.80,0.00) rectangle (4.80+5.80,4.80);

	\fill[red,draw=black,very thin] (0.00+5.80,4.20) rectangle (0.60+5.80,4.80);
	\fill[red,draw=black,very thin] (0.60+5.80,4.20) rectangle (1.20+5.80,4.80);

	\fill[green,draw=black,very thin] (1.20+5.80,4.20) rectangle (1.80+5.80,4.80);
	\fill[green,draw=black,very thin] (1.80+5.80,4.20) rectangle (2.40+5.80,4.80);
	\fill[green,draw=black,very thin] (2.40+5.80,3.60) rectangle (3.60+5.80,4.80);
	\fill[green,draw=black,very thin] (3.60+5.80,3.60) rectangle (4.80+5.80,4.80);
	
	\fill[red,draw=black,very thin] (0.60+5.80,3.60) rectangle (1.20+5.80,4.20);
	\fill[red,draw=black,very thin] (1.20+5.80,3.60) rectangle (1.80+5.80,4.20);
	\fill[red,draw=black,very thin] (0.00+5.80,3.60) rectangle (0.60+5.80,4.20);
	
	\fill[green,draw=black,very thin] (1.80+5.80,3.60) rectangle (2.40+5.80,4.20);
		
	\fill[red,draw=black,very thin] (1.20+5.80,3.00) rectangle (1.80+5.80,3.60);
	\fill[red,draw=black,very thin] (1.80+5.80,3.00) rectangle (2.40+5.80,3.60);
	\fill[red,draw=black,very thin] (0.60+5.80,3.00) rectangle (1.20+5.80,3.60);

	\fill[green,draw=black,very thin] (0.00+5.80,3.00) rectangle (0.60+5.80,3.60);
	\fill[green,draw=black,very thin] (2.40+5.80,3.00) rectangle (3.00+5.80,3.60);	
	\fill[green,draw=black,very thin] (3.00+5.80,3.00) rectangle (3.60+5.80,3.60);

	\fill[red,draw=black,very thin] (1.80+5.80,2.40) rectangle (2.40+5.80,3.00);
	\fill[red,draw=black,very thin] (2.40+5.80,2.40) rectangle (3.00+5.80,3.00);
	\fill[red,draw=black,very thin] (1.20+5.80,2.40) rectangle (1.80+5.80,3.00);
	
	\fill[green,draw=black,very thin] (0.00+5.80,2.40) rectangle (0.60+5.80,3.00);
	\fill[green,draw=black,very thin] (0.60+5.80,2.40) rectangle (1.20+5.80,3.00);
	\fill[green,draw=black,very thin] (3.00+5.80,2.40) rectangle (3.60+5.80,3.00);
	\fill[green,draw=black,very thin] (3.60+5.80,2.40) rectangle (4.80+5.80,3.60);
	
	\fill[red,draw=black,very thin] (2.40+5.80,1.80) rectangle (3.00+5.80,2.40);
	\fill[red,draw=black,very thin] (3.00+5.80,1.80) rectangle (3.60+5.80,2.40);
	\fill[red,draw=black,very thin] (1.80+5.80,1.80) rectangle (2.40+5.80,2.40);
	\fill[green,draw=black,very thin] (1.20+5.80,1.80) rectangle (1.80+5.80,2.40);
	\fill[green,draw=black,very thin] (3.60+5.80,1.80) rectangle (4.20+5.80,2.40);
	\fill[green,draw=black,very thin] (4.20+5.80,1.80) rectangle (4.80+5.80,2.40);

	\fill[green,draw=black,very thin] (0.00+5.80,1.20) rectangle (1.20+5.80,2.40);
	\fill[green,draw=black,very thin] (0.00+5.80,0.00) rectangle (1.20+5.80,1.20);
	\fill[green,draw=black,very thin] (1.20+5.80,0.00) rectangle (2.40+5.80,1.20);

	\fill[red,draw=black,very thin] (3.00+5.80,1.20) rectangle (3.60+5.80,1.80);
	\fill[red,draw=black,very thin] (3.60+5.80,1.20) rectangle (4.20+5.80,1.80);
	\fill[red,draw=black,very thin] (2.40+5.80,1.20) rectangle (3.00+5.80,1.80);

	\fill[green,draw=black,very thin] (1.20+5.80,1.20) rectangle (1.80+5.80,1.80);
	\fill[green,draw=black,very thin] (1.80+5.80,1.20) rectangle (2.40+5.80,1.80);
	\fill[green,draw=black,very thin] (4.20+5.80,1.20) rectangle (4.80+5.80,1.80);
	
	\fill[red,draw=black,very thin] (3.60+5.80,0.60) rectangle (4.20+5.80,1.20);
	\fill[red,draw=black,very thin] (4.20+5.80,0.60) rectangle (4.80+5.80,1.20);
	\fill[red,draw=black,very thin] (3.00+5.80,0.60) rectangle (3.60+5.80,1.20);
	
	\fill[green,draw=black,very thin] (2.40+5.80,0.60) rectangle (3.00+5.80,1.20);
	
	\fill[red,draw=black,very thin] (4.20+5.80,0.00) rectangle (4.80+5.80,0.60);
	\fill[red,draw=black,very thin] (3.60+5.80,0.00) rectangle (4.20+5.80,0.60);
	
	\fill[green,draw=black,very thin] (2.40+5.80,0.00) rectangle (3.00+5.80,0.60);
	\fill[green,draw=black,very thin] (3.00+5.80,0.00) rectangle (3.60+5.80,0.60);

	\node (2) at (12.50,-0.50) {(b)};

	\fill[green,draw=black] (13.00,1.80) rectangle (14.20,3.00);
		
	\node (3) at (14.50,2.40) {$=$};
	
	\fill[cyan,draw=black,very thin] (14.80,1.80) rectangle (15.00,3.00);
	\fill[orange,draw=black,very thin] (15.10,2.80) rectangle (15.30,3.00);
	\fill[cyan,draw=black,very thin] (15.40,2.80) rectangle (16.60,3.00);

	\end{tikzpicture}
\end{scaletikzpicturetowidth}
	\end{center}
	\caption{(a)~Approximation of the dense matrix by an $\mH^{2}$-matrix. (b)~Low rank approximation of an interaction between well-separated clusters: $\mathbf{A}^{(l)}_{ij} \simeq \mathbf{U}^{(l)}_{i} \mathbf{K}^{(l)}_{ij} \mathbf{V}^{(l)^{\mathrm{T}}}_{j}$.}
	\label{fig:approximation_A_tikz}
\end{figure}
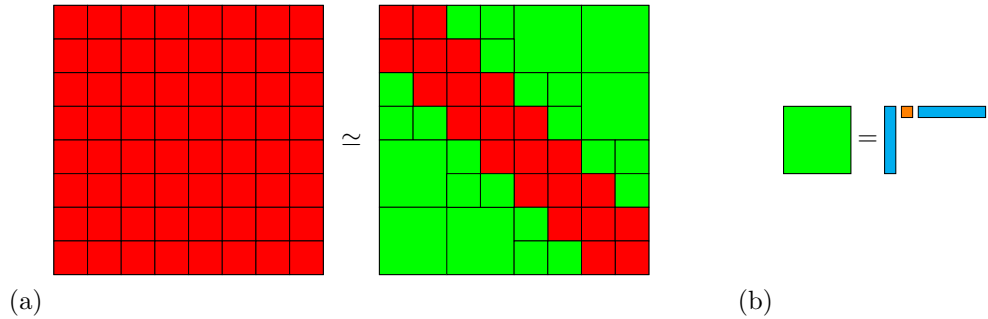

The sparsification concept is illustrated in the following paragraphs for a one-dimensional (1D) example, shown in \autoref{fig:FMM_1D_new}. A 1D domain~$\Omega$ is hierarchically subdivided into clusters $\Omega_{i}^{(l)}$ (with super- and subscripts~$^{(l)}$ and~$_{i}$ referring to the hierarchical level and the cluster index at that level, respectively) until the leaf level $l=L$ is reached (with $L=3$ in \autoref{fig:FMM_1D_new}). Blocks in the matrix $\mathbf{A}$ corresponding to interactions between well-separated clusters $\Omega_{i}^{(l)}$ and $\Omega_{j}^{(l)}$ are replaced by low-rank approximations of rank $r$, i.e., $\mathbf{A}^{(l)}_{ij} \simeq \mathbf{U}^{(l)}_{i} \mathbf{K}^{(l)}_{ij} \mathbf{V}^{(l)^{\mathrm{T}}}_{j}$ (see \autoref{fig:approximation_A_tikz}(b)). For the example of \autoref{fig:FMM_1D_new}, the resulting matrix structure is shown in \autoref{fig:approximation_A_tikz}(a). Note that the nested structure of the $\mH^{2}$-matrix is not shown in this figure. For readers familiar with the FMM, the matrices $\mathbf{U}^{(l)}_{i}$, $\mathbf{K}^{(l)}_{ij}$, and $\mathbf{V}^{(l)^{\mathrm{T}}}_{j}$ can be identified as \texttt{L2P}- (or \texttt{L2L}-), \texttt{M2L}-, and \texttt{P2M}- (or \texttt{M2M}-) operators, respectively.


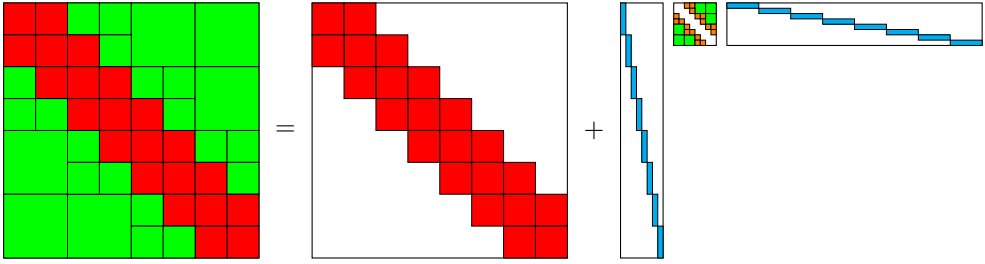
\begin{figure}[hbtp]
	\begin{center}
		\begin{scaletikzpicturetowidth}{\textwidth}
\begin{tikzpicture}[scale=\tikzscale]
	\fill[white,draw=black] (0.00,0.00) rectangle (4.80,4.80);

	\fill[red,draw=black,very thin] (0.00,4.20) rectangle (0.60,4.80);
	\fill[red,draw=black,very thin] (0.60,4.20) rectangle (1.20,4.80);

	\fill[green,draw=black,very thin] (1.20,4.20) rectangle (1.80,4.80);
	\fill[green,draw=black,very thin] (1.80,4.20) rectangle (2.40,4.80);
	\fill[green,draw=black,very thin] (2.40,3.60) rectangle (3.60,4.80);
	\fill[green,draw=black,very thin] (3.60,3.60) rectangle (4.80,4.80);
	
	\fill[red,draw=black,very thin] (0.60,3.60) rectangle (1.20,4.20);
	\fill[red,draw=black,very thin] (1.20,3.60) rectangle (1.80,4.20);
	\fill[red,draw=black,very thin] (0.00,3.60) rectangle (0.60,4.20);
	
	\fill[green,draw=black,very thin] (1.80,3.60) rectangle (2.40,4.20);
		
	\fill[red,draw=black,very thin] (1.20,3.00) rectangle (1.80,3.60);
	\fill[red,draw=black,very thin] (1.80,3.00) rectangle (2.40,3.60);
	\fill[red,draw=black,very thin] (0.60,3.00) rectangle (1.20,3.60);

	\fill[green,draw=black,very thin] (0.00,3.00) rectangle (0.60,3.60);
	\fill[green,draw=black,very thin] (2.40,3.00) rectangle (3.00,3.60);	
	\fill[green,draw=black,very thin] (3.00,3.00) rectangle (3.60,3.60);

	\fill[red,draw=black,very thin] (1.80,2.40) rectangle (2.40,3.00);
	\fill[red,draw=black,very thin] (2.40,2.40) rectangle (3.00,3.00);
	\fill[red,draw=black,very thin] (1.20,2.40) rectangle (1.80,3.00);
	
	\fill[green,draw=black,very thin] (0.00,2.40) rectangle (0.60,3.00);
	\fill[green,draw=black,very thin] (0.60,2.40) rectangle (1.20,3.00);
	\fill[green,draw=black,very thin] (3.00,2.40) rectangle (3.60,3.00);
	\fill[green,draw=black,very thin] (3.60,2.40) rectangle (4.80,3.60);
	
	\fill[red,draw=black,very thin] (2.40,1.80) rectangle (3.00,2.40);
	\fill[red,draw=black,very thin] (3.00,1.80) rectangle (3.60,2.40);
	\fill[red,draw=black,very thin] (1.80,1.80) rectangle (2.40,2.40);
	\fill[green,draw=black,very thin] (1.20,1.80) rectangle (1.80,2.40);
	\fill[green,draw=black,very thin] (3.60,1.80) rectangle (4.20,2.40);
	\fill[green,draw=black,very thin] (4.20,1.80) rectangle (4.80,2.40);

	\fill[green,draw=black,very thin] (0.00,1.20) rectangle (1.20,2.40);
	\fill[green,draw=black,very thin] (0.00,0.00) rectangle (1.20,1.20);
	\fill[green,draw=black,very thin] (1.20,0.00) rectangle (2.40,1.20);

	\fill[red,draw=black,very thin] (3.00,1.20) rectangle (3.60,1.80);
	\fill[red,draw=black,very thin] (3.60,1.20) rectangle (4.20,1.80);
	\fill[red,draw=black,very thin] (2.40,1.20) rectangle (3.00,1.80);

	\fill[green,draw=black,very thin] (1.20,1.20) rectangle (1.80,1.80);
	\fill[green,draw=black,very thin] (1.80,1.20) rectangle (2.40,1.80);
	\fill[green,draw=black,very thin] (4.20,1.20) rectangle (4.80,1.80);
	
	\fill[red,draw=black,very thin] (3.60,0.60) rectangle (4.20,1.20);
	\fill[red,draw=black,very thin] (4.20,0.60) rectangle (4.80,1.20);
	\fill[red,draw=black,very thin] (3.00,0.60) rectangle (3.60,1.20);
	
	\fill[green,draw=black,very thin] (2.40,0.60) rectangle (3.00,1.20);
	
	\fill[red,draw=black,very thin] (4.20,0.00) rectangle (4.80,0.60);
	\fill[red,draw=black,very thin] (3.60,0.00) rectangle (4.20,0.60);
	
	\fill[green,draw=black,very thin] (2.40,0.00) rectangle (3.00,0.60);
	\fill[green,draw=black,very thin] (3.00,0.00) rectangle (3.60,0.60);

	\node (0) at (5.30,2.40) {$=$};
	
	\fill[white,draw=black] (5.80,0.00) rectangle (10.60,4.80);
	\fill[red,draw=black,very thin] (0.00+5.80,4.20) rectangle (0.60+5.80,4.80);
	\fill[red,draw=black,very thin] (0.60+5.80,4.20) rectangle (1.20+5.80,4.80);	
	\fill[red,draw=black,very thin] (0.60+5.80,3.60) rectangle (1.20+5.80,4.20);
	\fill[red,draw=black,very thin] (1.20+5.80,3.60) rectangle (1.80+5.80,4.20);
	\fill[red,draw=black,very thin] (0.00+5.80,3.60) rectangle (0.60+5.80,4.20);
	\fill[red,draw=black,very thin] (1.20+5.80,3.00) rectangle (1.80+5.80,3.60);
	\fill[red,draw=black,very thin] (1.80+5.80,3.00) rectangle (2.40+5.80,3.60);
	\fill[red,draw=black,very thin] (0.60+5.80,3.00) rectangle (1.20+5.80,3.60);	
	\fill[red,draw=black,very thin] (1.80+5.80,2.40) rectangle (2.40+5.80,3.00);
	\fill[red,draw=black,very thin] (2.40+5.80,2.40) rectangle (3.00+5.80,3.00);
	\fill[red,draw=black,very thin] (1.20+5.80,2.40) rectangle (1.80+5.80,3.00);
	\fill[red,draw=black,very thin] (2.40+5.80,1.80) rectangle (3.00+5.80,2.40);
	\fill[red,draw=black,very thin] (3.00+5.80,1.80) rectangle (3.60+5.80,2.40);
	\fill[red,draw=black,very thin] (1.80+5.80,1.80) rectangle (2.40+5.80,2.40);
	\fill[red,draw=black,very thin] (3.00+5.80,1.20) rectangle (3.60+5.80,1.80);
	\fill[red,draw=black,very thin] (3.60+5.80,1.20) rectangle (4.20+5.80,1.80);
	\fill[red,draw=black,very thin] (2.40+5.80,1.20) rectangle (3.00+5.80,1.80);
	\fill[red,draw=black,very thin] (3.60+5.80,0.60) rectangle (4.20+5.80,1.20);
	\fill[red,draw=black,very thin] (4.20+5.80,0.60) rectangle (4.80+5.80,1.20);
	\fill[red,draw=black,very thin] (3.00+5.80,0.60) rectangle (3.60+5.80,1.20);	
	\fill[red,draw=black,very thin] (4.20+5.80,0.00) rectangle (4.80+5.80,0.60);
	\fill[red,draw=black,very thin] (3.60+5.80,0.00) rectangle (4.20+5.80,0.60);

	\node (1) at (11.10,2.40) {$+$};
	
	\fill[white,draw=black] (11.60,0.00) rectangle (12.40,4.80);
	\fill[cyan,draw=black,very thin] (0.00+11.60,4.20) rectangle (0.10+11.60,4.80);
	\fill[cyan,draw=black,very thin] (0.10+11.60,3.60) rectangle (0.20+11.60,4.20);	
	\fill[cyan,draw=black,very thin] (0.20+11.60,3.00) rectangle (0.30+11.60,3.60);	
	\fill[cyan,draw=black,very thin] (0.30+11.60,2.40) rectangle (0.40+11.60,3.00);	
	\fill[cyan,draw=black,very thin] (0.40+11.60,1.80) rectangle (0.50+11.60,2.40);	
	\fill[cyan,draw=black,very thin] (0.50+11.60,1.20) rectangle (0.60+11.60,1.80);	
	\fill[cyan,draw=black,very thin] (0.60+11.60,0.60) rectangle (0.70+11.60,1.20);	
	\fill[cyan,draw=black,very thin] (0.70+11.60,0.00) rectangle (0.80+11.60,0.60);	

	\fill[white,draw=black] (12.60,4.00) rectangle (13.40,4.80);
	
	\fill[orange,draw=black,very thin] (0.20+12.60,4.70) rectangle (0.30+12.60,4.80);
	\fill[orange,draw=black,very thin] (0.30+12.60,4.70) rectangle (0.40+12.60,4.80);
	\fill[green,draw=black,very thin] (0.40+12.60,4.60) rectangle (0.60+12.60,4.80);
	\fill[green,draw=black,very thin] (0.60+12.60,4.60) rectangle (0.80+12.60,4.80);
	\fill[orange,draw=black,very thin] (0.30+12.60,4.60) rectangle (0.40+12.60,4.70);
	
	\fill[orange,draw=black,very thin] (0.00+12.60,4.50) rectangle (0.10+12.60,4.60);
	\fill[orange,draw=black,very thin] (0.40+12.60,4.50) rectangle (0.50+12.60,4.60);
	\fill[orange,draw=black,very thin] (0.50+12.60,4.50) rectangle (0.60+12.60,4.60);
	\fill[orange,draw=black,very thin] (0.00+12.60,4.40) rectangle (0.10+12.60,4.50);
	\fill[orange,draw=black,very thin] (0.10+12.60,4.40) rectangle (0.20+12.60,4.50);
	\fill[orange,draw=black,very thin] (0.50+12.60,4.40) rectangle (0.60+12.60,4.50);
	\fill[green,draw=black,very thin] (0.60+12.60,4.40) rectangle (0.80+12.60,4.60);
	
	\fill[green,draw=black,very thin] (0.00+12.60,4.20) rectangle (0.20+12.60,4.40);
	\fill[green,draw=black,very thin] (0.00+12.60,4.00) rectangle (0.20+12.60,4.20);
	\fill[green,draw=black,very thin] (0.20+12.60,4.00) rectangle (0.40+12.60,4.20);

	\fill[orange,draw=black,very thin] (0.20+12.60,4.30) rectangle (0.30+12.60,4.40);
	\fill[orange,draw=black,very thin] (0.20+12.60,4.20) rectangle (0.30+12.60,4.30);
	\fill[orange,draw=black,very thin] (0.30+12.60,4.20) rectangle (0.40+12.60,4.30);
	\fill[orange,draw=black,very thin] (0.40+12.60,4.10) rectangle (0.50+12.60,4.20);
	\fill[orange,draw=black,very thin] (0.40+12.60,4.00) rectangle (0.50+12.60,4.10);
	\fill[orange,draw=black,very thin] (0.50+12.60,4.00) rectangle (0.60+12.60,4.10);
	\fill[orange,draw=black,very thin] (0.60+12.60,4.30) rectangle (0.70+12.60,4.40);
	\fill[orange,draw=black,very thin] (0.70+12.60,4.30) rectangle (0.80+12.60,4.40);
	\fill[orange,draw=black,very thin] (0.70+12.60,4.20) rectangle (0.80+12.60,4.30);

	\fill[white,draw=black] (13.60,4.00) rectangle (18.40,4.80);
	\fill[cyan,draw=black,very thin] (0.00+13.60,4.70) rectangle (0.60+13.60,4.80);
	\fill[cyan,draw=black,very thin] (0.60+13.60,4.60) rectangle (1.20+13.60,4.70);
	\fill[cyan,draw=black,very thin] (1.20+13.60,4.50) rectangle (1.80+13.60,4.60);
	\fill[cyan,draw=black,very thin] (1.80+13.60,4.40) rectangle (2.40+13.60,4.50);
	\fill[cyan,draw=black,very thin] (2.40+13.60,4.30) rectangle (3.00+13.60,4.40);
	\fill[cyan,draw=black,very thin] (3.00+13.60,4.20) rectangle (3.60+13.60,4.30);
	\fill[cyan,draw=black,very thin] (3.60+13.60,4.10) rectangle (4.20+13.60,4.20);
	\fill[cyan,draw=black,very thin] (4.20+13.60,4.00) rectangle (4.80+13.60,4.10);

	\end{tikzpicture}
\end{scaletikzpicturetowidth}
	\end{center}
	\caption{Decomposition of the $\mH^{2}$-matrix $\mathbf{A}$ at the leaf level $L$ (with $L=3$): $\mathbf{A} = \mathbf{S}^{(3)} + \mathbf{U}^{(3)} \mathbf{K}^{(3)}\mathbf{V}^{(3)^\mathrm{T}}$.}
	\label{fig:decomp_l3_tikz}
\end{figure}

\begin{figure}[hbtp]
	\centering\includegraphics[width=\textwidth]{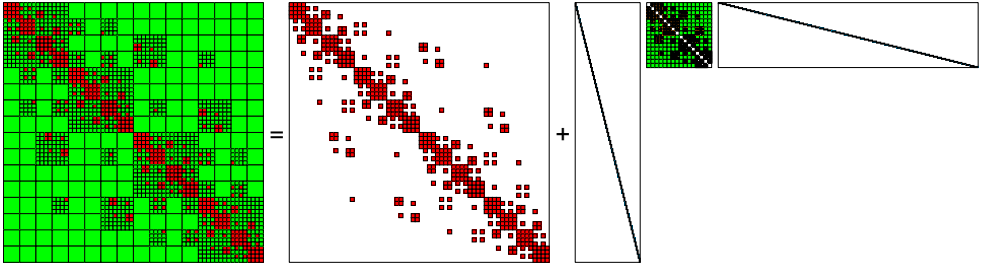}
	\caption{Quadtree decomposition of the matrix. This figure is showing the decomposition at level 3 for a 2D test case. This shows how the sketch in \autoref{fig:decomp_l3_tikz} is generalized in dimensions higher than 1.}
	\label{fig:decomp_2D_l3}
\end{figure}

The dense matrix~$\mathbf{A}$ can consequently be decomposed at the leaf level as follows:
\begin{align}
	\mathbf{A} = \mathbf{S}^{(3)} + \mathbf{U}^{(3)} \mathbf{K}^{(3)}\mathbf{V}^{(3)^\mathrm{T}}
	\label{eq:A_decomp_l3}
 \end{align}
The matrix $\mathbf{S}^{(3)}$ in \autoref{eq:A_decomp_l3} is a sparse matrix containing the interactions of leaf clusters with their neighbors (including the self-interactions). The full rank matrix blocks can be denoted as \texttt{P2P}-operators if the FMM terminology is employed. In 1D, $\mathbf{S}^{(3)}$ is a tri-diagonal matrix. However, for general point distributions, this matrix is a general sparse matrix, with no ``particular'' structure. The matrix $\mathbf{U}^{(3)} \mathbf{K}^{(3)}\mathbf{V}^{(3)^\mathrm{T}}$ represents the dense fill-ins arising from far field interactions. The decomposition presented in \autoref{eq:A_decomp_l3} is pictorially depicted in \autoref{fig:decomp_l3_tikz} for a 1D example and in \autoref{fig:decomp_2D_l3} for a 2D example.

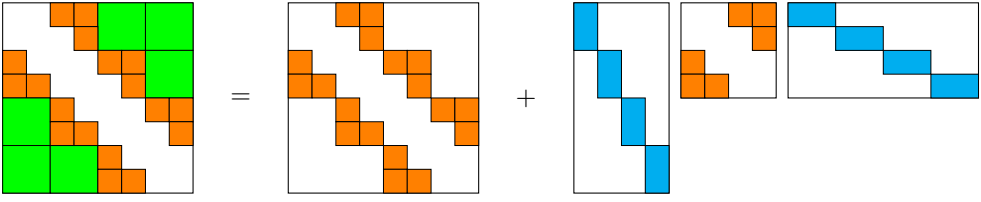
\begin{figure}[hbtp]
	\begin{center}
		\begin{scaletikzpicturetowidth}{\textwidth}
\begin{tikzpicture}[scale=\tikzscale]

	\fill[white,draw=black] (0.00,0.00) rectangle (0.80,0.80);
		
	\fill[orange,draw=black,very thin] (0.20,0.70) rectangle (0.30,0.80);
	\fill[orange,draw=black,very thin] (0.30,0.70) rectangle (0.40,0.80);
	\fill[green,draw=black,very thin] (0.40,0.60) rectangle (0.60,0.80);
	\fill[green,draw=black,very thin] (0.60,0.60) rectangle (0.80,0.80);
	\fill[orange,draw=black,very thin] (0.30,0.60) rectangle (0.40,0.70);
		
	\fill[orange,draw=black,very thin] (0.00,0.50) rectangle (0.10,0.60);
	\fill[orange,draw=black,very thin] (0.40,0.50) rectangle (0.50,0.60);
	\fill[orange,draw=black,very thin] (0.50,0.50) rectangle (0.60,0.60);
	\fill[orange,draw=black,very thin] (0.00,0.40) rectangle (0.10,0.50);
	\fill[orange,draw=black,very thin] (0.10,0.40) rectangle (0.20,0.50);
	\fill[orange,draw=black,very thin] (0.50,0.40) rectangle (0.60,0.50);
	\fill[green,draw=black,very thin] (0.60,0.40) rectangle (0.80,0.60);
		
	\fill[green,draw=black,very thin] (0.00,0.20) rectangle (0.20,0.40);
	\fill[green,draw=black,very thin] (0.00,0.00) rectangle (0.20,0.20);
	\fill[green,draw=black,very thin] (0.20,0.00) rectangle (0.40,0.20);
	
	\fill[orange,draw=black,very thin] (0.20,0.30) rectangle (0.30,0.40);
	\fill[orange,draw=black,very thin] (0.20,0.20) rectangle (0.30,0.30);
	\fill[orange,draw=black,very thin] (0.30,0.20) rectangle (0.40,0.30);
	\fill[orange,draw=black,very thin] (0.40,0.10) rectangle (0.50,0.20);
	\fill[orange,draw=black,very thin] (0.40,0.00) rectangle (0.50,0.10);
	\fill[orange,draw=black,very thin] (0.50,0.00) rectangle (0.60,0.10);
	\fill[orange,draw=black,very thin] (0.60,0.30) rectangle (0.70,0.40);
	\fill[orange,draw=black,very thin] (0.70,0.30) rectangle (0.80,0.40);
	\fill[orange,draw=black,very thin] (0.70,0.20) rectangle (0.80,0.30);

	\node (0) at (1.00,0.40) {$=$};

	\fill[white,draw=black] (1.20,0.00) rectangle (2.00,0.80);
	\fill[orange,draw=black,very thin] (0.20+1.20,0.70) rectangle (0.30+1.20,0.80);
	\fill[orange,draw=black,very thin] (0.30+1.20,0.70) rectangle (0.40+1.20,0.80);
	\fill[orange,draw=black,very thin] (0.30+1.20,0.60) rectangle (0.40+1.20,0.70);
	\fill[orange,draw=black,very thin] (0.00+1.20,0.50) rectangle (0.10+1.20,0.60);
	\fill[orange,draw=black,very thin] (0.40+1.20,0.50) rectangle (0.50+1.20,0.60);
	\fill[orange,draw=black,very thin] (0.50+1.20,0.50) rectangle (0.60+1.20,0.60);
	\fill[orange,draw=black,very thin] (0.00+1.20,0.40) rectangle (0.10+1.20,0.50);
	\fill[orange,draw=black,very thin] (0.10+1.20,0.40) rectangle (0.20+1.20,0.50);
	\fill[orange,draw=black,very thin] (0.50+1.20,0.40) rectangle (0.60+1.20,0.50);
	\fill[orange,draw=black,very thin] (0.20+1.20,0.30) rectangle (0.30+1.20,0.40);
	\fill[orange,draw=black,very thin] (0.20+1.20,0.20) rectangle (0.30+1.20,0.30);
	\fill[orange,draw=black,very thin] (0.30+1.20,0.20) rectangle (0.40+1.20,0.30);
	\fill[orange,draw=black,very thin] (0.40+1.20,0.10) rectangle (0.50+1.20,0.20);
	\fill[orange,draw=black,very thin] (0.40+1.20,0.00) rectangle (0.50+1.20,0.10);
	\fill[orange,draw=black,very thin] (0.50+1.20,0.00) rectangle (0.60+1.20,0.10);
	\fill[orange,draw=black,very thin] (0.60+1.20,0.30) rectangle (0.70+1.20,0.40);
	\fill[orange,draw=black,very thin] (0.70+1.20,0.30) rectangle (0.80+1.20,0.40);
	\fill[orange,draw=black,very thin] (0.70+1.20,0.20) rectangle (0.80+1.20,0.30);
	
	\node (1) at (2.20,0.40) {$+$};
	
	\fill[white,draw=black] (2.40,0.00) rectangle (2.80,0.80);
	\fill[cyan,draw=black,very thin] (0.00+2.40,0.60) rectangle (0.10+2.40,0.80);
	\fill[cyan,draw=black,very thin] (0.10+2.40,0.40) rectangle (0.20+2.40,0.60);	
	\fill[cyan,draw=black,very thin] (0.20+2.40,0.20) rectangle (0.30+2.40,0.40);	
	\fill[cyan,draw=black,very thin] (0.30+2.40,0.00) rectangle (0.40+2.40,0.20);

	\fill[white,draw=black] (2.85,0.40) rectangle (3.25,0.80);
	\fill[orange,draw=black,very thin] (0.20+2.85,0.70) rectangle (0.30+2.85,0.80);
	\fill[orange,draw=black,very thin] (0.30+2.85,0.70) rectangle (0.40+2.85,0.80);
	\fill[orange,draw=black,very thin] (0.30+2.85,0.60) rectangle (0.40+2.85,0.70);	
	\fill[orange,draw=black,very thin] (0.00+2.85,0.50) rectangle (0.10+2.85,0.60);
	\fill[orange,draw=black,very thin] (0.00+2.85,0.40) rectangle (0.10+2.85,0.50);
	\fill[orange,draw=black,very thin] (0.10+2.85,0.40) rectangle (0.20+2.85,0.50);

	\fill[white,draw=black] (3.30,0.40) rectangle (4.10,0.80);
	\fill[cyan,draw=black,very thin] (0.00+3.30,0.70) rectangle (0.20+3.30,0.80);
	\fill[cyan,draw=black,very thin] (0.20+3.30,0.60) rectangle (0.40+3.30,0.70);
	\fill[cyan,draw=black,very thin] (0.40+3.30,0.50) rectangle (0.60+3.30,0.60);
	\fill[cyan,draw=black,very thin] (0.60+3.30,0.40) rectangle (0.80+3.30,0.50);

	\end{tikzpicture}
\end{scaletikzpicturetowidth}
	\end{center}
	\caption{Decomposition of the matrix $\mathbf{K}^{(3)}$ at level $l$ (with $l=2$): $\mathbf{K}^{(3)} = \mathbf{S}^{(2)} + \mathbf{U}^{(2)} \mathbf{K}^{(2)}\mathbf{V}^{(2)^\mathrm{T}}$.}
	\label{fig:decomp_l2_tikz}
\end{figure}

\begin{figure}[hbtp]
	\centering\includegraphics[width=\textwidth]{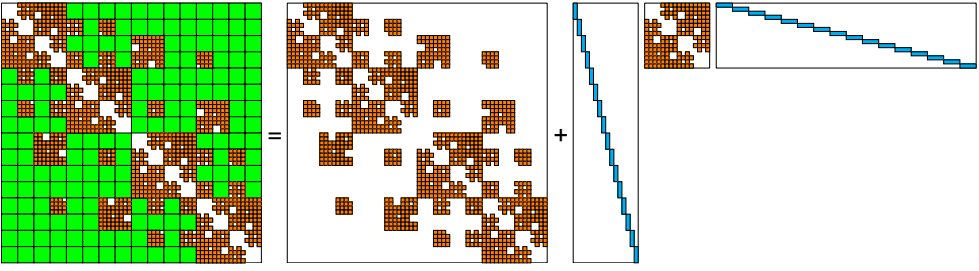}
	\caption{Quadtree decomposition at level 2 for a 2D test case. Compare with \autoref{fig:decomp_l2_tikz} for a 1D distribution of points.}
	\label{fig:decomp_2D_l2}
\end{figure}

The matrix~$\mathbf{K}^{(3)}$ in \autoref{eq:A_decomp_l3} is representative for all far field interactions. It is still rather dense and is therefore further decomposed using a similar decomposition as for $\mathbf{A}$:
\begin{align}
	\mathbf{A} =
	\mathbf{S}^{(3)} + \mathbf{U}^{(3)} \left( \mathbf{S}^{(2)} + \mathbf{U}^{(2)} \mathbf{K}^{(2)}\mathbf{V}^{(2)^\mathrm{T}} \right) \mathbf{V}^{(3)^\mathrm{T}}
	\label{eq:A_decomp_l2}
 \end{align}
This decomposition at level $l=2$ is illustrated in \autoref{fig:decomp_l2_tikz} and \ref{fig:decomp_2D_l2} for a 1D and 2D example, respectively. The sparse matrix $\mathbf{S}^{(2)}$ in \autoref{eq:A_decomp_l2} represents near field interactions at level $l=2$ (that have not been accounted for through $\mathbf{S}^{(3)}$), while $\mathbf{K}^{(2)}$ represents the remaining far field interactions. As no well-separated clusters exist at level $l=1$, the matrix $\mathbf{K}^{(2)}$ can not be further decomposed.

\autoref{eq:A_decomp_l2} is still a dense system, but can now be converted into an extended sparse system by introducing auxiliary variables at each level that contains well-separated clusters (i.e., $l \geq 2$). More specifically, multipole coefficients $\mathbf{y}^{(l)}$ are defined as:
\begin{align}
	\mathbf{y}^{(3)} &= \mathbf{V}^{(3)^\mathrm{T}}\mathbf{x} \label{eq:y3}\\
	\mathbf{y}^{(2)} &= \mathbf{V}^{(2)^\mathrm{T}}\mathbf{y}^{(3)} \label{eq:y2}
 \end{align}
This allows rewriting the equation $\mathbf{A} \mathbf{x} = \mathbf{b}$ as:
\begin{align}
	\mathbf{S}^{(3)}\mathbf{x} + \mathbf{U}^{(3)} \left( \mathbf{S}^{(2)}\mathbf{y}^{(3)} + \mathbf{U}^{(2)} \mathbf{K}^{(2)}\mathbf{y}^{(2)} \right) &= \mathbf{b}
	\label{eq:A_simplified}
 \end{align}
Furthermore, local coefficients $\mathbf{z}^{(l)}$ are identified at each level:
\begin{align}
	\mathbf{z}^{(2)} &= \mathbf{K}^{(2)}\mathbf{y}^{(2)} \label{eq:z2}\\
	\mathbf{z}^{(3)} &=  \mathbf{S}^{(2)}\mathbf{y}^{(3)} + \mathbf{U}^{(2)}\mathbf{z}^{(2)} \label{eq:z3}
 \end{align}
The original equation~$\mathbf{A} \mathbf{x} = \mathbf{b}$ thus becomes:
\begin{align}
	\mathbf{S}^{(3)}\mathbf{x} + \mathbf{U}^{(3)} \mathbf{z}^{(3)} &= \mathbf{b}
	\label{eq:decomposition_compact}
 \end{align}
\autoref{eq:y3}--\eqref{eq:decomposition_compact} provide an extended sparse system of equations: 
\begin{align}
	\begin{bmatrix}	
	\mathbf{S}^{(3)}  & \mathbf{U}^{(3)} & & &  \\
	\mathbf{V}^{(3)^\mathrm{T}}  & & -\mathbf{I} & &  \\
	  & -\mathbf{I} & \mathbf{S}^{(2)} & \mathbf{U}^{(2)} &  \\
	& & \mathbf{V}^{(2)^\mathrm{T}} & &-\mathbf{I}  \\
	& & & -\mathbf{I} & \mathbf{K}^{(2)}
	\end{bmatrix}
	\begin{Bmatrix}	
		\mathbf{x} \\ \mathbf{z}^{(3)} \\ \mathbf{y}^{(3)} \\ \mathbf{z}^{(2)} \\ \mathbf{y}^{(2)}
	\end{Bmatrix}
	&= 
	\begin{Bmatrix}	
		\mathbf{b} \\ \mathbf{0} \\ \mathbf{0} \\ \mathbf{0} \\ \mathbf{0}
	\end{Bmatrix}
	\label{eq:Ax=b_extended}
 \end{align}

\begin{figure}[hbtp]
	\begin{center}
		\input{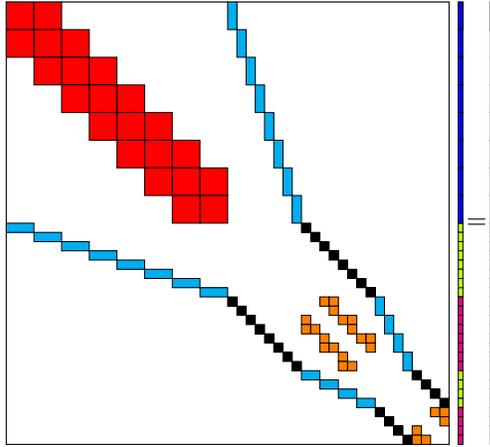}
	\end{center}
	\caption{A pictorial representation of the extended sparse system of \autoref{eq:Ax=b_extended} for a 1D distribution of points. Red blocks are \texttt{P2P}-operators, cyan blocks are \texttt{L2P}-, \texttt{P2M}-, \texttt{L2L}-, and \texttt{M2M}-operators, and orange blocks are \texttt{M2L}-operators. The black blocks represent negative identity matrices. The original unknowns $\mathbf{x}$ are shown in blue, the locals $\mathbf{z}^{(l)}$ in light green, and the multipoles $\mathbf{y}^{(l)}$ in magenta. The right hand side vector $\mathbf{b}$ is depicted in gray.}
	\label{fig:extended_matrix_new_tikz}
\end{figure}

\begin{figure}[hbtp]
	\centering\includegraphics[width=0.5\textwidth]{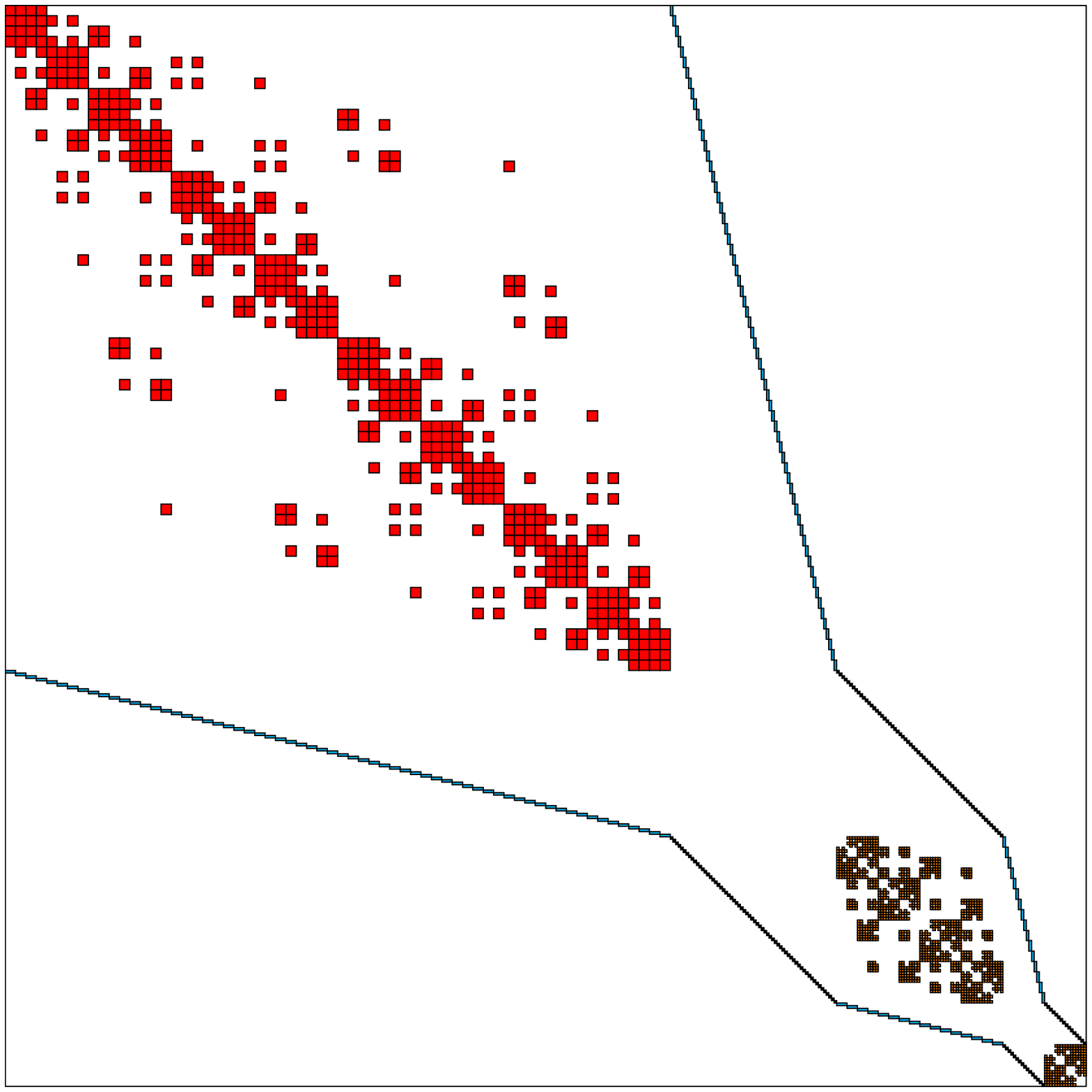}
	\caption{Extended sparse matrix for a 2D distribution of points. Compare with the sketch shown in \autoref{fig:extended_matrix_new_tikz} for a 1D distribution of points.}
	\label{fig:extended_matrix_2D}
\end{figure}

The extended sparse matrix in \autoref{eq:Ax=b_extended} (indicated as $\mathbf{E}$ in the following) contains matrices from the decompositions introduced in \autoref{eq:A_decomp_l3} and~\eqref{eq:A_decomp_l2} and exhibits a symmetric pattern in terms of fill-in, as is shown in \autoref{fig:extended_matrix_new_tikz} and \ref{fig:extended_matrix_2D}. In fact, for a symmetric matrix $\mathbf{A}$, the matrix $\mathbf{E}$ is symmetric as well. This is due to the specific ordering of the equations. The elimination and substitution processes are thus performed on this sparse system rather than on the original dense system. 

This sparsification procedure is a very general and powerful technique that is not restricted to 1D or 2D problems. In general, any dense $\mH^{2}$-matrix can be recursively decomposed into sparse matrices and low-rank approximations to obtain an extended sparsified matrix.

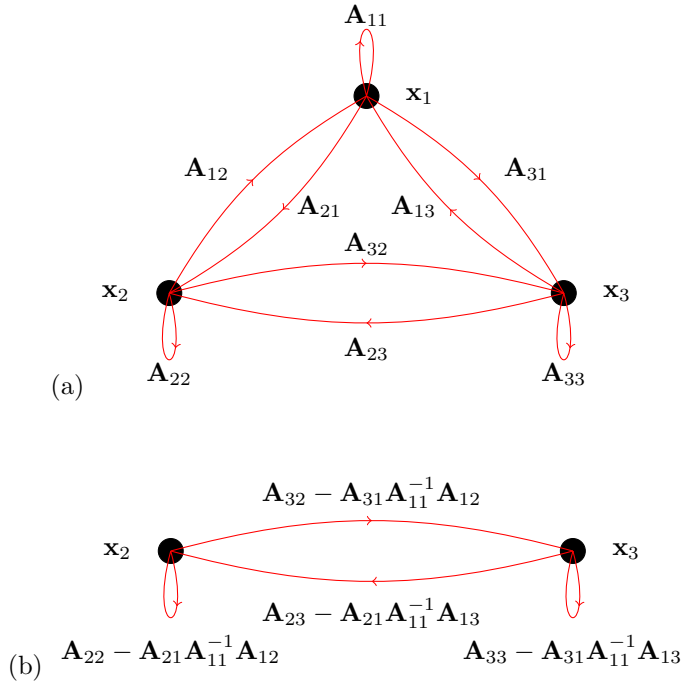
\begin{figure}[hbtp]
	\begin{center}
		(a)~\begin{tikzpicture}[scale=3.5]

		\node (0) at (0,0) {};
		
		\node[draw=black,fill=black,circle,minimum size=2mm,right=15mm of 0] (1) {};			
		\draw [red,->-=.4,loop above,looseness=10,min distance=3.5mm] (1.center) to (1.center);

		\node[draw=black,fill=black,circle,minimum size=2mm,below left=25mm and 25mm of 1.center] (2) {};			
		\draw [red,->-=.4,loop below,looseness=10,min distance=3.5mm] (2.center) to (2.center); 
		
		\node[draw=black,fill=black,circle,minimum size=2mm,below right=25mm and 25mm of 1.center] (3) {};	
		\draw [red,->-=.4,loop below,looseness=10,min distance=3.5mm] (3.center) to (3.center);

		\draw[red,bend left=15,->-=.5]  (2.center) to (1.center);
		\draw[red,bend left=15,->-=.5]  (1.center) to (2.center);

		\draw[red,bend left=15,->-=.5]  (3.center) to (1.center);
		\draw[red,bend left=15,->-=.5]  (1.center) to (3.center);

		\draw[red,bend left=15,->-=.5]  (2.center) to (3.center);
		\draw[red,bend left=15,->-=.5]  (3.center) to (2.center);

		\node[right=4mm of 1.center] () {\textcolor{black}{$\mathbf{x}_{1}$}};
		\node[left=4mm of 2.center] () {\textcolor{black}{$\mathbf{x}_{2}$}};
		\node[right=4mm of 3.center] () {\textcolor{black}{$\mathbf{x}_{3}$}};

		\node[above=8mm of 1.center] () {$\mathbf{A}_{11}$};
		\node[below=8mm of 2.center] () {$\mathbf{A}_{22}$};
		\node[below=8mm of 3.center] () {$\mathbf{A}_{33}$};

		\node[below=17mm of 1.center] () {$\mathbf{A}_{32}$};
		\node[below=31mm of 1.center] () {$\mathbf{A}_{23}$};

		\node[below right=7mm and 17mm of 1.center] () {$\mathbf{A}_{31}$};
		\node[below right=12mm and 2mm of 1.center] () {$\mathbf{A}_{13}$};	

		\node[below left=7mm and 17mm of 1.center] () {$\mathbf{A}_{12}$};
		\node[below left=12mm and 2mm of 1.center] () {$\mathbf{A}_{21}$};	
	

	\end{tikzpicture}\vspace{2.5em}\\
		(b)~\begin{tikzpicture}[scale=3.5]

		\node (0) at (0,0) {};

		\node[right=15mm of 0] (1) {};			
		
		\node[draw=black,fill=black,circle,minimum size=2mm,left=25mm of 1.center] (2) {};			
		\draw [red,->-=.4,loop below,looseness=10,min distance=3.5mm] (2.center) to (2.center); 
		
		\node[draw=black,fill=black,circle,minimum size=2mm,right=25mm of 1.center] (3) {};	
		\draw [red,->-=.4,loop below,looseness=10,min distance=3.5mm] (3.center) to (3.center);



		\draw[red,bend left=15,->-=.5]  (2.center) to (3.center);
		\draw[red,bend left=15,->-=.5]  (3.center) to (2.center);

		\node[left=4mm of 2.center] () {\textcolor{black}{$\mathbf{x}_{2}$}};
		\node[right=4mm of 3.center] () {\textcolor{black}{$\mathbf{x}_{3}$}};

		\node[below=10mm of 2.center] () {$\mathbf{A}_{22} - \mathbf{A}_{21} \mathbf{A}^{-1}_{11} \mathbf{A}_{12}$};
		\node[below=10mm of 3.center] () {$\mathbf{A}_{33} - \mathbf{A}_{31} \mathbf{A}^{-1}_{11} \mathbf{A}_{13}$};

		\node[below=-11mm of 1.center] () {$\mathbf{A}_{32} - \mathbf{A}_{31} \mathbf{A}^{-1}_{11} \mathbf{A}_{12}$};
		\node[below=5mm of 1.center] () {$\mathbf{A}_{23} - \mathbf{A}_{21} \mathbf{A}^{-1}_{11} \mathbf{A}_{13}$};


	

	\end{tikzpicture}
	\end{center}
	\caption{The graph of \autoref{eq:Ax=b} (a)~before and (b) after the elimination of the variable $\mathbf{x}_{1}$.}
	\label{fig:graph_example}
\end{figure}

In order to further illustrate the sparsification procedure, it is instructive to consider graph representations of the linear systems under concern. This will be particularly helpful for the interpretation of the second key step of the algorithm, as will be elucidated in \autoref{subsec:compression}. Before discussing the graph of \autoref{eq:Ax=b_extended}, the idea of representing a linear system as a graph is briefly recalled here by means of a small example. Consider the following generic $3 \times 3$ block system (this is not related to any of the systems discussed earlier):
\begin{align}
	\begin{bmatrix}	
	\mathbf{A}_{11}  & \mathbf{A}_{12} & \mathbf{A}_{13}  \\
	\mathbf{A}_{21}  & \mathbf{A}_{22} & \mathbf{A}_{23}  \\
	\mathbf{A}_{31}  & \mathbf{A}_{32} & \mathbf{A}_{33}
	\end{bmatrix}
	\begin{Bmatrix}	
		\mathbf{x}_{1} \\ \mathbf{x}_{2} \\ \mathbf{x}_{3}
	\end{Bmatrix}
	&= 
	\begin{Bmatrix}	
		\mathbf{b}_{1} \\ \mathbf{b}_{2} \\ \mathbf{b}_{3}
	\end{Bmatrix}
	\label{eq:Ax=b}
 \end{align}
A graph representation of \autoref{eq:Ax=b}, consisting of a set of vertices and edges, is depicted in \autoref{fig:graph_example}(a). Each vertex in the graph is associated with a variable $\mathbf{x}_{i}$, while the set of incoming edges into that vertex defines an equation that corresponds to row~$i$ of \autoref{eq:Ax=b}. The edges correspond to the matrix subblocks $\mathbf{A}_{ij}$ and represent the contribution of variable $\mathbf{x}_{j}$ to the equation of row~$i$. The arrows indicate the direction of each edge. The effect of eliminating a variable from \autoref{eq:Ax=b} on the corresponding graph is illustrated in \autoref{fig:graph_example}(b), highlighting how the edges are updated accordingly.

It is worth noting that the columns of the matrix have a clear meaning (i.e., they correspond to the unknowns), while the meaning of the rows is less clear as they correspond to the equations (and these can easily be permuted). Correspondingly in the graph, the tail of an arrow has a clear interpretation as it is associated with an unknown, while the head is somewhat arbitrary as it corresponds to an equation (and thus a specific ordering of the equations).

\begin{figure}[hbtp]
	\begin{center}
		\begin{scaletikzpicturetowidth}{\textwidth}
\begin{tikzpicture}[scale=\tikzscale]

		\node (0) at (0,0) {};

		\node[draw=red,fill=red,circle,inner sep=0pt,minimum size=0mm,right=00mm of 0] (1) {};	
		\node[draw=red,fill=red,circle,inner sep=0pt,minimum size=0mm,right=30mm of 0] (2) {};	
		\node[draw=red,fill=red,circle,inner sep=0pt,minimum size=0mm,right=60mm of 0] (3) {};	
		\node[draw=red,fill=red,circle,inner sep=0pt,minimum size=0mm,right=90mm of 0] (4) {};

		\node[draw=blue,fill=blue,circle,minimum size=2mm,below right=0mm and -8mm of 1.center] (5) {};	
		\node[draw=blue,fill=blue,circle,minimum size=2mm,below right=0mm and 8mm of 1.center] (6) {};	
		\node[draw=blue,fill=blue,circle,minimum size=2mm,below right=0mm and -8mm of 2.center] (7) {};	
		\node[draw=blue,fill=blue,circle,minimum size=2mm,below right=0mm and 8mm of 2.center] (8) {};	
		\node[draw=blue,fill=blue,circle,minimum size=2mm,below right=0mm and -8mm of 3.center] (9) {};	
		\node[draw=blue,fill=blue,circle,minimum size=2mm,below right=0mm and 8mm of 3.center] (10) {};	
		\node[draw=blue,fill=blue,circle,minimum size=2mm,below right=0mm and -8mm of 4.center] (11) {};	
		\node[draw=blue,fill=blue,circle,minimum size=2mm,below right=0mm and 8mm of 4.center] (12) {};

		\draw [red,->-=.4,loop below,looseness=10,min distance=5mm] (5.center) to (5.center); 
		\draw [red,->-=.4,loop below,looseness=10,min distance=5mm] (6.center) to (6.center); 
		\draw [red,->-=.4,loop below,looseness=10,min distance=5mm] (7.center) to (7.center); 
		\draw [red,->-=.4,loop below,looseness=10,min distance=5mm] (8.center) to (8.center); 
		\draw [red,->-=.4,loop below,looseness=10,min distance=5mm] (9.center) to (9.center); 
		\draw [red,->-=.4,loop below,looseness=10,min distance=5mm] (10.center) to (10.center); 
		\draw [red,->-=.4,loop below,looseness=10,min distance=5mm] (11.center) to (11.center); 
		\draw [red,->-=.4,loop below,looseness=10,min distance=5mm] (12.center) to (12.center);

		\draw[red,->-=.75] (5.center) -- (6.center);
		\draw[red,->-=.75] (6.center) -- (7.center);
		\draw[red,->-=.75] (7.center) -- (8.center);
		\draw[red,->-=.75] (8.center) -- (9.center);
		\draw[red,->-=.75] (9.center) -- (10.center);
		\draw[red,->-=.75] (10.center) -- (11.center);
		\draw[red,->-=.75] (11.center) -- (12.center);

		\draw[red,->-=.75] (6.center) -- (5.center);
		\draw[red,->-=.75] (7.center) -- (6.center);
		\draw[red,->-=.75] (8.center) -- (7.center);
		\draw[red,->-=.75] (9.center) -- (8.center);
		\draw[red,->-=.75] (10.center) -- (9.center);
		\draw[red,->-=.75] (11.center) -- (10.center);
		\draw[red,->-=.75] (12.center) -- (11.center);

		\draw[red,bend left=40,->-=.75] (5.center) to (7.center);
		\draw[red,bend left=40,-<-=.25] (5.center) to (7.center);
		\draw[red,bend left=40,->-=.75] (5.center) to (8.center);
		\draw[red,bend left=40,-<-=.25] (5.center) to (8.center);
		\draw[red,bend left=40,->-=.75] (5.center) to (9.center);
		\draw[red,bend left=40,-<-=.25] (5.center) to (9.center);
		\draw[red,bend left=40,->-=.75] (5.center) to (10.center);
		\draw[red,bend left=40,-<-=.25] (5.center) to (10.center);
		\draw[red,bend right=40,->-=.75] (5.center) to (11.center);
		\draw[red,bend right=40,-<-=.25] (5.center) to (11.center);
		\draw[red,bend right=40,->-=.75] (5.center) to (12.center);
		\draw[red,bend right=40,-<-=.25] (5.center) to (12.center);

		\draw[red,bend left=40,->-=.75] (6.center) to (8.center);
		\draw[red,bend left=40,-<-=.25] (6.center) to (8.center);
		\draw[red,bend left=40,->-=.75] (6.center) to (9.center);
		\draw[red,bend left=40,-<-=.25] (6.center) to (9.center);
		\draw[red,bend left=40,->-=.75] (6.center) to (10.center);
		\draw[red,bend left=40,-<-=.25] (6.center) to (10.center);
		\draw[red,bend right=40,->-=.75] (6.center) to (11.center);
		\draw[red,bend right=40,-<-=.25] (6.center) to (11.center);
		\draw[red,bend right=40,->-=.75] (6.center) to (12.center);
		\draw[red,bend right=40,-<-=.25] (6.center) to (12.center);
		
		\draw[red,bend right=40,->-=.75] (7.center) to (9.center);
		\draw[red,bend right=40,-<-=.25] (7.center) to (9.center);
		\draw[red,bend right=40,->-=.75] (7.center) to (10.center);
		\draw[red,bend right=40,-<-=.25] (7.center) to (10.center);
		\draw[red,bend left=40,->-=.75] (7.center) to (11.center);
		\draw[red,bend left=40,-<-=.25] (7.center) to (11.center);
		\draw[red,bend left=40,->-=.75] (7.center) to (12.center);
		\draw[red,bend left=40,-<-=.25] (7.center) to (12.center);

		\draw[red,bend right=40,->-=.75] (8.center) to (10.center);
		\draw[red,bend right=40,-<-=.25] (8.center) to (10.center);
		\draw[red,bend left=40,->-=.75] (8.center) to (11.center);
		\draw[red,bend left=40,-<-=.25] (8.center) to (11.center);
		\draw[red,bend left=40,->-=.75] (8.center) to (12.center);
		\draw[red,bend left=40,-<-=.25] (8.center) to (12.center);
		
		\draw[red,bend left=40,->-=.75] (9.center) to (11.center);
		\draw[red,bend left=40,-<-=.25] (9.center) to (11.center);
		\draw[red,bend left=40,->-=.75] (9.center) to (12.center);
		\draw[red,bend left=40,-<-=.25] (9.center) to (12.center);

		\draw[red,bend left=40,->-=.75] (10.center) to (12.center);		
		\draw[red,bend left=40,-<-=.25] (10.center) to (12.center);		
		
		\node[left=2mm of 5] () {\textcolor{blue}{$\mathbf{x}$}};

		\node[below left=15mm and 3mm of 5] () {(a)};

	\end{tikzpicture}
\end{scaletikzpicturetowidth}\vspace{2.5em}\\
		\begin{scaletikzpicturetowidth}{\textwidth}
\begin{tikzpicture}[scale=\tikzscale]

		\node (0) at (0,0) {};

		\node[draw=magenta,fill=magenta,circle,minimum size=2mm,right=00mm of 0] (1) {};	
		\node[draw=magenta,fill=magenta,circle,minimum size=2mm,right=30mm of 0] (2) {};	
		\node[draw=magenta,fill=magenta,circle,minimum size=2mm,right=60mm of 0] (3) {};	
		\node[draw=magenta,fill=magenta,circle,minimum size=2mm,right=90mm of 0] (4) {};

		\draw[orange,bend left=15,->-=.75]  (1.center) to (4.center);
		\draw[orange,bend left=15,-<-=.25]  (1.center) to (4.center);
		\draw[orange,bend right=15,->-=.75]  (1.center) to (3.center);
		\draw[orange,bend right=15,-<-=.25]  (1.center) to (3.center);
		\draw[orange,bend right=15,->-=.75]  (2.center) to (4.center);
		\draw[orange,bend right=15,-<-=.25]  (2.center) to (4.center);

		\node[draw=lime,fill=lime,circle,minimum size=2mm,below=10mm of 1.center] (5) {};	
		\node[draw=lime,fill=lime,circle,minimum size=2mm,below=10mm of 2.center] (6) {};	
		\node[draw=lime,fill=lime,circle,minimum size=2mm,below=10mm of 3.center] (7) {};	
		\node[draw=lime,fill=lime,circle,minimum size=2mm,below=10mm of 4.center] (8) {};	

		\draw[black,->-=.75] (1.center) -- (5.center);
		\draw[black,->-=.75] (2.center) -- (6.center);
		\draw[black,->-=.75] (3.center) -- (7.center);
		\draw[black,->-=.75] (4.center) -- (8.center);

		\draw[black,->-=.75] (5.center) -- (1.center);
		\draw[black,->-=.75] (6.center) -- (2.center);
		\draw[black,->-=.75] (7.center) -- (3.center);
		\draw[black,->-=.75] (8.center) -- (4.center);

		\node[draw=magenta,fill=magenta,circle,minimum size=2mm,below right=10mm and -8mm of 5.center] (9) {};	
		\node[draw=magenta,fill=magenta,circle,minimum size=2mm,below right=10mm and 8mm of 5.center] (10) {};	
		\node[draw=magenta,fill=magenta,circle,minimum size=2mm,below right=10mm and -8mm of 6.center] (11) {};	
		\node[draw=magenta,fill=magenta,circle,minimum size=2mm,below right=10mm and 8mm of 6.center] (12) {};	
		\node[draw=magenta,fill=magenta,circle,minimum size=2mm,below right=10mm and -8mm of 7.center] (13) {};	
		\node[draw=magenta,fill=magenta,circle,minimum size=2mm,below right=10mm and 8mm of 7.center] (14) {};	
		\node[draw=magenta,fill=magenta,circle,minimum size=2mm,below right=10mm and -8mm of 8.center] (15) {};	
		\node[draw=magenta,fill=magenta,circle,minimum size=2mm,below right=10mm and 8mm of 8.center] (16) {};
		
		\draw[cyan,->-=.75] (5.center) -- (9.center);
		\draw[cyan,->-=.75] (5.center) -- (10.center);
		\draw[cyan,->-=.75] (6.center) -- (11.center);
		\draw[cyan,->-=.75] (6.center) -- (12.center);	
		\draw[cyan,->-=.75] (7.center) -- (13.center);
		\draw[cyan,->-=.75] (7.center) -- (14.center);	
		\draw[cyan,->-=.75] (8.center) -- (15.center);
		\draw[cyan,->-=.75] (8.center) -- (16.center);	

		\draw[cyan,->-=.75] (9.center) -- (5.center);
		\draw[cyan,->-=.75] (10.center) -- (5.center);
		\draw[cyan,->-=.75] (11.center) -- (6.center);
		\draw[cyan,->-=.75] (12.center) -- (6.center);	
		\draw[cyan,->-=.75] (13.center) -- (7.center);
		\draw[cyan,->-=.75] (14.center) -- (7.center);	
		\draw[cyan,->-=.75] (15.center) -- (8.center);
		\draw[cyan,->-=.75] (16.center) -- (8.center);	
		
		\draw[orange,bend left=15,->-=.75]  (9.center) to (12.center);
		\draw[orange,bend left=15,-<-=.25]  (9.center) to (12.center);
		\draw[orange,bend right=15,->-=.75] (9.center) to (11.center);
		\draw[orange,bend right=15,-<-=.25] (9.center) to (11.center);
		\draw[orange,bend right=15,->-=.75] (10.center) to (12.center);
		\draw[orange,bend right=15,-<-=.25] (10.center) to (12.center);
		\draw[orange,bend left=15,->-=.75]  (11.center) to (14.center);
		\draw[orange,bend left=15,-<-=.25]  (11.center) to (14.center);
		\draw[orange,bend right=15,->-=.75] (11.center) to (13.center);
		\draw[orange,bend right=15,-<-=.25] (11.center) to (13.center);
		\draw[orange,bend right=15,->-=.75] (12.center) to (14.center);
		\draw[orange,bend right=15,-<-=.25] (12.center) to (14.center);
		\draw[orange,bend left=15,->-=.75]  (13.center) to (16.center);
		\draw[orange,bend left=15,-<-=.25]  (13.center) to (16.center);
		\draw[orange,bend right=15,->-=.75] (13.center) to (15.center);
		\draw[orange,bend right=15,-<-=.25] (13.center) to (15.center);
		\draw[orange,bend right=15,->-=.75] (14.center) to (16.center);
		\draw[orange,bend right=15,-<-=.25] (14.center) to (16.center);

		\node[draw=lime,fill=lime,circle,minimum size=2mm,below=10mm of 9.center] (17) {};	
		\node[draw=lime,fill=lime,circle,minimum size=2mm,below=10mm of 10.center] (18) {};	
		\node[draw=lime,fill=lime,circle,minimum size=2mm,below=10mm of 11.center] (19) {};	
		\node[draw=lime,fill=lime,circle,minimum size=2mm,below=10mm of 12.center] (20) {};	
		\node[draw=lime,fill=lime,circle,minimum size=2mm,below=10mm of 13.center] (21) {};	
		\node[draw=lime,fill=lime,circle,minimum size=2mm,below=10mm of 14.center] (22) {};	
		\node[draw=lime,fill=lime,circle,minimum size=2mm,below=10mm of 15.center] (23) {};	
		\node[draw=lime,fill=lime,circle,minimum size=2mm,below=10mm of 16.center] (24) {};	

		\draw[black,->-=.75] (9.center) -- (17.center);
		\draw[black,->-=.75] (10.center) -- (18.center);
		\draw[black,->-=.75] (11.center) -- (19.center);
		\draw[black,->-=.75] (12.center) -- (20.center);	
		\draw[black,->-=.75] (13.center) -- (21.center);
		\draw[black,->-=.75] (14.center) -- (22.center);	
		\draw[black,->-=.75] (15.center) -- (23.center);
		\draw[black,->-=.75] (16.center) -- (24.center);

		\draw[black,->-=.75] (17.center) -- (9.center);
		\draw[black,->-=.75] (18.center) -- (10.center);
		\draw[black,->-=.75] (19.center) -- (11.center);
		\draw[black,->-=.75] (20.center) -- (12.center);	
		\draw[black,->-=.75] (21.center) -- (13.center);
		\draw[black,->-=.75] (22.center) -- (14.center);	
		\draw[black,->-=.75] (23.center) -- (15.center);
		\draw[black,->-=.75] (24.center) -- (16.center);
		
		\node[draw=blue,fill=blue,circle,minimum size=2mm,below=10mm of 17.center] (25) {};	
		\node[draw=blue,fill=blue,circle,minimum size=2mm,below=10mm of 18.center] (26) {};	
		\node[draw=blue,fill=blue,circle,minimum size=2mm,below=10mm of 19.center] (27) {};	
		\node[draw=blue,fill=blue,circle,minimum size=2mm,below=10mm of 20.center] (28) {};	
		\node[draw=blue,fill=blue,circle,minimum size=2mm,below=10mm of 21.center] (29) {};	
		\node[draw=blue,fill=blue,circle,minimum size=2mm,below=10mm of 22.center] (30) {};	
		\node[draw=blue,fill=blue,circle,minimum size=2mm,below=10mm of 23.center] (31) {};	
		\node[draw=blue,fill=blue,circle,minimum size=2mm,below=10mm of 24.center] (32) {};	

		\draw[cyan,->-=.75] (17.center) -- (25.center);
		\draw[cyan,->-=.75] (18.center) -- (26.center);
		\draw[cyan,->-=.75] (19.center) -- (27.center);
		\draw[cyan,->-=.75] (20.center) -- (28.center);	
		\draw[cyan,->-=.75] (21.center) -- (29.center);
		\draw[cyan,->-=.75] (22.center) -- (30.center);	
		\draw[cyan,->-=.75] (23.center) -- (31.center);
		\draw[cyan,->-=.75] (24.center) -- (32.center);

		\draw[cyan,->-=.75] (25.center) -- (17.center);
		\draw[cyan,->-=.75] (26.center) -- (18.center);
		\draw[cyan,->-=.75] (27.center) -- (19.center);
		\draw[cyan,->-=.75] (28.center) -- (20.center);	
		\draw[cyan,->-=.75] (29.center) -- (21.center);
		\draw[cyan,->-=.75] (30.center) -- (22.center);	
		\draw[cyan,->-=.75] (31.center) -- (23.center);
		\draw[cyan,->-=.75] (32.center) -- (24.center);
		
		\draw [red,->-=.4,loop below,looseness=10,min distance=5mm] (25.center) to (25.center); 
		\draw [red,->-=.4,loop below,looseness=10,min distance=5mm] (26.center) to (26.center); 
		\draw [red,->-=.4,loop below,looseness=10,min distance=5mm] (27.center) to (27.center); 
		\draw [red,->-=.4,loop below,looseness=10,min distance=5mm] (28.center) to (28.center); 
		\draw [red,->-=.4,loop below,looseness=10,min distance=5mm] (29.center) to (29.center); 
		\draw [red,->-=.4,loop below,looseness=10,min distance=5mm] (30.center) to (30.center); 
		\draw [red,->-=.4,loop below,looseness=10,min distance=5mm] (31.center) to (31.center); 
		\draw [red,->-=.4,loop below,looseness=10,min distance=5mm] (32.center) to (32.center); 
		
		\draw[red,->-=.75] (25.center) -- (26.center);
		\draw[red,->-=.75] (26.center) -- (27.center);
		\draw[red,->-=.75] (27.center) -- (28.center);
		\draw[red,->-=.75] (28.center) -- (29.center);
		\draw[red,->-=.75] (29.center) -- (30.center);
		\draw[red,->-=.75] (30.center) -- (31.center);
		\draw[red,->-=.75] (31.center) -- (32.center);

		\draw[red,->-=.75] (26.center) -- (25.center);
		\draw[red,->-=.75] (27.center) -- (26.center);
		\draw[red,->-=.75] (28.center) -- (27.center);
		\draw[red,->-=.75] (29.center) -- (28.center);
		\draw[red,->-=.75] (30.center) -- (29.center);
		\draw[red,->-=.75] (31.center) -- (30.center);
		\draw[red,->-=.75] (32.center) -- (31.center);
		
		\node[left=2mm of 25] () {\textcolor{blue}{$\mathbf{x}$}};
		\node[left=2mm of 17] () {\textcolor{lime}{$\mathbf{z}^{(3)}$}};
		\node[left=2mm of 9] () {\textcolor{blue}{$\tilde{\mathbf{x}}^{(2)}$} = \textcolor{magenta}{$\mathbf{y}^{(3)}$}};
		
		\node[left=9mm of 5] () {\textcolor{lime}{$\mathbf{z}^{(2)}$}};
		\node[left=9mm of 1] () {\textcolor{blue}{$\tilde{\mathbf{x}}^{(1)}$} = \textcolor{magenta}{$\mathbf{y}^{(2)}$}};

		\node[below left=5mm and 3mm of 25] () {(b)};
	
	\end{tikzpicture}
\end{scaletikzpicturetowidth}
	\end{center}
	\caption{The graph of (a)~the original dense system and (b)~the extended sparse system, corresponding to the 1D problem of \autoref{fig:FMM_1D_new}. The same color convention as in \autoref{fig:extended_matrix_new_tikz} is adopted: red edges in between $\mathbf{x}$ correspond to $\mathbf{S}^{(3)}$, cyan edges from $\mathbf{z}^{(3)}$ to $\mathbf{x}$ correspond to $\mathbf{U}^{(3)}$, cyan edges from $\mathbf{x}$ to $\mathbf{z}^{(3)}$ correspond to $\mathbf{V}^{(3)^\mathrm{T}}$, black edges from $\mathbf{y}^{(3)}$ to $\mathbf{z}^{(3)}$ correspond to $-\mathbf{I}$, black edges from $\mathbf{z}^{(3)}$ to $\mathbf{y}^{(3)}$ correspond to $-\mathbf{I}$, orange edges in between $\mathbf{y}^{(3)}$ correspond to $\mathbf{S}^{(2)}$, cyan edges from $\mathbf{z}^{(2)}$ to $\mathbf{y}^{(3)}$  correspond to $ \mathbf{U}^{(2)}$, cyan edges from $\mathbf{y}^{(3)}$ to $\mathbf{z}^{(2)}$ correspond to $ \mathbf{V}^{(2)^\mathrm{T}}$,  black edges from $\mathbf{y}^{(2)}$ to $\mathbf{z}^{(2)}$ correspond to $-\mathbf{I}$, black edges from $\mathbf{z}^{(2)}$ to $\mathbf{y}^{(2)}$ correspond to $-\mathbf{I}$, and orange edges in between $\mathbf{y}^{(2)}$ correspond to $\mathbf{K}^{(2)}$.}
	\label{fig:initial_graphs}
\end{figure}
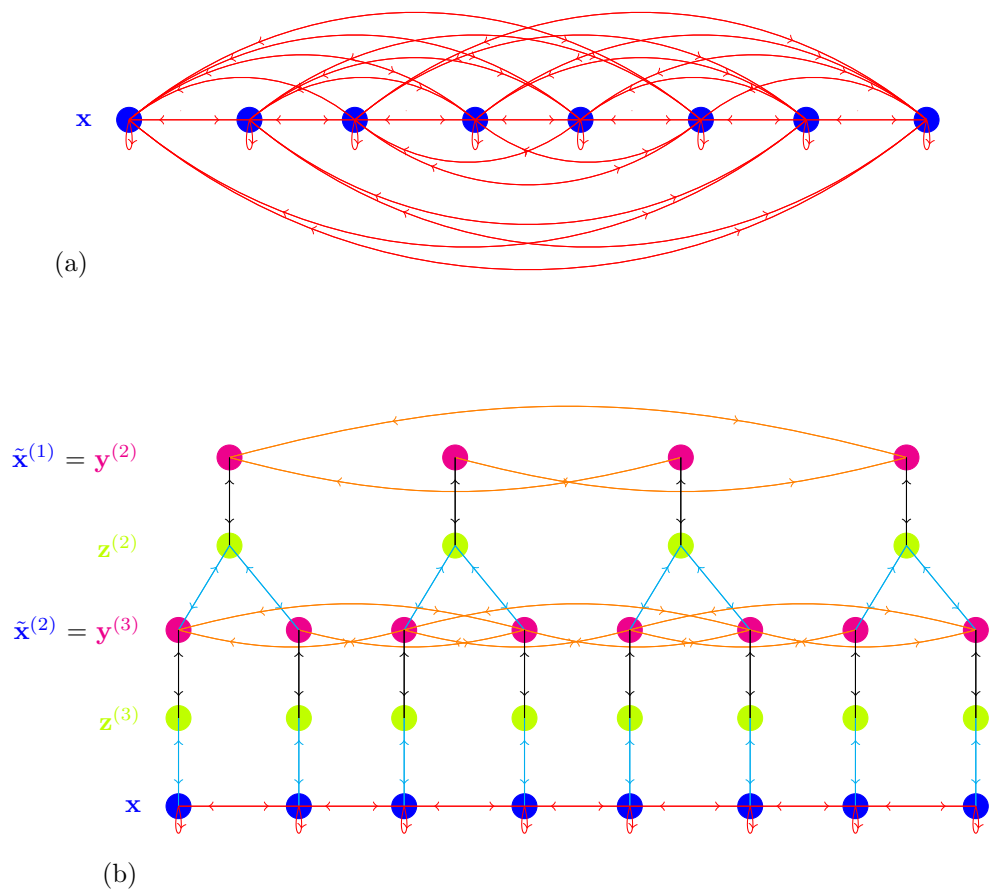

Coming back to the problem under concern, \autoref{fig:initial_graphs}(a) and~\ref{fig:initial_graphs}(b) show the graphs of the original dense system and the extended sparse system corresponding to the 1D problem introduced in \autoref{fig:FMM_1D_new}, respectively. Note that in- and outgoing edges have been combined into a single edge (with two arrows) in order to limit the number of edges in these figures. The layout of the graph is dictated by the specific ordering of equations that we have chosen in \autoref{eq:Ax=b_extended}, leading to a symmetric extended sparse matrix. The same color convention as in \autoref{fig:extended_matrix_new_tikz} is employed for the vertices and edges. The transformation from \autoref{fig:initial_graphs}(a) to \ref{fig:initial_graphs}(b) clearly demonstrates how one can take advantage of the $\mH^{2}$-structure of the matrix to obtain a sparsified system. \autoref{fig:initial_graphs}(b) also illustrates that the variables $\mathbf{y}^{(l)}$ at level $l$ can be considered as unknowns $\tilde{\mathbf{x}}^{(l-1)}$ at the parent level $l-1$, hence revealing the hierarchical structure of the $\mH^2$-matrix (with $\tilde{\mathbf{x}}^{(L)}$ corresponding to $\mathbf{x}$).


%
%

\subsection{Compression and redirection of fill-ins}\label{subsec:compression}
Transforming a dense $\mH^{2}$-matrix into an extended sparse matrix is a crucial but insufficient step to obtain a fast solver. Indeed, feeding a matrix such as the ones depicted in \autoref{fig:extended_matrix_new_tikz} and \autoref{fig:extended_matrix_2D} to a conventional sparse solver leads to the creation of dense fill-ins throughout the elimination phase (or, equivalently, additional edges in the graph of \autoref{fig:initial_graphs}(b)), which is harmful for the scaling of the algorithm (see for example the thesis of Pals~\cite{pals04a}). In order to preserve the initial sparsity pattern, a second key idea is therefore introduced: fill-ins corresponding to interactions between well-separated clusters are compressed into a low-rank representation and are subsequently redirected through other operators/edges. This procedure is discussed step by step in the next paragraphs, focusing on a graph representation of the methodology. To this end, the one-dimensional example introduced in \autoref{subsec:extended_sparsification} is abandoned and a more general case is now considered. References to the original FMM-terminology (e.g., particles, multipoles, and locals) are omitted to stress the generality of the method. 

\begin{figure}[hbtp]
	\begin{center}
		\begin{tikzpicture}[scale=2.5]

		\node (0) at (0,0) {};
		
		\node[draw=lime,fill=lime,circle,minimum size=2mm,right=15mm of 0] (1) {};			
		\draw[black,-<-=.75,dashed] (1.center) -- ++(90:2mm);
		\draw[cyan,-<-=.75,dashed] (1.center) -- ++(270:2mm);
		\draw[cyan,-<-=.75,dashed] (1.center) -- ++(315:2mm);
		
		\node[draw=magenta,fill=magenta,circle,minimum size=2mm,below left=10mm and 10mm of 1.center] (2) {};			
		\draw[cyan,->-=.75] (1.center) -- (2.center);
		\draw[cyan,->-=.75] (2.center) -- (1.center);
		
		\draw[orange,-<-=.75,dashed] (2.center) -- ++(22.5:2mm);
		\draw[orange,-<-=.75,dashed] (2.center) -- ++(157.5:2mm);
		\draw[orange,-<-=.75,dashed] (2.center) -- ++(202.5:2mm);
		\draw[orange,-<-=.75,dashed] (2.center) -- ++(337.5:2mm);
		
		\node[draw=lime,fill=lime,circle,minimum size=2mm,below=10mm of 2.center] (3) {};	
		\draw[black,->-=.75] (2.center) -- (3.center);
		\draw[black,->-=.75] (3.center) -- (2.center);

		\node[draw=blue,fill=blue,circle,minimum size=2mm,below=10mm of 3.center] (4) {};
		\draw[cyan,->-=.75] (3.center) -- (4.center);
		\draw[cyan,->-=.75] (4.center) -- (3.center);
		\draw [red,->-=.4,loop below,looseness=10,min distance=3mm] (4.center) to (4.center);
		\draw[red,-<-=.75,dashed] (4.center) -- ++(0:2mm);
		\draw[red,-<-=.75,dashed] (4.center) -- ++(22.5:2mm);
		\draw[red,-<-=.75,dashed] (4.center) -- ++(180:2mm);
		\draw[red,-<-=.75,dashed] (4.center) -- ++(157.5:2mm);
		\draw[red,-<-=.75,dashed] (4.center) -- ++(202.5:2mm);
		\draw[red,-<-=.75,dashed] (4.center) -- ++(337.5:2mm);

		\node[below left=-3mm and 25mm of 1.center] () {\textcolor{lime}{$\mathbf{z}_{\mathrm{p}}^{(l-1)}$}};
		\node[below left=8mm and 25mm of 1.center] () {\textcolor{magenta}{$\mathbf{y}_{i}^{(l)}$}};
		\node[below left=20mm and 25mm of 1.center] () {\textcolor{lime}{$\mathbf{z}_{i}^{(l)}$}};
		\node[below left=31mm and 25mm of 1.center] () {\textcolor{blue}{$\tilde{\mathbf{x}}_{i}^{(l)}$}};

		\node[below right=-3mm and 15mm of 1.center] () {$\texttt{i}^{\dagger^{(l)}} = \texttt{p}^{(l-1)}$};
		\node[below right=8mm and 15mm of 1.center] () {$\overline{\texttt{i}}^{(l)}$};
		\node[below right=20mm and 15mm of 1.center] () {$\texttt{i}^{(l)}$};
		\node[below right=32mm and 15mm of 1.center] () {$\underline{\texttt{i}}^{(l)}$};

		\node[below = 7mm of 4.center] () {$\Omega_{i}^{(l)}$};
	
	\end{tikzpicture}
	\end{center}
	\caption{Alternative notation for the nodes related to cluster~$\Omega_{i}^{(l)}$.}
	\label{fig:definitions_tikz}
\end{figure}
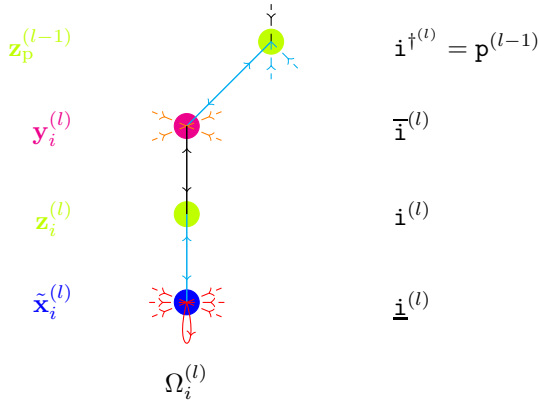

The following simplifying notations are introduced (see \autoref{fig:definitions_tikz}): nodes in a graph corresponding to the variables~$\tilde{\mathbf{x}}_{i}^{(l)}$, $\mathbf{z}_{i}^{(l)}$, and $\mathbf{y}_{i}^{(l)}$ of a cluster $\Omega_{i}^{(l)}$ are referred to as nodes $\underline{\texttt{i}}^{(l)}$, $\texttt{i}^{(l)}$, and $\overline{\texttt{i}}^{(l)}$, respectively (the position of the bar is a mnemonic corresponding to the location of the node in the tree, either below or above). The symbol $\texttt{i}^{\dagger^{(l)}}$ represents the variable $\mathbf{z}_{\mathrm{p}}^{(l-1)}$ of its parent cluster $\Omega_{\mathrm{p}}^{(l-1)}$. The latter can also be denoted as $\texttt{p}^{(l-1)}$. If no confusion is possible, the superscripts referring to the hierarchical level of a cluster will be omitted in the following for sake of readability. With this notation in place, operators such as the inter- and anterpolation\footnote{Interpolation = from coarse to fine. Anterpolation = from fine to coarse.} matrices $\mathbf{U}_{i}$ and $\mathbf{V}^{\mathrm{T}}_{i}$ of a cluster $\Omega_{i}$ are written as $ \mathbf{U} ( \underline{\texttt{i}}, \texttt{i} )$ and $\mathbf{V}^{\mathrm{T}}(\texttt{i}, \underline{\texttt{i}} )$, respectively.

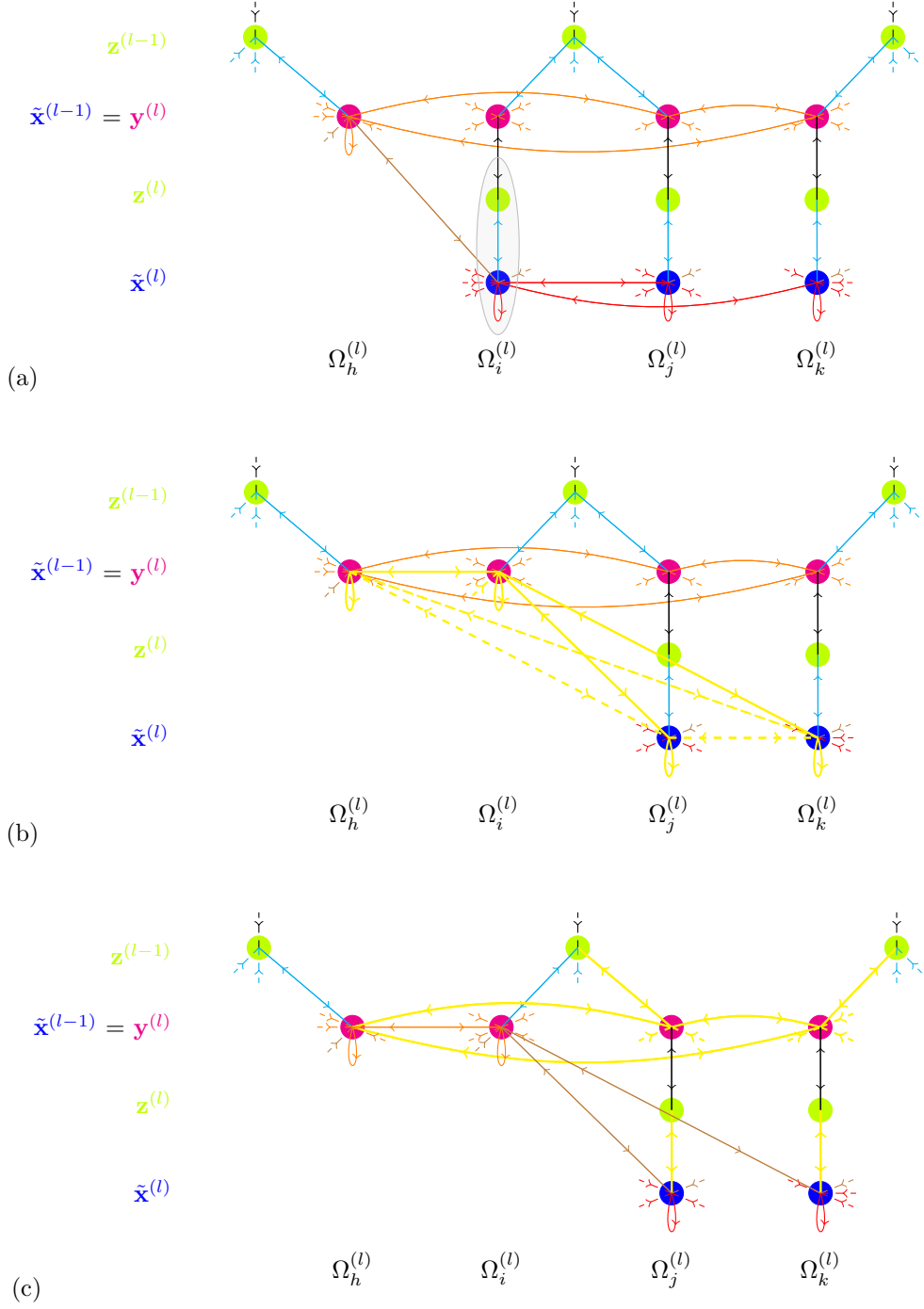
\begin{figure}[hbtp]
	\begin{center}
		(a)~\hspace{-1em}\begin{tikzpicture}[scale=2.5]

		\node (0) at (0,0) {};

		\node[draw=lime,fill=lime,circle,minimum size=2mm,right=00mm of 0] (1) {};	
		\node[draw=lime,fill=lime,circle,minimum size=2mm,right=45mm of 0] (2) {};	
		\node[draw=lime,fill=lime,circle,minimum size=2mm,right=90mm of 0] (3) {};	

		\draw[black,-<-=.75,dashed] (1.center) -- ++(90:2mm);
		\draw[black,-<-=.75,dashed] (2.center) -- ++(90:2mm);
		\draw[black,-<-=.75,dashed] (3.center) -- ++(90:2mm);

		\node[draw=magenta,fill=magenta,circle,minimum size=2mm,below right=10mm and 12mm of 1.center] (4) {};	
		\node[draw=magenta,fill=magenta,circle,minimum size=2mm,below right=10mm and -12mm of 2.center] (5) {};	
		\node[draw=magenta,fill=magenta,circle,minimum size=2mm,below right=10mm and 12mm of 2.center] (6) {};	
		\node[draw=magenta,fill=magenta,circle,minimum size=2mm,below right=10mm and -12mm of 3.center] (7) {};	
		
		\draw[cyan,->-=.75] (1.center) -- (4.center);
		\draw[cyan,->-=.75] (2.center) -- (5.center);
		\draw[cyan,->-=.75] (2.center) -- (6.center);
		\draw[cyan,->-=.75] (3.center) -- (7.center);	

		\draw[cyan,->-=.75] (4.center) -- (1.center);
		\draw[cyan,->-=.75] (5.center) -- (2.center);
		\draw[cyan,->-=.75] (6.center) -- (2.center);
		\draw[cyan,->-=.75] (7.center) -- (3.center);	

		\draw[cyan,-<-=.75,dashed] (1.center) -- ++(270:2mm);
		\draw[cyan,-<-=.75,dashed] (1.center) -- ++(225:2mm);	
		\draw[cyan,-<-=.75,dashed] (2.center) -- ++(270:2mm);
		\draw[cyan,-<-=.75,dashed] (3.center) -- ++(270:2mm);
		\draw[cyan,-<-=.75,dashed] (3.center) -- ++(315:2mm);
		
		\draw[orange,bend left=15,->-=.75]  (4.center) to (6.center);
		\draw[orange,bend left=15,-<-=.25]  (4.center) to (6.center);
		\draw[orange,bend left=15,->-=.75] (6.center) to (7.center);
		\draw[orange,bend left=15,-<-=.25] (6.center) to (7.center);
		\draw[orange,bend right=15,->-=.75] (4.center) to (7.center);
		\draw[orange,bend right=15,-<-=.25] (4.center) to (7.center);
		
		\draw[orange,-<-=.75,dashed] (4.center) -- ++(337.5:2mm);
		\draw[orange,-<-=.75,dashed] (4.center) -- ++(157.5:2mm);
		\draw[orange,-<-=.75,dashed] (4.center) -- ++(180:2mm);
		\draw[orange,->-=.4,loop below,looseness=10,min distance=3mm] (4.center) to (4.center);
		\draw[brown,-<-=.75,dashed] (4.center) -- ++(225:2mm);
	
		\draw[orange,-<-=.75,dashed] (5.center) -- ++(337.5:2mm);
		\draw[orange,-<-=.75,dashed] (5.center) -- ++(157.5:2mm);
		\draw[orange,-<-=.75,dashed] (5.center) -- ++(202.5:2mm);
		\draw[orange,-<-=.75,dashed] (5.center) -- ++(22.5:2mm);

		\draw[orange,-<-=.75,dashed] (6.center) -- ++(202.5:2mm);
		\draw[orange,-<-=.75,dashed] (6.center) -- ++(337.5:2mm);

		\draw[orange,-<-=.75,dashed] (7.center) -- ++(22.5:2mm);
		\draw[orange,-<-=.75,dashed] (7.center) -- ++(202.5:2mm);
		\draw[orange,-<-=.75,dashed] (7.center) -- ++(337.5:2mm);

		\node[draw=lime,fill=lime,circle,minimum size=2mm,below=10mm of 5.center] (8) {};	
		\node[draw=lime,fill=lime,circle,minimum size=2mm,below=10mm of 6.center] (9) {};	
		\node[draw=lime,fill=lime,circle,minimum size=2mm,below=10mm of 7.center] (10) {};	

		\draw[black,->-=.75] (5.center) -- (8.center);
		\draw[black,->-=.75] (6.center) -- (9.center);
		\draw[black,->-=.75] (7.center) -- (10.center);

		\draw[black,->-=.75] (8.center) -- (5.center);
		\draw[black,->-=.75] (9.center) -- (6.center);
		\draw[black,->-=.75] (10.center) -- (7.center);
		
		\node[draw=blue,fill=blue,circle,minimum size=2mm,below=10mm of 8.center] (11) {};	
		\node[draw=blue,fill=blue,circle,minimum size=2mm,below=10mm of 9.center] (12) {};	
		\node[draw=blue,fill=blue,circle,minimum size=2mm,below=10mm of 10.center] (13) {};	

		\draw[cyan,->-=.75] (8.center) -- (11.center);
		\draw[cyan,->-=.75] (9.center) -- (12.center);
		\draw[cyan,->-=.75] (10.center) -- (13.center);

		\draw[cyan,->-=.75] (11.center) -- (8.center);
		\draw[cyan,->-=.75] (12.center) -- (9.center);
		\draw[cyan,->-=.75] (13.center) -- (10.center);
	
		\draw[brown,->-=.75] (4.center) -- (11.center);
		\draw[brown,->-=.75] (11.center) -- (4.center);
		
		\draw [red,->-=.4,loop below,looseness=10,min distance=3mm] (11.center) to (11.center); 
		\draw [red,->-=.4,loop below,looseness=10,min distance=3mm] (12.center) to (12.center); 
		\draw [red,->-=.4,loop below,looseness=10,min distance=3mm] (13.center) to (13.center);

		\draw[red,->-=.75] (11.center) -- (12.center);
		\draw[red,->-=.75] (12.center) -- (11.center);
		\draw[red,bend right=15,->-=.75]  (11.center) to (13.center);
		\draw[red,bend right=15,-<-=.25]  (11.center) to (13.center);

		\draw[brown,-<-=.75,dashed] (11.center) -- ++(22.5:2mm);
		\draw[red,-<-=.75,dashed] (11.center) -- ++(157.5:2mm);
		\draw[red,-<-=.75,dashed] (11.center) -- ++(180:2mm);
		\draw[red,-<-=.75,dashed] (11.center) -- ++(202.5:2mm);

		\draw[brown,-<-=.75,dashed] (12.center) -- ++(22.5:2mm);
		\draw[red,-<-=.75,dashed] (12.center) -- ++(157.5:2mm);
		\draw[red,-<-=.75,dashed] (12.center) -- ++(202.5:2mm);
		\draw[red,-<-=.75,dashed] (12.center) -- ++(337.5:2mm);
	
		\draw[red,-<-=.75,dashed] (13.center) -- ++(0:2mm);
		\draw[brown,-<-=.75,dashed] (13.center) -- ++(22.5:2mm);
		\draw[red,-<-=.75,dashed] (13.center) -- ++(157.5:2mm);
		\draw[red,-<-=.75,dashed] (13.center) -- ++(337.5:2mm);

		
		\node[below left=-2mm and 11mm of 1.center] () {\textcolor{lime}{$\mathbf{z}^{(l-1)}$}};
		\node[below left=8mm and 11mm of 1.center] () {\textcolor{blue}{$\tilde{\mathbf{x}}^{(l-1)}$} = \textcolor{magenta}{$\mathbf{y}^{(l)}$}};
		\node[below left=19mm and 11mm of 1.center] () {\textcolor{lime}{$\mathbf{z}^{(l)}$}};
		\node[below left=31.5mm and 11mm of 1.center] () {\textcolor{blue}{$\tilde{\mathbf{x}}^{(l)}$}};

		\node[draw=lightgray,fill=lightgray,ellipse,minimum height=25mm,minimum width=6mm,fill opacity=0.1,below=-6mm of 8.center] (14) {};	

%
		\node[below = 30.5mm of 4.center] () {$\Omega_{h}^{(l)}$};
		\node[below = 7mm of 11.center] () {$\Omega_{i}^{(l)}$};
		\node[below = 7mm of 12.center] () {$\Omega_{j}^{(l)}$};
		\node[below = 7mm of 13.center] () {$\Omega_{k}^{(l)}$};
		
	
	\end{tikzpicture}\vspace{2.5em}\\
		(b)~\hspace{-1em}\begin{tikzpicture}[scale=2.5]

		\node (0) at (0,0) {};

		\node[draw=lime,fill=lime,circle,minimum size=2mm,right=00mm of 0] (1) {};	
		\node[draw=lime,fill=lime,circle,minimum size=2mm,right=45mm of 0] (2) {};	
		\node[draw=lime,fill=lime,circle,minimum size=2mm,right=90mm of 0] (3) {};	

		\draw[black,-<-=.75,dashed] (1.center) -- ++(90:2mm);
		\draw[black,-<-=.75,dashed] (2.center) -- ++(90:2mm);
		\draw[black,-<-=.75,dashed] (3.center) -- ++(90:2mm);

		\node[draw=magenta,fill=magenta,circle,minimum size=2mm,below right=10mm and 12mm of 1.center] (4) {};	
		\node[draw=magenta,fill=magenta,circle,minimum size=2mm,below right=10mm and -12mm of 2.center] (5) {};	
		\node[draw=magenta,fill=magenta,circle,minimum size=2mm,below right=10mm and 12mm of 2.center] (6) {};	
		\node[draw=magenta,fill=magenta,circle,minimum size=2mm,below right=10mm and -12mm of 3.center] (7) {};	
		
		\draw[cyan,->-=.75] (1.center) -- (4.center);
		\draw[cyan,->-=.75] (2.center) -- (5.center);
		\draw[cyan,->-=.75] (2.center) -- (6.center);
		\draw[cyan,->-=.75] (3.center) -- (7.center);

		\draw[cyan,->-=.75] (4.center) -- (1.center);
		\draw[cyan,->-=.75] (5.center) -- (2.center);
		\draw[cyan,->-=.75] (6.center) -- (2.center);
		\draw[cyan,->-=.75] (7.center) -- (3.center);
		
		\draw[cyan,-<-=.75,dashed] (1.center) -- ++(270:2mm);
		\draw[cyan,-<-=.75,dashed] (1.center) -- ++(225:2mm);	
		\draw[cyan,-<-=.75,dashed] (2.center) -- ++(270:2mm);
		\draw[cyan,-<-=.75,dashed] (3.center) -- ++(270:2mm);
		\draw[cyan,-<-=.75,dashed] (3.center) -- ++(315:2mm);
		
		\draw[orange,bend left=15,->-=.75]  (4.center) to (6.center);
		\draw[orange,bend left=15,-<-=.25]  (4.center) to (6.center);
		\draw[orange,bend left=15,->-=.75] (6.center) to (7.center);
		\draw[orange,bend left=15,-<-=.25] (6.center) to (7.center);
		\draw[orange,bend right=15,->-=.75] (4.center) to (7.center);
		\draw[orange,bend right=15,-<-=.25] (4.center) to (7.center);
		
		\draw[orange,-<-=.75,dashed] (4.center) -- ++(337.5:2mm);
		\draw[orange,-<-=.75,dashed] (4.center) -- ++(157.5:2mm);
		\draw[orange,-<-=.75,dashed] (4.center) -- ++(180:2mm);
		\draw[yellow,thick,->-=.4,loop below,looseness=10,min distance=3mm] (4.center) to (4.center);
		\draw[brown,-<-=.75,dashed] (4.center) -- ++(225:2mm);
	
		\draw[orange,-<-=.75,dashed] (5.center) -- ++(337.5:2mm);
		\draw[orange,-<-=.75,dashed] (5.center) -- ++(157.5:2mm);
		\draw[orange,-<-=.75,dashed] (5.center) -- ++(202.5:2mm);
		\draw[orange,-<-=.75,dashed] (5.center) -- ++(22.5:2mm);
		\draw[yellow,thick,-<-=.75,dashed] (5.center) -- ++(222.5:2mm);

		\draw[orange,-<-=.75,dashed] (6.center) -- ++(202.5:2mm);
		\draw[orange,-<-=.75,dashed] (6.center) -- ++(337.5:2mm);

		\draw[orange,-<-=.75,dashed] (7.center) -- ++(22.5:2mm);
		\draw[orange,-<-=.75,dashed] (7.center) -- ++(202.5:2mm);
		\draw[orange,-<-=.75,dashed] (7.center) -- ++(337.5:2mm);

		\node[draw=lime,fill=lime,circle,minimum size=2mm,below=10mm of 6.center] (9) {};	
		\node[draw=lime,fill=lime,circle,minimum size=2mm,below=10mm of 7.center] (10) {};	

		\draw[black,->-=.75] (6.center) -- (9.center);
		\draw[black,->-=.75] (7.center) -- (10.center);

		\draw[black,->-=.75] (9.center) -- (6.center);
		\draw[black,->-=.75] (10.center) -- (7.center);
		
		\node[draw=blue,fill=blue,circle,minimum size=2mm,below=10mm of 9.center] (12) {};	
		\node[draw=blue,fill=blue,circle,minimum size=2mm,below=10mm of 10.center] (13) {};	

		\draw[cyan,->-=.75] (9.center) -- (12.center);
		\draw[cyan,->-=.75] (10.center) -- (13.center);

		\draw[cyan,->-=.75] (12.center) -- (9.center);
		\draw[cyan,->-=.75] (13.center) -- (10.center);
	
		
		\draw [yellow,thick,->-=.4,loop below,looseness=10,min distance=3mm] (12.center) to (12.center); 
		\draw [yellow,thick,->-=.4,loop below,looseness=10,min distance=3mm] (13.center) to (13.center);



		\draw[brown,-<-=.75,dashed] (12.center) -- ++(22.5:2mm);
		\draw[red,-<-=.75,dashed] (12.center) -- ++(157.5:2mm);
		\draw[red,-<-=.75,dashed] (12.center) -- ++(202.5:2mm);
		\draw[red,-<-=.75,dashed] (12.center) -- ++(337.5:2mm);
	
		\draw[red,-<-=.75,dashed] (13.center) -- ++(0:2mm);
		\draw[brown,-<-=.75,dashed] (13.center) -- ++(22.5:2mm);
		\draw[red,-<-=.75,dashed] (13.center) -- ++(157.5:2mm);
		\draw[red,-<-=.75,dashed] (13.center) -- ++(337.5:2mm);

		
		\node[below left=-2mm and 11mm of 1.center] () {\textcolor{lime}{$\mathbf{z}^{(l-1)}$}};
		\node[below left=8mm and 11mm of 1.center] () {\textcolor{blue}{$\tilde{\mathbf{x}}^{(l-1)}$} = \textcolor{magenta}{$\mathbf{y}^{(l)}$}};
		\node[below left=19mm and 11mm of 1.center] () {\textcolor{lime}{$\mathbf{z}^{(l)}$}};
		\node[below left=31.5mm and 11mm of 1.center] () {\textcolor{blue}{$\tilde{\mathbf{x}}^{(l)}$}};


%
		\draw[yellow,thick,->-=.75] (4.center) -- (5.center);
		\draw[yellow,thick,->-=.75] (5.center) -- (4.center);
		\draw[yellow,thick,->-=.4,loop below,looseness=10,min distance=3mm] (5.center) to (5.center); 
		\draw[yellow,thick,->-=.75] (5.center) -- (12.center);
		\draw[yellow,thick,->-=.75] (12.center) -- (5.center);
		\draw[yellow,thick,->-=.75] (5.center) -- (13.center);
		\draw[yellow,thick,->-=.75] (13.center) -- (5.center);
		
		\draw[yellow,thick,->-=.75,dashed] (12.center) -- (13.center);
		\draw[yellow,thick,->-=.75,dashed] (4.center) -- (12.center);
		\draw[yellow,thick,->-=.75,dashed] (4.center) -- (13.center);

		\draw[yellow,thick,-<-=.25,dashed] (12.center) -- (13.center);
		\draw[yellow,thick,-<-=.25,dashed] (4.center) -- (12.center);
		\draw[yellow,thick,->-=.75,dashed] (13.center) -- (4.center);

		\node[below = 30.5mm of 4.center] () {$\Omega_{h}^{(l)}$};
		\node[below = 7mm of 11.center] () {$\Omega_{i}^{(l)}$};
		\node[below = 7mm of 12.center] () {$\Omega_{j}^{(l)}$};
		\node[below = 7mm of 13.center] () {$\Omega_{k}^{(l)}$};
	
	\end{tikzpicture}\vspace{2.5em}\\
		(c)~\hspace{-1em}\begin{tikzpicture}[scale=2.5]

		\node (0) at (0,0) {};

		\node[draw=lime,fill=lime,circle,minimum size=2mm,right=00mm of 0] (1) {};	
		\node[draw=lime,fill=lime,circle,minimum size=2mm,right=45mm of 0] (2) {};	
		\node[draw=lime,fill=lime,circle,minimum size=2mm,right=90mm of 0] (3) {};	

		\draw[black,-<-=.75,dashed] (1.center) -- ++(90:2mm);
		\draw[black,-<-=.75,dashed] (2.center) -- ++(90:2mm);
		\draw[black,-<-=.75,dashed] (3.center) -- ++(90:2mm);

		\node[draw=magenta,fill=magenta,circle,minimum size=2mm,below right=10mm and 12mm of 1.center] (4) {};	
		\node[draw=magenta,fill=magenta,circle,minimum size=2mm,below right=10mm and -12mm of 2.center] (5) {};	
		\node[draw=magenta,fill=magenta,circle,minimum size=2mm,below right=10mm and 12mm of 2.center] (6) {};	
		\node[draw=magenta,fill=magenta,circle,minimum size=2mm,below right=10mm and -12mm of 3.center] (7) {};	
		
		\draw[cyan,->-=.75] (1.center) -- (4.center);
		\draw[cyan,->-=.75] (2.center) -- (5.center);
		\draw[yellow,thick,->-=.75] (2.center) -- (6.center);
		\draw[yellow,thick,->-=.75] (3.center) -- (7.center);	
		
		\draw[cyan,->-=.75] (4.center) -- (1.center);
		\draw[cyan,->-=.75] (5.center) -- (2.center);
		\draw[yellow,thick,->-=.75] (6.center) -- (2.center);
		\draw[yellow,thick,->-=.75] (7.center) -- (3.center);

		\draw[cyan,-<-=.75,dashed] (1.center) -- ++(270:2mm);
		\draw[cyan,-<-=.75,dashed] (1.center) -- ++(225:2mm);	
		\draw[cyan,-<-=.75,dashed] (2.center) -- ++(270:2mm);
		\draw[cyan,-<-=.75,dashed] (3.center) -- ++(270:2mm);
		\draw[cyan,-<-=.75,dashed] (3.center) -- ++(315:2mm);
		
		\draw[yellow,thick,bend left=15,->-=.75]  (4.center) to (6.center);
		\draw[yellow,thick,bend left=15,-<-=.25]  (4.center) to (6.center);
		\draw[yellow,thick,bend left=15,->-=.75] (6.center) to (7.center);
		\draw[yellow,thick,bend left=15,-<-=.25] (6.center) to (7.center);
		\draw[yellow,thick,bend right=15,->-=.75] (4.center) to (7.center);
		\draw[yellow,thick,bend right=15,-<-=.25] (4.center) to (7.center);

		\draw[orange,-<-=.75,dashed] (4.center) -- ++(337.5:2mm);
		\draw[orange,-<-=.75,dashed] (4.center) -- ++(157.5:2mm);
		\draw[orange,-<-=.75,dashed] (4.center) -- ++(180:2mm);
		\draw[orange,->-=.4,loop below,looseness=10,min distance=3mm] (4.center) to (4.center);
		\draw[brown,-<-=.75,dashed] (4.center) -- ++(225:2mm);
	
		\draw[orange,-<-=.75,dashed] (5.center) -- ++(337.5:2mm);
		\draw[orange,-<-=.75,dashed] (5.center) -- ++(157.5:2mm);
		\draw[orange,-<-=.75,dashed] (5.center) -- ++(202.5:2mm);
		\draw[orange,-<-=.75,dashed] (5.center) -- ++(22.5:2mm);
		\draw[brown,-<-=.75,dashed] (5.center) -- ++(222.5:2mm);

		\draw[yellow,thick,-<-=.75,dashed] (6.center) -- ++(202.5:2mm);
		\draw[yellow,thick,-<-=.75,dashed] (6.center) -- ++(337.5:2mm);

		\draw[yellow,thick,-<-=.75,dashed] (7.center) -- ++(22.5:2mm);
		\draw[yellow,thick,-<-=.75,dashed] (7.center) -- ++(202.5:2mm);
		\draw[yellow,thick,-<-=.75,dashed] (7.center) -- ++(337.5:2mm);

		\node[draw=lime,fill=lime,circle,minimum size=2mm,below=10mm of 6.center] (9) {};	
		\node[draw=lime,fill=lime,circle,minimum size=2mm,below=10mm of 7.center] (10) {};	

		\draw[black,->-=.75] (6.center) -- (9.center);
		\draw[black,->-=.75] (7.center) -- (10.center);

		\draw[black,->-=.75] (9.center) -- (6.center);
		\draw[black,->-=.75] (10.center) -- (7.center);

		\node[draw=blue,fill=blue,circle,minimum size=2mm,below=10mm of 9.center] (12) {};	
		\node[draw=blue,fill=blue,circle,minimum size=2mm,below=10mm of 10.center] (13) {};	

		\draw[yellow,thick,->-=.75] (9.center) -- (12.center);
		\draw[yellow,thick,->-=.75] (10.center) -- (13.center);

		\draw[yellow,thick,->-=.75] (12.center) -- (9.center);
		\draw[yellow,thick,->-=.75] (13.center) -- (10.center);
	
		
		\draw [red,->-=.4,loop below,looseness=10,min distance=3mm] (12.center) to (12.center); 
		\draw [red,->-=.4,loop below,looseness=10,min distance=3mm] (13.center) to (13.center);



		\draw[brown,-<-=.75,dashed] (12.center) -- ++(22.5:2mm);
		\draw[red,-<-=.75,dashed] (12.center) -- ++(157.5:2mm);
		\draw[red,-<-=.75,dashed] (12.center) -- ++(202.5:2mm);
		\draw[red,-<-=.75,dashed] (12.center) -- ++(337.5:2mm);
	
		\draw[red,-<-=.75,dashed] (13.center) -- ++(0:2mm);
		\draw[brown,-<-=.75,dashed] (13.center) -- ++(22.5:2mm);
		\draw[red,-<-=.75,dashed] (13.center) -- ++(157.5:2mm);
		\draw[red,-<-=.75,dashed] (13.center) -- ++(337.5:2mm);

		
		\node[below left=-2mm and 11mm of 1.center] () {\textcolor{lime}{$\mathbf{z}^{(l-1)}$}};
		\node[below left=8mm and 11mm of 1.center] () {\textcolor{blue}{$\tilde{\mathbf{x}}^{(l-1)}$} = \textcolor{magenta}{$\mathbf{y}^{(l)}$}};
		\node[below left=19mm and 11mm of 1.center] () {\textcolor{lime}{$\mathbf{z}^{(l)}$}};
		\node[below left=31.5mm and 11mm of 1.center] () {\textcolor{blue}{$\tilde{\mathbf{x}}^{(l)}$}};


%
		\draw[orange,->-=.75] (4.center) -- (5.center);
		\draw[orange,->-=.75] (5.center) -- (4.center);
		\draw[orange,->-=.4,loop below,looseness=10,min distance=3mm] (5.center) to (5.center); 
		\draw[brown,->-=.75] (5.center) -- (12.center);
		\draw[brown,->-=.75] (5.center) -- (13.center);
		\draw[brown,->-=.75] (12.center) -- (5.center);
		\draw[brown,->-=.75] (13.center) -- (5.center);
		

		\node[below = 30.5mm of 4.center] () {$\Omega_{h}^{(l)}$};
		\node[below = 7mm of 11.center] () {$\Omega_{i}^{(l)}$};
		\node[below = 7mm of 12.center] () {$\Omega_{j}^{(l)}$};
		\node[below = 7mm of 13.center] () {$\Omega_{k}^{(l)}$};
	
	\end{tikzpicture}
	\end{center}
	\caption{Subgraph (a)~before and (b)~after the elimination of the variables $\tilde{\mathbf{x}}_{i}^{(l)}$ and $\mathbf{z}_{i}^{(l)}$. New and modified edges are highlighted in yellow in~(b). Operators associated with edges between well-separated clusters (indicated by dashed yellow lines in (b)) are compressed and redirected, resulting in the subgraph shown in (c). The edges that are updated after compression are highlighted in yellow in (c). This strategy of compressing and redirecting dense fill-ins is similar to the initial sparsification approach introduced in \autoref{subsec:extended_sparsification}, as low-rank interactions are represented at the parent level.}
	\label{fig:subgraph_tikz}
\end{figure}

\autoref{fig:subgraph_tikz}(a) shows part of a graph of an extended sparse matrix at a certain stage of the elimination phase of level $l$, i.e., after elimination of levels $L$, $L-1$, \ldots, $l+1$. This graph shows a ``generic'' description of the process of removing fill-in. The nodes shown are generic nodes in the graph of the matrix.

The variables $\tilde{\mathbf{x}}^{(l)}$ and $\mathbf{z}^{(l)}$ of some clusters at level $l$ have already been eliminated at this point, e.g., $\tilde{\mathbf{x}}_{h}^{(l)}$ and $\mathbf{z}_{h}^{(l)}$ of cluster $\Omega_{h}^{(l)}$ (or, equivalently, the nodes $\underline{\texttt{h}}$ and $\texttt{h}$). This has created edges in the graph between neighboring clusters, e.g., connections $\left( \overline{\texttt{h}},\underline{\texttt{i}} \right)$ and $\left( \underline{\texttt{i}},\overline{\texttt{h}} \right)$ between $\mathbf{y}_{h}^{(l)}$ and $\tilde{\mathbf{x}}_{i}^{(l)}$ (and vice versa). Additional edges such as $\left( \overline{\texttt{h}}, \overline{\texttt{h}} \right)$ (i.e., between $\mathbf{y}_{h}^{(l)}$ and itself) have been established as well. The variables $\tilde{\mathbf{x}}_{i}^{(l)}$ and $\mathbf{z}_{i}^{(l)}$ of cluster $\Omega_{i}^{(l)}$ that are to be eliminated next are highlighted in \autoref{fig:subgraph_tikz}(a). Performing this elimination results in new interactions between all neighbors of $\Omega_{i}^{(l)}$ (including itself), as is illustrated in \autoref{fig:subgraph_tikz}(b)\footnote{Eliminating a node also implies that the right hand side of the system needs to be updated.}. Some of these fill-ins establish direct connections between clusters that did not exist prior to this elimination step, however. For example, fill-ins $\mathbf{E}^{\prime}(\underline{\texttt{k}}, \underline{\texttt{j}})$ and $\mathbf{E}^{\prime}( \underline{\texttt{j}}, \underline{\texttt{k}} )$ between $\Omega_{j}^{(l)}$ and $\Omega_{k}^{(l)}$ arise (where $\mathbf{E}^{\prime}$ denotes a fill-in in the extended sparse matrix). As the aforementioned clusters are well-separated (i.e., their existing interactions are accounted for through the edges between $\overline{\texttt{k}}$ and $\overline{\texttt{j}}$ rather than between $\underline{\texttt{k}}$ and $\underline{\texttt{j}}$), these fill-ins are expected to be low-rank and can be approximated as follows:
\begin{align}
	 \mathbf{E}^{\prime}(\underline{\texttt{k}}, \underline{\texttt{j}}) & \simeq  \mathbf{U}^{\prime}( \underline{\texttt{k}}, \texttt{k} ) \, \mathbf{K}^{\prime}( \overline{\texttt{k}}, \overline{\texttt{j}} ) \,
	 \mathbf{V}^{\prime^{\mathrm{T}}} ( \texttt{j}, \underline{\texttt{j}} ) \label{eq:P2P_kj}\\
 	\mathbf{E}^{\prime}(\underline{\texttt{j}}, \underline{\texttt{k}}) & \simeq  \mathbf{U}^{\prime} ( \underline{\texttt{j}}, \texttt{j} ) \,  \mathbf{K}^{\prime}( \overline{\texttt{j}}, \overline{\texttt{k}} ) \, \mathbf{V}^{\prime^{\mathrm{T}}}( \texttt{k}, \underline{\texttt{k}} ) \label{eq:P2P_jk}
\end{align}
It is assumed in the following that \autoref{eq:P2P_kj} and \eqref{eq:P2P_jk} represent partial singular value decompositions (SVD), implying that the matrices $\mathbf{U}^{\prime}$ and $\mathbf{V}^{\prime}$ are orthonormal and that $\mathbf{K}^{\prime}$ is a diagonal matrix containing the singular values. The matrices  $\mathbf{K}^{\prime}( \overline{\texttt{k}}, \overline{\texttt{j}} )$ and  $\mathbf{K}^{\prime}( \overline{\texttt{j}}, \overline{\texttt{k}} )$ are therefore also denoted as $\mathbf{\Sigma}^{\prime}( \overline{\texttt{k}}, \overline{\texttt{j}} )$ and $\mathbf{\Sigma}^{\prime}( \overline{\texttt{j}}, \overline{\texttt{k}} )$, respectively. The implications for the algorithm if this is not the case will be discussed in \autoref{subsec:remarks}. It is also assumed that a criterion is available for determining an appropriate rank for approximations~\eqref{eq:P2P_kj} and~\eqref{eq:P2P_jk} (see \autoref{subsec:remarks} for more details).

These low-rank approximations for the fill-in enable the calculation of new inter- and anterpolation operators for $\Omega_{j}^{(l)}$ and $\Omega_{k}^{(l)}$, respectively. For example, a new interpolation operator $\widehat{\mathbf{U}}(\underline{\texttt{j}}, \texttt{j})$ for $\Omega_{j}^{(l)}$ is obtained by concatenating the existing scaled interpolation operator $\mathbf{U}(\underline{\texttt{j}}, \texttt{j}) \, \mathbf{\Sigma}_\mathbf{U}$ (where $\mathbf{\Sigma}_\mathbf{U}$ contains the weights associated with each column of $\mathbf{U}(\underline{\texttt{j}}, \texttt{j})$)\footnote{Initial weights $\mathbf{\Sigma}_\mathbf{U}$ for $\mathbf{U}(\underline{\texttt{j}}, \texttt{j})$ (i.e., upon initialization of the IFMM operators) can be obtained rigorously by applying an SVD to $\left[ \cdots \, | \mathbf{U}(\underline{\texttt{j}}, \texttt{j}) \,  \mathbf{K}( \overline{\texttt{j}}, \overline{\texttt{q}} ) \, \mathbf{V}^{\mathrm{T}}( \texttt{q}, \underline{\texttt{q}} )     \, | \cdots \, \right]$ ($\forall \Omega_{q}^{(l)} \in \mathcal{I}_{j}^{(l)}$, with $\mathcal{I}_{j}^{(l)}$ the interaction list of $\Omega_{j}^{(l)}$), or can be roughly estimated by sampling some matrix entries. The latter approach is less expensive in terms of computational cost.} and the scaled fill-in operator $\mathbf{U}^{\prime} ( \underline{\texttt{j}}, \texttt{j} ) \, \mathbf{\Sigma}^{\prime}( \overline{\texttt{j}}, \overline{\texttt{k}} )$ and subsequently applying a low-rank recompression (thin SVD):
\begin{align}
	\left[\mathbf{U}(\underline{\texttt{j}}, \texttt{j}) \, \mathbf{\Sigma}_\mathbf{U} \,  | \,  \mathbf{U}^{\prime} ( \underline{\texttt{j}}, \texttt{j} ) \, \mathbf{\Sigma}^{\prime}( \overline{\texttt{j}}, \overline{\texttt{k}} ) \right] & \simeq  \widehat{\mathbf{U}}(\underline{\texttt{j}}, \texttt{j}) \, \widehat{\mathbf{\Sigma}}_\mathbf{U} \, \mathbf{\Phi}_{j} \label{eq:U_j_recomp}
\end{align}
$\mathbf{\Phi}_{j} = \left[ {\boldsymbol{\phi}}_{j} \,  | \, \boldsymbol{\phi}^{\prime}_{j} \right]$ represents a fat, row orthogonal, matrix. It is important to note that the recompression in \autoref{eq:U_j_recomp} is applied to the scaled interpolation operators rather than to $\left[\mathbf{U}(\underline{\texttt{j}}, \texttt{j}) \,  | \,  \mathbf{U}^{\prime} ( \underline{\texttt{j}}, \texttt{j} ) \right]$, as this allows to account for the correct ``weight'' of each vector when obtaining a new basis.\footnote{The proof of this result, which is elementary, is omitted.} \autoref{eq:U_j_recomp} allows expressing the existing and fill-in operators $\mathbf{U}(\underline{\texttt{j}}, \texttt{j})$ and $\mathbf{U}^{\prime}(\underline{\texttt{j}}, \texttt{j})$ in terms of the new basis $\widehat{\mathbf{U}}(\underline{\texttt{j}}, \texttt{j})$:
\begin{align}
	\mathbf{U}(\underline{\texttt{j}}, \texttt{j}) & \simeq  \widehat{\mathbf{U}}(\underline{\texttt{j}}, \texttt{j}) \, \widehat{\mathbf{\Sigma}}_\mathbf{U} \, {\boldsymbol{\phi}}_{j} \, \mathbf{\Sigma}^{-1}_\mathbf{U} = \widehat{\mathbf{U}}(\underline{\texttt{j}}, \texttt{j}) \, \mathbf{r}_{j}\\
	\mathbf{U}^{\prime}(\underline{\texttt{j}}, \texttt{j}) & \simeq  \widehat{\mathbf{U}}(\underline{\texttt{j}}, \texttt{j}) \, \widehat{\mathbf{\Sigma}}_\mathbf{U} \, {\boldsymbol{\phi}}^{\prime}_{j} \, \mathbf{\Sigma}^{\prime}( \overline{\texttt{j}}, \overline{\texttt{k}} )^{-1}  = \widehat{\mathbf{U}}(\underline{\texttt{j}}, \texttt{j}) \, \mathbf{r}^{\prime}_{j}
\end{align}
The weights $\widehat{\Sigma}_\mathbf{U}$ of the new interpolation basis $\widehat{\mathbf{U}}(\underline{\texttt{j}}, \texttt{j})$ are stored, as they are needed in potential future recompressions. A new anterpolation operator $\widehat{\mathbf{V}}(\texttt{j},\underline{\texttt{j}})$ for $\Omega_{j}^{(l)}$ is obtained in a similar way:
\begin{align}
	\left[ \mathbf{V}(\texttt{j}, \underline{\texttt{j}})  \, \mathbf{\Sigma}_\mathbf{V} \, | \,  \mathbf{V}^{\prime}(\texttt{j}, \underline{\texttt{j}}) \, \mathbf{\Sigma}^{\prime}( \overline{\texttt{k}}, \overline{\texttt{j}} ) \right] & \simeq  \widehat{\mathbf{V}}(\texttt{j},\underline{\texttt{j}}) \, \widehat{\mathbf{\Sigma}}_\mathbf{V} \, \mathbf{\Psi}_{j} \label{eq:V_j_recomp}
\end{align}
with $\mathbf{\Psi}_{j} = \left[ \boldsymbol{\psi}_{j} \,  | \, \boldsymbol{\psi}^{\prime}_{j} \right]$, such that $\mathbf{V}(\texttt{j}, \underline{\texttt{j}})$ and $\mathbf{V}^{\prime}(\texttt{j}, \underline{\texttt{j}})$ can be written as:
\begin{align}
	\mathbf{V}(\underline{\texttt{j}}, \texttt{j}) & \simeq  \widehat{\mathbf{V}}(\underline{\texttt{j}}, \texttt{j}) \, \widehat{\mathbf{\Sigma}}_\mathbf{V} \, {\boldsymbol{\psi}}_{j} \, \mathbf{\Sigma}^{-1}_\mathbf{V} = \widehat{\mathbf{V}}(\underline{\texttt{j}}, \texttt{j}) \, \mathbf{t}_{j}\\
	\mathbf{V}^{\prime}(\underline{\texttt{j}}, \texttt{j}) & \simeq  \widehat{\mathbf{V}}(\underline{\texttt{j}}, \texttt{j}) \, \widehat{\mathbf{\Sigma}}_\mathbf{V} \, {\boldsymbol{\psi}}^{\prime}_{j} \, \mathbf{\Sigma}^{\prime}( \overline{\texttt{k}}, \overline{\texttt{j}} )^{-1}  = \widehat{\mathbf{V}}(\underline{\texttt{j}}, \texttt{j}) \, \mathbf{t}^{\prime}_{j}\label{eq:V_prime}
\end{align}
The same strategy as outlined in \autoref{eq:U_j_recomp} to \eqref{eq:V_prime} is applicable for obtaining new operators $\widehat{\mathbf{U}} ( \underline{\texttt{k}}, \texttt{k} )$ and $\widehat{\mathbf{V}}(\texttt{k},\underline{\texttt{k}})$ for $\Omega_{k}^{(l)}$.

The crucial benefit of compressing the dense fill-ins and updating the inter- and anterpolation operators of the corresponding clusters is that the former can now be accounted for through the operators $\mathbf{K}(\overline{\texttt{k}}, \overline{\texttt{j}})$ and $\mathbf{K}(\overline{\texttt{j}}, \overline{\texttt{k}})$ between $\Omega_{j}^{(l)}$ and $\Omega_{k}^{(l)}$. For example, the fill-in $\mathbf{E}^{\prime} ( \underline{\texttt{k}}, \underline{\texttt{j}} )$ is redirected through the $\mathbf{K}(\overline{\texttt{k}}, \overline{\texttt{j}})$-operator, which becomes:
\begin{align}
	\mathbf{K}(\overline{\texttt{k}}, \overline{\texttt{j}}) &\leftarrow 
	\mathbf{r}_{k} \, \mathbf{K}(\overline{\texttt{k}}, \overline{\texttt{j}}) \, \mathbf{t}^{\mathrm{T}}_{j} + \mathbf{r}^{\prime}_{k} \, \mathbf{K}^{\prime}(\overline{\texttt{k}},\overline{\texttt{j}}) \, \mathbf{t}^{\prime^{\mathrm{T}}}_{j} \label{eq:K_kj_update}
\end{align}
Updating the inter- and anterpolation operators also implies that all the other edges in the graph associated with $\Omega_{j}^{(l)}$ need to be adjusted. More specifically, the remaining operators $\mathbf{K}(\overline{\texttt{q}}, \overline{\texttt{j}})$ and $\mathbf{K}( \overline{\texttt{j}}, \overline{\texttt{q}})$ (with $\Omega_{q}^{(l)} \in \mathcal{N}_{j}^{(l)} \cup \mathcal{I}_{j}^{(l)} \backslash \Omega_{k}^{(l)} $) are replaced by:
\begin{align}
	\mathbf{K}(\overline{\texttt{q}}, \overline{\texttt{j}}) &\leftarrow \mathbf{K}(\overline{\texttt{q}}, \overline{\texttt{j}}) \, \mathbf{t}^{\mathrm{T}}_{j} \label{eq:K_qj_update}\\
	\mathbf{K}(\overline{\texttt{j}}, \overline{\texttt{q}}) &\leftarrow \mathbf{r}_{j} \, \mathbf{K}(\overline{\texttt{j}}, \overline{\texttt{q}}) \label{eq:K_jq_update}
\end{align}
where $\mathcal{N}_{j}^{(l)}$ and $\mathcal{I}_{j}^{(l)}$ represent the neighbor list (including itself) and interaction list (i.e., well-separated clusters that are children of its parent's neighbors) of $\Omega_{j}^{(l)}$, respectively. Furthermore, the inter- and anterpolation operators of the parent cluster of $\Omega_{j}^{(l)}$ are updated as follows:
\begin{align}
	\mathbf{U}( \overline{\texttt{j}},\texttt{j}^{\dagger} ) &\leftarrow \mathbf{r}_{j} \,  \mathbf{U}( \overline{\texttt{j}},\texttt{j}^{\dagger} ) \label{eq:U_parent_j_update}\\
	 \mathbf{V}( \texttt{j}^{\dagger}, \overline{\texttt{j}} ) &\leftarrow \mathbf{t}_{j}  \mathbf{V}( \texttt{j}^{\dagger}, \overline{\texttt{j}} )  \label{eq:V_parent_j_update}
\end{align}
Similar updates as in \autoref{eq:K_kj_update}--\eqref{eq:V_parent_j_update} are required for all edges associated with the cluster $\Omega_{k}^{(l)}$. Edges that need to be updated are highlighted in yellow in \autoref{fig:subgraph_tikz}(c). 

The technique outlined above leads to a significant reduction of the amount of fill-in and ensures that the sparsity pattern of the extended sparse matrix is preserved throughout the elimination phase. This leads to an algorithm with linear complexity $\mO(N)$. This assumes that the rank remains bounded independent of the matrix size. This is expected to be the case, as only interactions between well-separated clusters are replaced by low-rank approximations.

This strategy of compressing and redirecting dense fill-ins is very similar to the initial sparsification approach introduced in \autoref{subsec:extended_sparsification}, as low-rank interactions are represented at the parent level. Apart from fill-ins between the $\tilde{\mathbf{x}}$-variables (e.g., between $\tilde{\mathbf{x}}_{j}^{(l)}$ and $\tilde{\mathbf{x}}_{k}^{(l)}$), fill-ins between the variables $\tilde{\mathbf{x}}^{(l)}$ and $\mathbf{y}^{(l)}$ are created as well throughout the elimination phase. If these occur between well-separated clusters (for example $\mathbf{E}^{\prime}(\overline{\texttt{h}},\underline{\texttt{j}})$, $\mathbf{E}^{\prime}(\underline{\texttt{j}}, \overline{\texttt{h}})$ and $\mathbf{E}^{\prime}( \overline{\texttt{h}},\underline{\texttt{k}})$, $\mathbf{E}^{\prime}(\underline{\texttt{k}}, \overline{\texttt{h}})$ in \autoref{fig:subgraph_tikz}(b)), the fill-ins are compressed and redirected through other operators at the parent level.

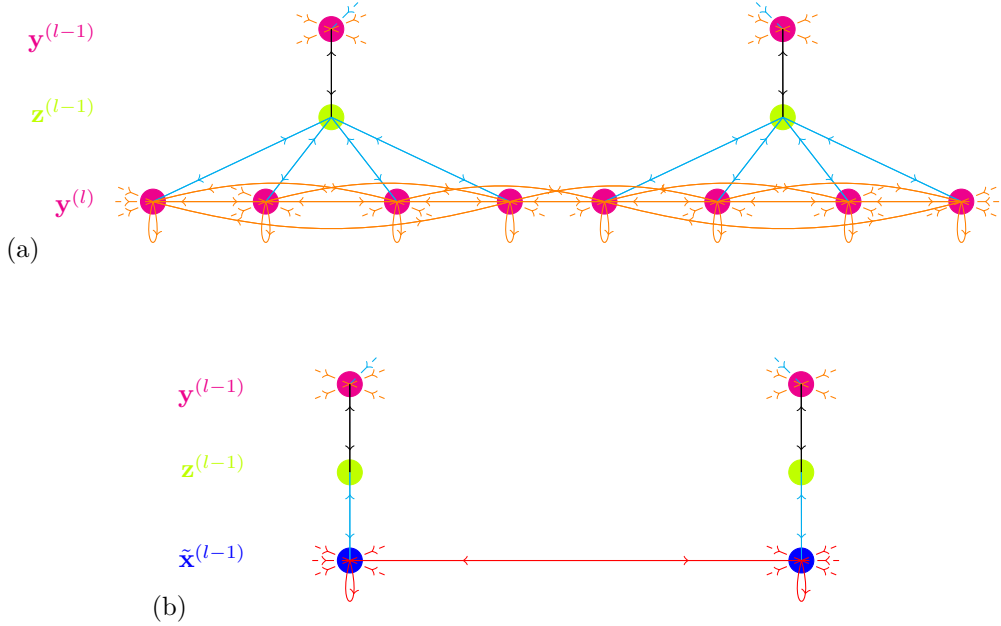
\begin{figure}[hbtp]
	\begin{center}	
		(a)~\hspace{-1em}\begin{tikzpicture}[scale=2.5]

		\node (0) at (0,0) {};

		\node[draw=magenta,fill=magenta,circle,minimum size=2mm,right=15mm of 0] (1) {};	
		\node[draw=magenta,fill=magenta,circle,minimum size=2mm,right=75mm of 0] (2) {};	
		
		\draw[cyan,-<-=.75,dashed] (1.center) -- ++(45:2mm);
		\draw[cyan,-<-=.75,dashed] (2.center) -- ++(135:2mm);
		
		\draw[orange,-<-=.75,dashed] (1.center) -- ++(22.5:2mm);
		\draw[orange,-<-=.75,dashed] (1.center) -- ++(157.5:2mm);
		\draw[orange,-<-=.75,dashed] (1.center) -- ++(202.5:2mm);
		\draw[orange,-<-=.75,dashed] (1.center) -- ++(337.5:2mm);
		
		\draw[orange,-<-=.75,dashed] (2.center) -- ++(22.5:2mm);
		\draw[orange,-<-=.75,dashed] (2.center) -- ++(157.5:2mm);
		\draw[orange,-<-=.75,dashed] (2.center) -- ++(202.5:2mm);
		\draw[orange,-<-=.75,dashed] (2.center) -- ++(337.5:2mm);
				
		\node[draw=lime,fill=lime,circle,minimum size=2mm,below=10mm of 1.center] (3) {};	
		\node[draw=lime,fill=lime,circle,minimum size=2mm,below=10mm of 2.center] (4) {};	
		
		\draw[black,->-=.75] (1.center) -- (3.center);
		\draw[black,->-=.75] (2.center) -- (4.center);

		\draw[black,->-=.75] (3.center) -- (1.center);
		\draw[black,->-=.75] (4.center) -- (2.center);
				
%
%
%
%
%
%

		\node[draw=magenta,fill=magenta,circle,minimum size=2mm,below left=10mm and 7.5mm of 3.center] (6) {};	
		\node[draw=magenta,fill=magenta,circle,minimum size=2mm,below right=10mm and 7.5mm of 3.center] (7) {};	
		\node[draw=magenta,fill=magenta,circle,minimum size=2mm,below left=10mm and 22.5mm of 3.center] (5) {};	
		\node[draw=magenta,fill=magenta,circle,minimum size=2mm,below right=10mm and 22.5mm of 3.center] (8) {};	

		\node[draw=magenta,fill=magenta,circle,minimum size=2mm,below left=10mm and 7.5mm of 4.center] (10) {};	
		\node[draw=magenta,fill=magenta,circle,minimum size=2mm,below right=10mm and 7.5mm of 4.center] (11) {};	
		\node[draw=magenta,fill=magenta,circle,minimum size=2mm,below left=10mm and 22.5mm of 4.center] (9) {};	
		\node[draw=magenta,fill=magenta,circle,minimum size=2mm,below right=10mm and 22.5mm of 4.center] (12) {};

		\draw[cyan,->-=.75] (5.center) -- (3.center);
		\draw[cyan,->-=.75] (6.center) -- (3.center);
		\draw[cyan,->-=.75] (7.center) -- (3.center);
		\draw[cyan,->-=.75] (8.center) -- (3.center);

		\draw[cyan,->-=.75] (9.center) -- (4.center);
		\draw[cyan,->-=.75] (10.center) -- (4.center);
		\draw[cyan,->-=.75] (11.center) -- (4.center);
		\draw[cyan,->-=.75] (12.center) -- (4.center);

		\draw[cyan,->-=.75] (3.center) -- (5.center);
		\draw[cyan,->-=.75] (3.center) -- (6.center);
		\draw[cyan,->-=.75] (3.center) -- (7.center);
		\draw[cyan,->-=.75] (3.center) -- (8.center);

		\draw[cyan,->-=.75] (4.center) -- (9.center);
		\draw[cyan,->-=.75] (4.center) -- (10.center);
		\draw[cyan,->-=.75] (4.center) -- (11.center);
		\draw[cyan,->-=.75] (4.center) -- (12.center);

		\draw [orange,->-=.4,loop below,looseness=10,min distance=3mm] (5.center) to (5.center);
		\draw [orange,->-=.4,loop below,looseness=10,min distance=3mm] (6.center) to (6.center);
		\draw [orange,->-=.4,loop below,looseness=10,min distance=3mm] (7.center) to (7.center);
		\draw [orange,->-=.4,loop below,looseness=10,min distance=3mm] (8.center) to (8.center);
		
		\draw [orange,->-=.4,loop below,looseness=10,min distance=3mm] (9.center) to (9.center);
		\draw [orange,->-=.4,loop below,looseness=10,min distance=3mm] (10.center) to (10.center);
		\draw [orange,->-=.4,loop below,looseness=10,min distance=3mm] (11.center) to (11.center);
		\draw [orange,->-=.4,loop below,looseness=10,min distance=3mm] (12.center) to (12.center);	

		\draw[orange,->-=.75] (5.center) -- (6.center);
		\draw[orange,->-=.75] (6.center) -- (7.center);
		\draw[orange,->-=.75] (7.center) -- (8.center);
		\draw[orange,->-=.75] (8.center) -- (9.center);
		\draw[orange,->-=.75] (9.center) -- (10.center);
		\draw[orange,->-=.75] (10.center) -- (11.center);
		\draw[orange,->-=.75] (11.center) -- (12.center);

		\draw[orange,->-=.75] (6.center) -- (5.center);
		\draw[orange,->-=.75] (7.center) -- (6.center);
		\draw[orange,->-=.75] (8.center) -- (7.center);
		\draw[orange,->-=.75] (9.center) -- (8.center);
		\draw[orange,->-=.75] (10.center) -- (9.center);
		\draw[orange,->-=.75] (11.center) -- (10.center);
		\draw[orange,->-=.75] (12.center) -- (11.center);

		\draw[orange,-<-=.75,dashed] (5.center) -- ++(157.5:2mm);
		\draw[orange,-<-=.75,dashed] (5.center) -- ++(180:2mm);
		\draw[orange,-<-=.75,dashed] (5.center) -- ++(202.5:2mm);
		\draw[orange,-<-=.75,dashed] (6.center) -- ++(157.5:2mm);
		\draw[orange,-<-=.75,dashed] (6.center) -- ++(202.5:2mm);
		\draw[orange,-<-=.75,dashed] (6.center) -- ++(337.5:2mm);
		\draw[orange,-<-=.75,dashed] (7.center) -- ++(202.5:2mm);
		\draw[orange,-<-=.75,dashed] (7.center) -- ++(337.5:2mm);
		\draw[orange,-<-=.75,dashed] (8.center) -- ++(337.5:2mm);
		\draw[orange,-<-=.75,dashed] (9.center) -- ++(202.5:2mm);
		\draw[orange,-<-=.75,dashed] (10.center) -- ++(202.5:2mm);
		\draw[orange,-<-=.75,dashed] (10.center) -- ++(337.5:2mm);
		\draw[orange,-<-=.75,dashed] (11.center) -- ++(22.5:2mm);
		\draw[orange,-<-=.75,dashed] (11.center) -- ++(202.5:2mm);
		\draw[orange,-<-=.75,dashed] (11.center) -- ++(337.5:2mm);
		\draw[orange,-<-=.75,dashed] (12.center) -- ++(0:2mm);	
		\draw[orange,-<-=.75,dashed] (12.center) -- ++(22.5:2mm);
		\draw[orange,-<-=.75,dashed] (12.center) -- ++(202.5:2mm);
		\draw[orange,-<-=.75,dashed] (12.center) -- ++(337.5:2mm);

		\draw[orange,bend left=15,->-=.75]  (5.center) to (7.center);
		\draw[orange,bend left=15,-<-=.25]  (5.center) to (7.center);
		\draw[orange,bend left=15,->-=.75]  (6.center) to (8.center);
		\draw[orange,bend left=15,-<-=.25]  (6.center) to (8.center);
		\draw[orange,bend left=15,->-=.75]  (7.center) to (9.center);
		\draw[orange,bend left=15,-<-=.25]  (7.center) to (9.center);
		\draw[orange,bend left=15,->-=.75]  (8.center) to (10.center);
		\draw[orange,bend left=15,-<-=.25]  (8.center) to (10.center);
		\draw[orange,bend left=15,->-=.75]  (9.center) to (11.center);
		\draw[orange,bend left=15,-<-=.25]  (9.center) to (11.center);
		\draw[orange,bend left=15,->-=.75]  (10.center) to (12.center);
		\draw[orange,bend left=15,-<-=.25]  (10.center) to (12.center);
		
		\draw[orange,bend right=15,->-=.75]  (5.center) to (8.center);
		\draw[orange,bend right=15,-<-=.25]  (5.center) to (8.center);
		\draw[orange,bend right=15,->-=.75]  (9.center) to (12.center);
		\draw[orange,bend right=15,-<-=.25]  (9.center) to (12.center);

		\node[below left=-2mm and 30mm of 1.center] () {\textcolor{magenta}{$\mathbf{y}^{(l-1)}$}};
		\node[below left=8mm and 30mm of 1.center] () {\textcolor{lime}{$\mathbf{z}^{(l-1)}$}};
		\node[below left=20mm and 30mm of 1.center] () {\textcolor{magenta}{$\mathbf{y}^{(l)}$}};
%
%
%
%
%
	
	\end{tikzpicture}\vspace{2.5em}\\
		(b)~\hspace{-1em}\begin{tikzpicture}[scale=2.5]

		\node (0) at (0,0) {};

		\node[draw=magenta,fill=magenta,circle,minimum size=2mm,right=15mm of 0] (1) {};	
		\node[draw=magenta,fill=magenta,circle,minimum size=2mm,right=75mm of 0] (2) {};	
		
		\draw[cyan,-<-=.75,dashed] (1.center) -- ++(45:2mm);
		\draw[cyan,-<-=.75,dashed] (2.center) -- ++(135:2mm);
		
		\draw[orange,-<-=.75,dashed] (1.center) -- ++(22.5:2mm);
		\draw[orange,-<-=.75,dashed] (1.center) -- ++(157.5:2mm);
		\draw[orange,-<-=.75,dashed] (1.center) -- ++(202.5:2mm);
		\draw[orange,-<-=.75,dashed] (1.center) -- ++(337.5:2mm);
		
		\draw[orange,-<-=.75,dashed] (2.center) -- ++(22.5:2mm);
		\draw[orange,-<-=.75,dashed] (2.center) -- ++(157.5:2mm);
		\draw[orange,-<-=.75,dashed] (2.center) -- ++(202.5:2mm);
		\draw[orange,-<-=.75,dashed] (2.center) -- ++(337.5:2mm);
				
		\node[draw=lime,fill=lime,circle,minimum size=2mm,below=10mm of 1.center] (3) {};	
		\node[draw=lime,fill=lime,circle,minimum size=2mm,below=10mm of 2.center] (4) {};	
		
		\draw[black,->-=.75] (1.center) -- (3.center);
		\draw[black,->-=.75] (3.center) -- (1.center);
		\draw[black,->-=.75] (2.center) -- (4.center);
		\draw[black,->-=.75] (4.center) -- (2.center);

		\node[draw=blue,fill=blue,circle,minimum size=2mm,below=10mm of 3.center] (5) {};
		\node[draw=blue,fill=blue,circle,minimum size=2mm,below=10mm of 4.center] (6) {};


		\draw[cyan,->-=.75] (5.center) -- (3.center);
		\draw[cyan,->-=.75] (6.center) -- (4.center);

		\draw[cyan,->-=.75] (3.center) -- (5.center);
		\draw[cyan,->-=.75] (4.center) -- (6.center);

		\draw [red,->-=.4,loop below,looseness=10,min distance=3mm] (5.center) to (5.center);
		\draw [red,->-=.4,loop below,looseness=10,min distance=3mm] (6.center) to (6.center);

		\draw[red,->-=.75] (5.center) -- (6.center);
		\draw[red,->-=.75] (6.center) -- (5.center);

		\draw[red,-<-=.75,dashed] (5.center) -- ++(22.5:2mm);
		\draw[red,-<-=.75,dashed] (5.center) -- ++(180:2mm);
		\draw[red,-<-=.75,dashed] (5.center) -- ++(157.5:2mm);
		\draw[red,-<-=.75,dashed] (5.center) -- ++(202.5:2mm);
		\draw[red,-<-=.75,dashed] (5.center) -- ++(337.5:2mm);
	
		\draw[red,-<-=.75,dashed] (6.center) -- ++(0:2mm);
		\draw[red,-<-=.75,dashed] (6.center) -- ++(22.5:2mm);
		\draw[red,-<-=.75,dashed] (6.center) -- ++(157.5:2mm);
		\draw[red,-<-=.75,dashed] (6.center) -- ++(202.5:2mm);
		\draw[red,-<-=.75,dashed] (6.center) -- ++(337.5:2mm);

		\node[below left=-2mm and 12.5mm of 1.center] () {\textcolor{magenta}{$\mathbf{y}^{(l-1)}$}};
		\node[below left=8mm and 12.5mm of 1.center] () {\textcolor{lime}{$\mathbf{z}^{(l-1)}$}};
		\node[below left=20mm and 12.5mm of 1.center] () {\textcolor{blue}{$\tilde{\mathbf{x}}^{(l-1)}$}};
%
%
%
%
%
	
	\end{tikzpicture}	
	\end{center}
	\caption{Subgraph (a)~before and (b)~after merging the variables $\mathbf{y}^{(l)}$ of level $l$ into unknowns $\tilde{\mathbf{x}}^{(l-1)}$ of the parent level $l-1$. The associated operators are combined accordingly.}
	\label{fig:transfer_tikz}
\end{figure}

Once all unknowns $\tilde{\mathbf{x}}^{(l)}$ and variables $\mathbf{z}^{(l)}$ at level $l$ have been eliminated, the next level is considered. From a conceptual point of view, variables $\mathbf{y}^{(l)}$ can be interpreted as unknowns $\tilde{\mathbf{x}}^{(l-1)}$ of the parent level $l-1$. More specifically, a set of variables $\mathbf{y}_{i}^{(l)}$ is combined to form $\tilde{\mathbf{x}}_{\mathrm{p}}^{(l-1)}$, where the clusters $\Omega_{i}^{(l)}$ are children of the parent cluster $\Omega_{\mathrm{p}}^{(l-1)}$:
\begin{align}
	\tilde{\mathbf{x}}_{\mathrm{p}}^{(l-1)} &= \begin{Bmatrix} \cdots \, \mathbf{y}_{i}^{(l)} \, \cdots  \end{Bmatrix}^{\mathrm{T}} \qquad \forall \, \Omega_{i}^{(l)} \in \mathcal{C}_{\mathrm{p}}^{(l-1)} 
	\label{eq:y_merge}
\end{align}
or, equivalently, the nodes $\overline{\texttt{i}}^{(l)}$ are joined to form $\underline{\texttt{p}}^{(l-1)}$:
\begin{align}
	\underline{\texttt{p}}^{(l-1)} &= \bigcup_{\Omega_{i}^{(l)} \in \mathcal{C}_{\mathrm{p}}^{(l-1)}} \overline{\texttt{i}}^{(l)}
	\label{eq:y_merge_2}
\end{align}
$\mathcal{C}_{\mathrm{p}}^{(l-1)}$ in \autoref{eq:y_merge} and~\eqref{eq:y_merge_2} denotes the children list of $\Omega_{\mathrm{p}}^{(l-1)}$. This merging procedure is illustrated in \autoref{fig:transfer_tikz}. Furthermore, the corresponding operators are combined accordingly to form the inter- and anterpolation operators $\mathbf{U} ( \underline{\texttt{p}}^{(l-1)}, \texttt{p}^{(l-1)} )$ and $\mathbf{V} ( \texttt{p}^{(l-1)}, \underline{\texttt{p}}^{(l-1)} )$ of $\Omega_{\mathrm{p}}^{(l-1)}$:
\begin{align}
\mathbf{U} ( \underline{\texttt{p}}^{(l-1)}, \texttt{p}^{(l-1)} ) &= 
\begin{bmatrix} \vdots \\ \mathbf{U}(\overline{\texttt{i}}^{(l)}, 
{\texttt{i}^{\dagger}}^{(l)}) \\ \vdots \end{bmatrix} 
\qquad \forall \, \Omega_{i}^{(l)} \in \mathcal{C}_{\mathrm{p}}^{(l-1)} \label{eq:U_merge}\\
\mathbf{V} ( \texttt{p}^{(l-1)}, \underline{\texttt{p}}^{(l-1)} ) &=
\begin{bmatrix} \vdots \\ \mathbf{V}({\texttt{i}^{\dagger}}^{(l)}), \overline{\texttt{i}}^{(l)}) \\ \vdots \end{bmatrix} 
\qquad \forall \, \Omega_{i}^{(l)} \in \mathcal{C}_{\mathrm{p}}^{(l-1)} \label{eq:V_merge}
\end{align}
A similar operation is performed to merge the interactions between variables $\mathbf{y}^{(l)}$ into interactions between unknowns $\tilde{\mathbf{x}}^{(l-1)}$. For example, the interaction between $\tilde{\mathbf{x}}_{\mathrm{p}}^{(l-1)}$ and $\tilde{\mathbf{x}}_{\mathrm{q}}^{(l-1)}$ of clusters $\Omega_{\mathrm{p}}^{(l-1)}$ and $\Omega_{\mathrm{q}}^{(l-1)}$ is assembled as follows:
\begin{align}
	\mathbf{K}( \underline{\texttt{p}}^{(l-1)}, \underline{\texttt{q}}^{(l-1)}) &= \begin{bmatrix} & \vdots & \\ \cdots & \mathbf{K}( \overline{\texttt{i}}^{(l)}, \overline{\texttt{j}}^{(l)} ) &  \cdots \\ & \vdots & \end{bmatrix} \qquad \forall \, \Omega_{i}^{(l)} \in \mathcal{C}_{\mathrm{p}}^{(l-1)}, \,  \forall \, \Omega_{j}^{(l)} \in \mathcal{C}_{\mathrm{q}}^{(l-1)} \label{eq:K_merge}
\end{align}
The resulting graph has exactly the same structure as the initial one, except for the fact that it consists of one level less. The elimination and merging procedure are repeated consecutively for each level $l \geq 2$ until only a small dense system for $\mathbf{y}^{(2)}$ (or, equivalently, $\tilde{\mathbf{x}}^{(1)}$) remains to be solved. Once this is done, the elimination phase is terminated and all the other variables can be retrieved through substitution, marching down the graph in the reverse order. In the end, the solution $\mathbf{x}$ is obtained.

%

\subsection{The algorithm}\label{subsec:algorithm}
Based on the key ideas introduced in \autoref{subsec:extended_sparsification} and \ref{subsec:compression}, the general IFMM algorithm can now be presented in a formal way. In the following, it is assumed that the $\mH^2$-structure of the matrix $\mathbf{A}$ is based on a three-dimensional octree decomposition, implying that every cluster $\Omega_{i}^{(l)}$ can have up to $3^3 = 27$ clusters in its neighbor list $\mathcal{N}_{i}^{(l)}$ (including itself) and up to $6^3-3^3 = 189$ clusters in its interaction list $\mathcal{I}_{i}^{(l)}$.

\autoref{alg:IFMM_general} presents the general outline of the IFMM. It consists of initializing the IFMM-operators, performing an upward pass to eliminate each level~$l$ (and forming the appropriate operators of the parent level~$l-1$), solving for $\mathbf{y}^{(2)} (= \tilde{\mathbf{x}}^{(1)})$ at the top level, and substituting through a downward pass to finally retrieve the solution~$\mathbf{x}$. The elimination procedure at each level $l$ is detailed in \autoref{alg:IFMM_elimination}. Eliminating a cluster $\Omega_{i}^{(l)}$ leads to new interactions between its neighbors $\Omega_{j}^{(l)}$ and $\Omega_{k}^{(l)}$ (with $\Omega_{j}^{(l)},\Omega_{k}^{(l)} \in \mathcal{N}_{i}^{(l)}$): $\mathbf{E}^{\prime}(\overline{\texttt{k}}, \overline{\texttt{j}} )$ and $\mathbf{E}^{\prime} (\overline{\texttt{j}}, \overline{\texttt{k}})$ if both $\Omega_{j}^{(l)}$ and $\Omega_{k}^{(l)}$ have already been eliminated, $\mathbf{E}^{\prime}( \underline{\texttt{k}}, \overline{\texttt{j}} )$ and $\mathbf{E}^{\prime}(\overline{\texttt{j}}, \underline{\texttt{k}})$ if $\Omega_{j}^{(l)}$ has been eliminated but $\Omega_{k}^{(l)}$ has not (and vice versa), and $\mathbf{E}^{\prime}(\underline{\texttt{k}}, \underline{\texttt{j}})$ and $\mathbf{E}^{\prime} (\underline{\texttt{j}}, \underline{\texttt{k}})$ if neither $\Omega_{j}^{(l)}$ nor $\Omega_{k}^{(l)}$ has been eliminated yet. If $\Omega_{j}^{(l)}$ and $\Omega_{k}^{(l)}$ are well-separated, these fill-ins are compressed into a low-rank approximation and subsequently redirected through other operators, as discussed in \autoref{subsec:compression}. This redirection procedure is described in \autoref{alg:IFMM_update}. If this is not the case, the additional interactions are added to the existing ones. The substitution phase proceeds level per level in the reverse elimination order and allows determining $\mathbf{z}^{(l)}$ and $\tilde{\mathbf{x}}^{(l)}$ based on $\mathbf{y}^{(l)}$. If $l$ is not the leaf level, the variables $\tilde{\mathbf{x}}^{(l)}$ are nothing else than the variables $\mathbf{y}^{(l+1)}$ of the child level, cf., \autoref{eq:y_merge}. At the leaf level, $\tilde{\mathbf{x}}^{(L)}$ corresponds to the desired solution~$\mathbf{x}$. As the elimination phase is decoupled from the substitution phase (as in a standard direct solve), the algorithm is well suited to handle multiple right hand sides $\mathbf{b}$.

The algorithm, as described in this Subsection, is for a general matrix $\mathbf{A}$; it is not specialized for the symmetric and positive definite case. It is possible to take advantage of the symmetry and positive-definiteness of a matrix to reduce the computational cost with almost a factor of two (similar to the Cholesky factorization that is almost twice as efficient as the LU factorization). This does not alter the fundamental ideas of the proposed algorithm.

\begin{algorithm}
	\begin{algorithmic}[1]
	\STATE Construct an octree decomposition of the domain $\Omega$
	\FOR{$l=L,L-1,\ldots,2$}
		\STATE Initialize the IFMM operators at level $l$
	\ENDFOR

	\FOR{$l=L,L-1,\ldots,2$}
		\STATE Eliminate level $l$ \textit{(\autoref{alg:IFMM_elimination})} and update the right hand side
		\IF{$l>2$}
			\STATE Form the appropriate operators of the parent level $l-1$
		\ENDIF
	\ENDFOR
	
	\STATE Solve for $\mathbf{y}^{(2)} (= \tilde{\mathbf{x}}^{(1)})$
	\FOR{$l=2,\ldots,L-1,L$}
		\STATE Substitute at level $l$ for $\mathbf{z}^{(l)}$ and $\tilde{\mathbf{x}}^{(l)}$
		\IF{$l<L$}
			\STATE	Form $\mathbf{y}^{(l+1)}$ based on  $\tilde{\mathbf{x}}^{(l)}$ \textit{(\autoref{eq:y_merge})}
		\ELSE
			\STATE	Set $\mathbf{x} = \tilde{\mathbf{x}}^{(L)}$
		\ENDIF
	\ENDFOR
	\end{algorithmic}
	\caption{The inverse fast multipole method: overall algorithm}
	\label{alg:IFMM_general}
\end{algorithm}
\begin{algorithm}
	\algsetup{indent=0.55em}
	\begin{algorithmic}[1] 
		\FORALL{clusters $\Omega_{i}^{(l)}$ at level $l$}
			\STATE Eliminate $\underline{\texttt{i}}$ and $\texttt{i}$ of $\Omega_{i}^{(l)}$
			\FORALL{unique pairs $\{\Omega_{j}^{(l)}, \Omega_{k}^{(l)}\}$ with $\Omega_{j}^{(l)},\Omega_{k}^{(l)} \in \mathcal{N}_{i}^{(l)}$}
					\IF{$\underline{\texttt{j}}$ and $\texttt{j}$ of $\Omega_{j}^{(l)}$ have been eliminated}
						\IF{$\underline{\texttt{k}}$ and $\texttt{k}$ of $\Omega_{k}^{(l)}$ have been eliminated}
							\STATE Compute the fill-ins $\mathbf{E}^{\prime}(\overline{\texttt{k}}, \overline{\texttt{j}})$ and $\mathbf{E}^{\prime}(\overline{\texttt{j}}, \overline{\texttt{k}})$
							\STATE Add the fill-ins: $\mathbf{E}(\overline{\texttt{k}}, \overline{\texttt{j}}) \leftarrow \mathbf{E}(\overline{\texttt{k}}, \overline{\texttt{j}}) + \mathbf{E}^{\prime}(\overline{\texttt{k}}, \overline{\texttt{j}})$ and $\mathbf{E}(\overline{\texttt{j}}, \overline{\texttt{k}}) \leftarrow \mathbf{E}(\overline{\texttt{j}}, \overline{\texttt{k}}) + \mathbf{E}^{\prime}(\overline{\texttt{j}}, \overline{\texttt{k}})$ 
						\ELSE
							\STATE Compute the fill-ins $\mathbf{E}^{\prime} ( \underline{\texttt{k}}, \overline{\texttt{j}})$ and $\mathbf{E}^{\prime} (\overline{\texttt{j}}, \underline{\texttt{k}})$	
							\IF{$\Omega_{j}^{(l)} \not \in \mathcal{N}_{k}^{(l)}$}
								\STATE Compress $\mathbf{E}^{\prime} (\underline{\texttt{k}}, \overline{\texttt{j}}) \simeq  \mathbf{U}^{\prime}(\underline{\texttt{k}}, \texttt{k}) \, \mathbf{K}^{\prime} (\overline{\texttt{k}}, \overline{\texttt{j}})$ and $\mathbf{E}^{\prime}(\overline{\texttt{j}}, \underline{\texttt{k}}) \simeq \mathbf{K}^{\prime} (\overline{\texttt{j}}, \overline{\texttt{k}}) \, \mathbf{V}^{\prime^{\mathrm{T}}}( \texttt{k}, \underline{\texttt{k}}) $
								\STATE Update the edges associated with $\Omega_{k}^{(l)}$ \textit{(\autoref{alg:IFMM_update})}
							\ELSE
								\STATE Add the fill-ins: $\mathbf{E} ( \underline{\texttt{k}}, \overline{\texttt{j}}) \leftarrow \mathbf{E} ( \underline{\texttt{k}}, \overline{\texttt{j}})+ \mathbf{E}^{\prime} ( \underline{\texttt{k}}, \overline{\texttt{j}})$ and $\mathbf{E}(\overline{\texttt{j}}, \underline{\texttt{k}}) \leftarrow \mathbf{E}(\overline{\texttt{j}}, \underline{\texttt{k}}) + \mathbf{E}^{\prime} (\overline{\texttt{j}}, \underline{\texttt{k}})$	
							\ENDIF
						\ENDIF	
					\ELSE
						\IF{$\underline{\texttt{k}}$ and $\texttt{k}$ of $\Omega_{k}^{(l)}$ have been eliminated}
							\STATE Compute the fill-ins: $\mathbf{E}^{\prime}(\overline{\texttt{k}}, \underline{\texttt{j}})$ and $\mathbf{E}^{\prime}(\underline{\texttt{j}}, \overline{\texttt{k}})$	
							\IF{$\Omega_{j}^{(l)} \not \in \mathcal{N}_{k}^{(l)}$}
								\STATE Compress $\mathbf{E}^{\prime}(\overline{\texttt{k}}, \underline{\texttt{j}} ) \simeq \mathbf{K}^{\prime}( \overline{\texttt{k}}, \overline{\texttt{j}}) \, \mathbf{V}^{\prime^{\mathrm{T}}}(\texttt{j},\underline{\texttt{j}}) $ and $\mathbf{E}^{\prime}(\underline{\texttt{j}}, \overline{\texttt{k}}) \simeq \mathbf{U}^{\prime}(\underline{\texttt{j}}, \texttt{j}) \, \mathbf{K}^{\prime} (\overline{\texttt{j}}, \overline{\texttt{k}})$
								\STATE Update the edges associated with $\Omega_{j}^{(l)}$ \textit{(\autoref{alg:IFMM_update})}
							\ELSE
								\STATE Add the fill-ins: $\mathbf{E}( \overline{\texttt{k}}, \underline{\texttt{j}}) \leftarrow \mathbf{E} ( \overline{\texttt{k}}, \underline{\texttt{j}})+ \mathbf{E}^{\prime} ( \overline{\texttt{k}}, \underline{\texttt{j}})$ and $\mathbf{E}(\underline{\texttt{j}}, \overline{\texttt{k}}) \leftarrow \mathbf{E}(\underline{\texttt{j}}, \overline{\texttt{k}}) + \mathbf{E}^{\prime} (\underline{\texttt{j}}, \overline{\texttt{k}})$
							\ENDIF
						\ELSE
							\STATE Compute the fill-ins $\mathbf{E}^{\prime}( \underline{\texttt{k}}, \underline{\texttt{j}})$ and $\mathbf{E}^{\prime}(\underline{\texttt{j}}, \underline{\texttt{k}})$	
							\IF{$\Omega_{j}^{(l)} \not \in \mathcal{N}_{k}^{(l)}$}
								\STATE Compress $\mathbf{E}^{\prime}(\underline{\texttt{k}}, \underline{\texttt{j}})  \simeq  \mathbf{U}^{\prime}( \underline{\texttt{k}}, \texttt{k}) \, \mathbf{K}^{\prime}( \overline{\texttt{k}}, \overline{\texttt{j}} ) \, \mathbf{V}^{\prime^{\mathrm{T}}} ( \texttt{j}, \underline{\texttt{j}} )$ and  $\mathbf{E}^{\prime}(\underline{\texttt{j}}, \underline{\texttt{k}}) \simeq  \mathbf{U}^{\prime} ( \underline{\texttt{j}}, \texttt{j} ) \,  \mathbf{K}^{\prime}( \overline{\texttt{j}}, \overline{\texttt{k}} ) \, \mathbf{V}^{\prime^{\mathrm{T}}}( \texttt{k}, \underline{\texttt{k}})$
								\STATE Update the edges associated with $\Omega_{j}^{(l)}$ \textit{(\autoref{alg:IFMM_update})}
								\STATE Update the edges associated with $\Omega_{k}^{(l)}$ \textit{(\autoref{alg:IFMM_update})}
							\ELSE
								\STATE Add the fill-ins: $\mathbf{E}( \underline{\texttt{k}}, \underline{\texttt{j}}) \leftarrow \mathbf{E}( \underline{\texttt{k}}, \underline{\texttt{j}}) + \mathbf{E}^{\prime}( \underline{\texttt{k}}, \underline{\texttt{j}})$ and $\mathbf{E}( \underline{\texttt{j}}, \underline{\texttt{k}}) \leftarrow \mathbf{E}( \underline{\texttt{j}}, \underline{\texttt{k}}) + \mathbf{E}^{\prime}( \underline{\texttt{j}}, \underline{\texttt{k}})$
							\ENDIF
						\ENDIF
					\ENDIF
			\ENDFOR
		\ENDFOR
	\end{algorithmic}
	\caption{The inverse fast multipole method: elimination of level $l$. In order to simplify the notations as much as possible, we made the following choices. The letter $\mathbf{E}$ denotes the sparse matrix on which we are doing the Gaussian elimination. The fill-in introduced is denoted by $\mathbf{E}^{\prime}$. This fill-in is progressively compressed (if well-separated) and folded into existing edges (non-zero entries) of $\mathbf{E}$. This compression process is done using matrices denoted by $\mathbf{U}^{\prime}$, $\mathbf{V}^{\prime}$, and $\mathbf{K}^{\prime}$. $\mathbf{U}^{\prime}$ is an interpolation operator ($\texttt{i}$ to $\underline{\texttt{i}}$), while $\mathbf{V}^{\prime}$ is an anterpolation operator ($\underline{\texttt{i}}$ to $\texttt{i}$). The entries denoted by $\mathbf{K}^{\prime}$ correspond to edges between nodes at the same level.}
	\label{alg:IFMM_elimination}
\end{algorithm}
\begin{algorithm}
	\begin{algorithmic}[1]
	\STATE $\mathbf{U}(\underline{\texttt{j}}, \texttt{j}) \leftarrow \widehat{\mathbf{U}}(\underline{\texttt{j}}, \texttt{j})$, $\mathbf{\Sigma}_\mathbf{U} \leftarrow \widehat{\mathbf{\Sigma}}_\mathbf{U}$
	\STATE $\mathbf{V}(\texttt{j},\underline{\texttt{j}}) \leftarrow \widehat{\mathbf{V}}(\texttt{j}, \underline{\texttt{j}})$, $\mathbf{\Sigma}_\mathbf{V} \leftarrow \widehat{\mathbf{\Sigma}}_\mathbf{V}$
	\STATE $\mathbf{U}( \overline{\texttt{j}},\texttt{j}^{\dagger} ) \leftarrow \mathbf{r}_{j} \,  \mathbf{U}( \overline{\texttt{j}},\texttt{j}^{\dagger} )$
	\STATE $\mathbf{V}( \texttt{j}^{\dagger}, \overline{\texttt{j}} ) \leftarrow \mathbf{t}_{j}  \mathbf{V}( \texttt{j}^{\dagger}, \overline{\texttt{j}} )$

	\FORALL{$\Omega_{q}^{(l)} \in \mathcal{N}_{j}^{(l)} \cup \mathcal{I}_{j}^{(l)} \backslash \Omega_{k}^{(l)}$}
		\STATE $\mathbf{K}(\overline{\texttt{q}}, \overline{\texttt{j}}) \leftarrow \mathbf{K}(\overline{\texttt{q}}, \overline{\texttt{j}}) \, \mathbf{t}^{\mathrm{T}}_{j}$
		\STATE $\mathbf{K}(\overline{\texttt{j}}, \overline{\texttt{q}}) \leftarrow \mathbf{r}_{j} \, \mathbf{K}(\overline{\texttt{j}}, \overline{\texttt{q}})$
	\ENDFOR
	
	\IF{$\underline{\texttt{k}}$ and $\texttt{k}$ of $\Omega_{k}^{(l)}$ have been eliminated}
		\STATE $\mathbf{K}(\overline{\texttt{k}}, \overline{\texttt{j}}) \leftarrow \mathbf{K}(\overline{\texttt{k}}, \overline{\texttt{j}}) \, \mathbf{t}^{\mathrm{T}}_{j} + \mathbf{K}^{\prime}(\overline{\texttt{k}},\overline{\texttt{j}}) \, \mathbf{t}^{\prime^{\mathrm{T}}}_{j}$
		\STATE $\mathbf{K}(\overline{\texttt{j}}, \overline{\texttt{k}}) \leftarrow \mathbf{r}_{j} \, \mathbf{K}(\overline{\texttt{j}}, \overline{\texttt{k}}) + \mathbf{r}^{\prime}_{j} \, \mathbf{K}^{\prime}(\overline{\texttt{j}},\overline{\texttt{k}})$
	\ELSE
		\STATE 	$\mathbf{K}(\overline{\texttt{k}}, \overline{\texttt{j}}) \leftarrow \mathbf{r}_{k} \, \mathbf{K}(\overline{\texttt{k}}, \overline{\texttt{j}}) \, \mathbf{t}^{\mathrm{T}}_{j} + \mathbf{r}^{\prime}_{k} \, \mathbf{K}^{\prime}(\overline{\texttt{k}},\overline{\texttt{j}}) \, \mathbf{t}^{\prime^{\mathrm{T}}}_{j} $
	\ENDIF
	\end{algorithmic}
	\caption{The inverse fast multipole method: updating of the edges associated with $\Omega_{j}^{(l)}$ after compressing the interactions between $\Omega_{j}^{(l)}$ and $\Omega_{k}^{(l)}$}
	\label{alg:IFMM_update}
\end{algorithm}
%

\subsection{Some important remarks}\label{subsec:remarks}
Several techniques can be employed for obtaining a suitable low-rank approximation for the dense fill-ins~$\mathbf{E}^{\prime}$ arising throughout the elimination phase (cf., \autoref{eq:P2P_kj} and \eqref{eq:P2P_jk}). The best possible approximation of rank $k$ of a matrix (minimizing the error norm) is given by its $k^{\mathrm{th}}$ partial singular value decomposition (SVD)~\cite{golu96a}. Unfortunately, the calculation of a complete SVD requires $\mO\left(\min \left(m n^2, m^2 n \right)\right)$ arithmetical operations for a $(m \times n)$ matrix, making it rather unattractive from a computational point of view. The randomized SVD (rSVD) algorithm~\cite{halk11a} provides a very efficient alternative for computing (approximations of) the first few singular values and vectors. This method exploits the power of randomization to identify a subspace that captures most of the action of a matrix~\cite{halk11a}. It involves generating random matrices, performing matrix-matrix products, applying orthonormalizations, and computing the SVD of a much smaller matrix. Although it uses random matrices, the rSVD algorithm gives very accurate results with a high probability, provided that the singular values of the matrix decay sufficiently fast. Rather than specifying a fixed rank for the low-rank approximation of a dense fill-in, the rank can also be determined by prescribing a relative accuracy $\varepsilon$. The relative importance of a fill-in is hence evaluated by weighting its singular values (obtained, for example, with the rSVD algorithm) with respect to some reference value. In this paper, it is proposed to employ the largest singular value $\sigma_{0}(\mathbf{E})$ of the extended sparse matrix~$\mathbf{E}$ for this purpose. The retained rank of a fill-in $\mathbf{E}^{\prime}$ is thus given by\footnote{$k$ is also bounded by the size of $\mathbf{E}^{\prime}$.}:
\begin{equation}
	k(\varepsilon) = \min \left\{ k \in \mathbb{N} : \sigma_{k+1}(\mathbf{E}^{\prime}) \leq  \varepsilon \sigma_{0}(\mathbf{E}) \right\}
	\label{eq:methods_SVD_rank_spectralnorm}
\end{equation}
$\sigma_{0}(\mathbf{E})$ can easily be estimated by application of the rSVD algorithm to the extended sparse matrix~$\mathbf{E}$.

A low-rank recompression is also required for calculating new inter- and anterpolation operators of a cluster after compression of the fill-ins (cf., \autoref{eq:U_j_recomp} and~\eqref{eq:V_j_recomp}). Once again, computing the SVD will provide the optimal recompression, but at a relatively high computational cost. Applying the rSVD algorithm will not be very efficient either in this case, as the rank of the new operators $\widehat{\mathbf{U}}(\underline{\texttt{j}}, \texttt{j})$ and $\widehat{\mathbf{V}}(\texttt{j},\underline{\texttt{j}})$ is at least as large as the rank of the original operators $\mathbf{U}(\underline{\texttt{j}}, \texttt{j})$ and $\mathbf{V}(\texttt{j},\underline{\texttt{j}})$ (as the latter consist of orthogonal columns). This recompression hence differs from the compression of the fill-ins, where only the first few singular values are (expected to be) of importance. The fully pivoted adaptive cross approximation (ACA) algorithm~\cite{bebe03a,rjas07a}\footnote{This approach is the same as an LU factorization with full pivoting.}, followed by an SVD recompression (i.e., two skinny QR decompositions and an SVD performed on the product of the resulting triangular matrices), provides a more efficient way for performing the recompression and is therefore favored.

A final remark is related to the weighting of the basis vectors when computing new inter- and anterpolation operators, as presented in \autoref{eq:U_j_recomp} and \eqref{eq:V_j_recomp}. The advantage of using rSVD and ACA+SVD is that the obtained singular values straightforwardly provide the weights for each vector in the basis. If the low-rank approximations in \autoref{eq:P2P_kj}, \eqref{eq:P2P_jk}, \eqref{eq:U_j_recomp} and \eqref{eq:V_j_recomp} do not correspond to singular value decompositions (e.g., when using a rank-revealing LU decomposition), the weighting needs to be chosen in a different manner. The pivots of the LU factorization could, for example, be employed for this purpose. Although this approach is less accurate than using the singular values, it can still provide a reasonable measure of the scale of each basis vector when performing the recompression.

\subsection{Computational cost}
We will derive the computational cost in the case of a uniform tree decomposition of the matrix (as opposed to adaptive), although the proof extends to general adaptive trees. 

Let's assume that we have a matrix $\mathbf{A}$ in $\mH^2$-format. Given a tolerance $\varepsilon$, we assume that there is an upper bound $r$ on the rank of all compressed blocks such that $r \in \mO(\log 1/\varepsilon)$. One can prove that for the original matrix $\mathbf{A}$, in the case of smooth kernels (e.g., Green's functions), the rank $r$ of the compressed blocks does satisfy this property~\cite{fong09a}. However, in our case, as we proceed through the LU factorization, recompression is required. We have not proved that the rank does not increase in an unbounded fashion (for a fixed $\varepsilon$, as $N$ varies) as a result of this process. However, we did observe experimentally that the rank is in fact bounded. See in particular \autoref{fig:benchmark_ii_3D_sphere_paper_rank_vs_N_n_4}.

Under these assumptions, the cost of operating on a node at a given level is $\mO(r^3)$. The cost of a leaf node elimination is $\mO( N_0^3 )$, where $N_0$ is an upper bound on the size of all leaf clusters. Choose $N_0 \in \Omega(r)$. This requirement sets the number of levels in the tree, and the number of leaf nodes. In that case, the cost of the elimination of all leaf nodes is $\mO( N r^2 )$. The cost of compression of the far-away fill-in blocks, and computation of the updated parent-level matrix blocks is also $\mO( N r^2 )$. Since at the next level, the number of clusters is divided at least by a fixed factor (e.g., 2 for a binary tree or 8 for an octree), the total computational cost of the method is $\mO( N r^2 )$ or $\mO( N \log(1/\varepsilon)^2 )$.

\section{Numerical examples}\label{sec:NumericalExamples}
The IFMM algorithm has been implemented in \texttt{C++}, using the \texttt{Eigen}-library~\cite{eigen-3} for performing dense linear algebra operations. In the following subsections, numerical examples are considered to validate the implementation and to demonstrate the efficiency of the methodology. All computations have been performed on Intel\textsuperscript{\textregistered} Xeon\textsuperscript{\textregistered} E5-2650 v2 (2.60~GHz) CPUs. In each of these examples, the initial low-rank approximations in the $\mH^2$-matrix are created using interpolation based on Chebyshev polynomials~\cite{fong09a}, followed by an additional SVD to reduce the rank. A uniform octree decomposition of the domain is employed; the number of levels is adjusted to ensure that the leaf clusters contain approximately $100$~points. This value was found to provide a reasonable trade-off between the time spent at the leaf level and the time spent at higher levels in the tree.\footnote{Less points per leaf cluster makes the elimination of the leaves faster, but increases the number of levels and thus the time spent in the tree, and vice versa.} A relative accuracy $\varepsilon = 10^{-3}$ is used, both for compressing the fill-ins (with rSVD) and for updating the inter- and anterpolation operators (with ACA).

\subsection{The IFMM as a direct solver}\label{subsec:ifmm_direct}

\begin{figure}[hbtp]
	\begin{center}
		\includegraphics*[width=0.475\linewidth,clip=true]{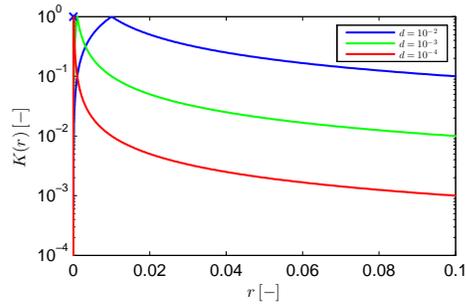}
	\end{center}
	\caption{Kernel function $K(r)$ for three values of the parameter $d$.}
\label{fig:kernel_benchmark_ii_semilogy_paper}
\end{figure}

\begin{figure}[hbtp]
	\begin{center}
	(a)~\hspace{-1.5em}\includegraphics*[width=0.475\linewidth,clip=true]{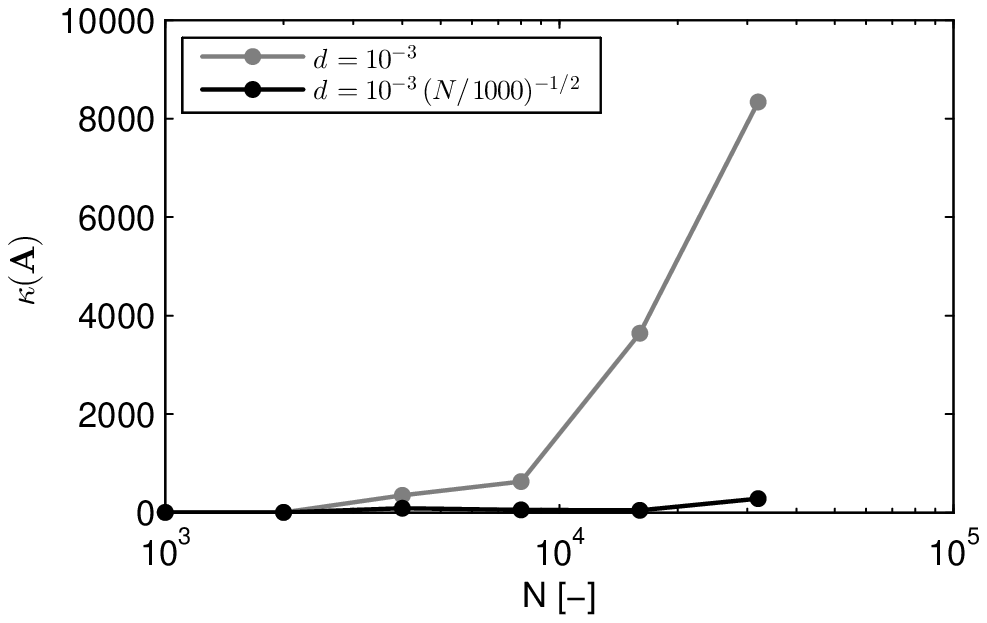}
	(b)~\hspace{-1.5em}\includegraphics*[width=0.475\linewidth,clip=true]{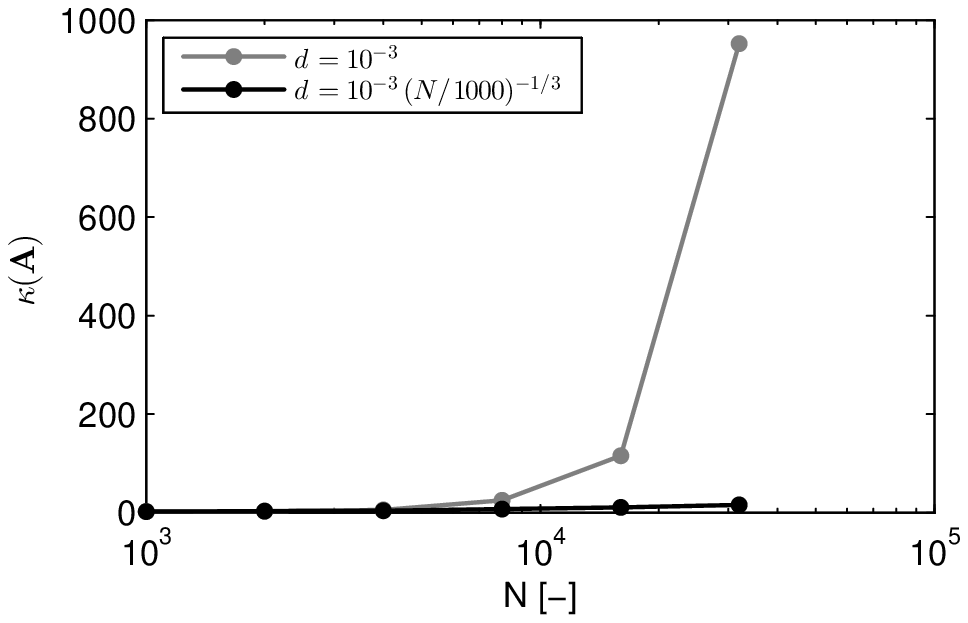}
	\end{center}
	\caption{Condition number $\kappa(\mathbf{A})$ as a function of the matrix size $N$, for points randomly distributed (a)~on the unit sphere and (b)~in the cube $\left[-1,1\right]^3$.}
\label{fig:benchmark_ii_conditionnumber}
\end{figure}

Consider the linear system of equations $\mathbf{A} \mathbf{x} = \mathbf{b}$ with matrix entries $a_{ij} = K(\|\mathbf{r}_i - \mathbf{r}_j \|)$ depending on the kernel $K(r)$:
\begin{equation}
  K(r)=
  \begin{cases}
    1 & \text{if} \quad r =0, \\
    \frac{r}{d} & \text{if} \quad 0<r<d, \\
    \frac{d}{r} & \text{if} \quad r \geq d.
  \end{cases}
  \label{eq:K_ij}
\end{equation}
with $\mathbf{r}_i$ and $\mathbf{r}_j$ indicating 3D position vectors and $r = \|\mathbf{r}_i - \mathbf{r}_j \|$. \autoref{fig:kernel_benchmark_ii_semilogy_paper} shows the kernel function $K(r)$ for three values of the parameter $d$ in \autoref{eq:K_ij}. This parameter strongly affects the condition number $\kappa(\mathbf{A})$ of the matrix $\mathbf{A}$: the larger $d$, the larger the condition number (for a fixed matrix size $N$). The condition number also increases with $N$ for a fixed value of $d$, as illustrated in \autoref{fig:benchmark_ii_conditionnumber}. Note that the matrix entries are directly defined by the kernel, rather than by a boundary integral formulation that involves integration of the kernel. Instead of using a kernel $K(r) = 1/r$, the kernel is modified near $r=0$ according to \autoref{eq:K_ij} in order to ensure that bounded matrix entries are obtained on/near the diagonal of $\mathbf{A}$.

In the following, a known vector $\mathbf{x}$ is chosen and the matrix-vector product $\mathbf{A} \mathbf{x}$ is evaluated for obtaining $\mathbf{b}$. The IFMM is subsequently used as a direct solver, with the aforementioned vector $\mathbf{b}$ as right hand side. The retrieved solution is denoted as~$\widehat{\mathbf{x}}$.

\begin{figure}[hbtp]
	\begin{center}
	(a)~\hspace{-1.5em}\includegraphics*[width=0.475\linewidth,clip=true]{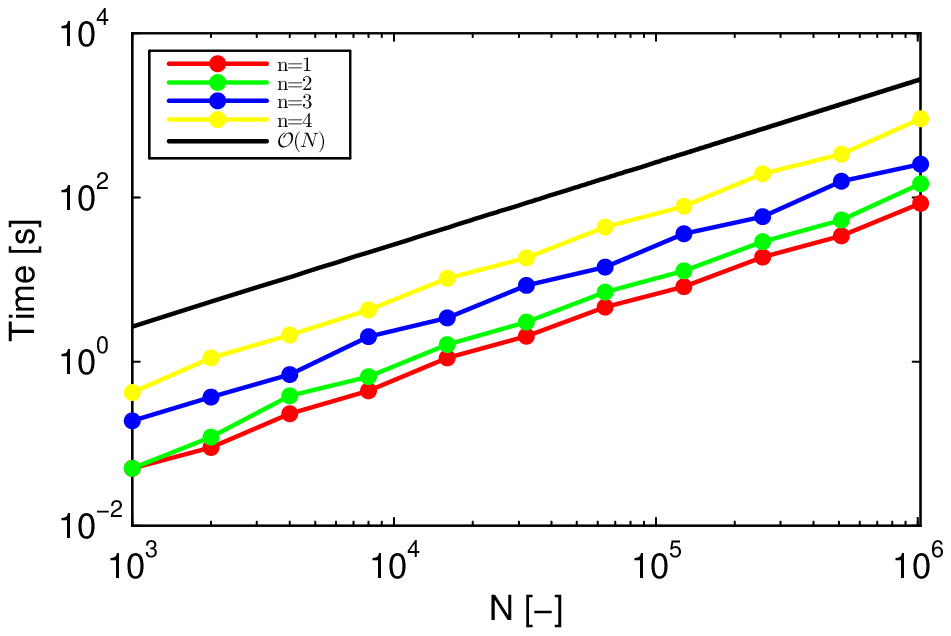}
	(b)~\hspace{-1.5em}\includegraphics*[width=0.475\linewidth,clip=true]{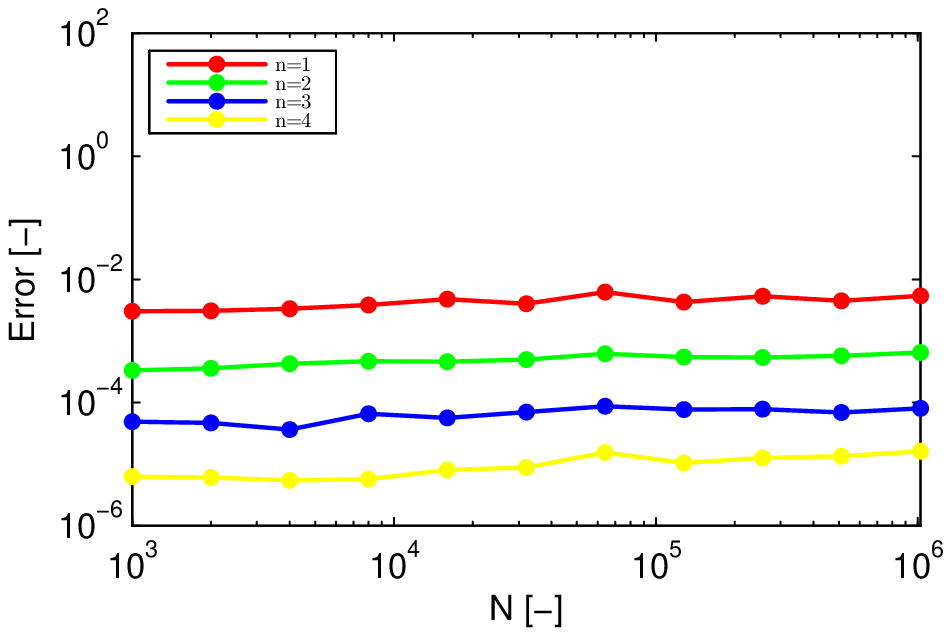}
	\end{center}
	\caption{(a)~Total CPU time and (b)~the relative error $\|\mathbf{x} -  \widehat{\mathbf{x}} \|_{2}/ \| \mathbf{x} \|_{2}$ as a function of the matrix size~$N$, for points randomly distributed on the unit sphere.}
\label{fig:benchmark_ii_3D_sphere}
\end{figure}

In a first example, $N$ points $\mathbf{r}$ are randomly distributed on the unit sphere. The parameter $d$ is scaled according to $d =10^{-3} \, (N/1000)^{-1/2}$ in order to ensure that the condition number of $\mathbf{A}$ remains nearly constant for all $N$ (see \autoref{fig:benchmark_ii_conditionnumber}(a)). \autoref{fig:benchmark_ii_3D_sphere}(a) shows the total time taken by the IFMM solver as a function of the number of points $N$, with the number of Chebyshev nodes~$n$ (in each dimension) varying from $1$ to $4$ (the total number of Chebyshev interpolation nodes employed for constructing the initial low-rank approximations hence equals $n^3$). These curves demonstrate that, as expected, a linear complexity $\mO(N)$ is achieved. The computation time increases for larger values of~$n$. The relative error $\|\mathbf{x} -  \widehat{\mathbf{x}} \|_{2}/ \| \mathbf{x} \|_{2}$ for each of these cases is depicted in \autoref{fig:benchmark_ii_3D_sphere}(b). The figure clearly shows a decrease of the error as $n$ is increased. The errors are reasonably small for all problem sizes~$N$.

\begin{figure}[hbtp]
	\begin{center}
	\includegraphics*[width=0.475\linewidth,clip=true]{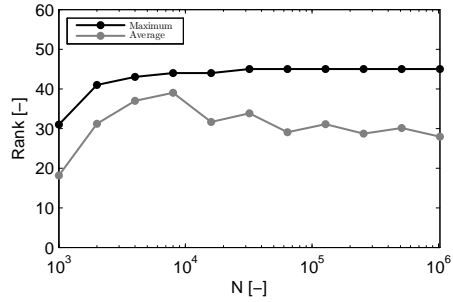}
	\end{center}
	\caption{Maximum (black line) and average (grey line) rank of the matrix blocks corresponding to far field interactions as a function of the matrix size~$N$, for points randomly distributed on the unit sphere (with $n=4$).}
\label{fig:benchmark_ii_3D_sphere_paper_rank_vs_N_n_4}
\end{figure}

A crucial feature for achieving the linear complexity demonstrated in \autoref{fig:benchmark_ii_3D_sphere}(a) is the fact that the rank of the matrix blocks corresponding to far field interactions remains (almost) constant if the matrix size~$N$ is increased (due to the strong admissibility criterion employed for $\mH^2$-matrices). This is in contrast with HODLR and HSS matrices, where the rank grows like $\mO(\sqrt{N})$ for 2D problems~\cite{ambi13b}. \autoref{fig:benchmark_ii_3D_sphere_paper_rank_vs_N_n_4} shows the maximum and average rank of the aforementioned matrix blocks as a function of $N$ for the case of $n=4$, confirming the previous statement (especially for $N > 10^{4}$). The kinks in the curve of the average rank are caused by varying the number of levels in the octree.

\begin{figure}[hbtp]
	\begin{center}
	(a)~\hspace{-1.5em}\includegraphics*[width=0.475\linewidth,clip=true]{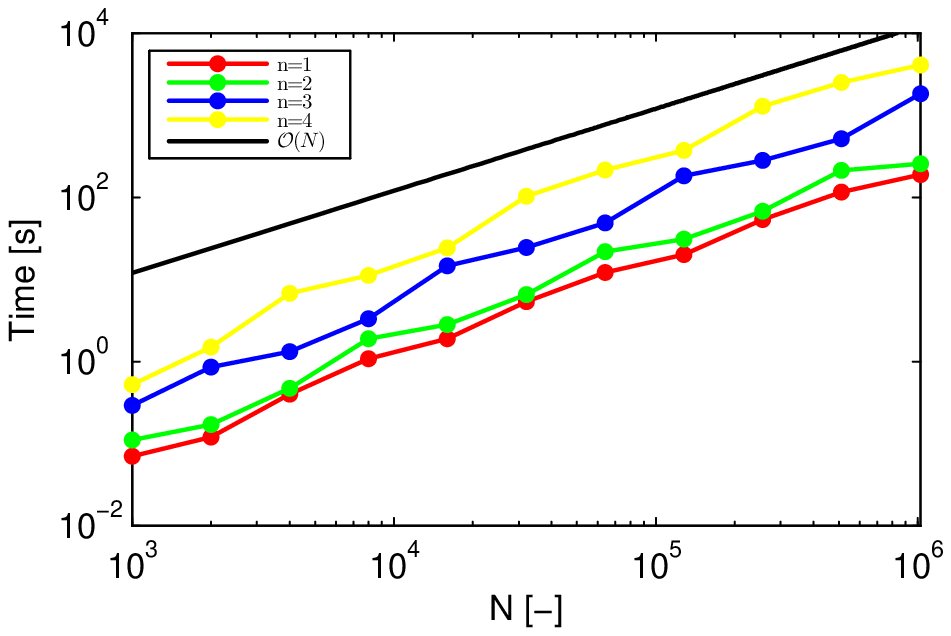}
	(b)~\hspace{-1.5em}\includegraphics*[width=0.475\linewidth,clip=true]{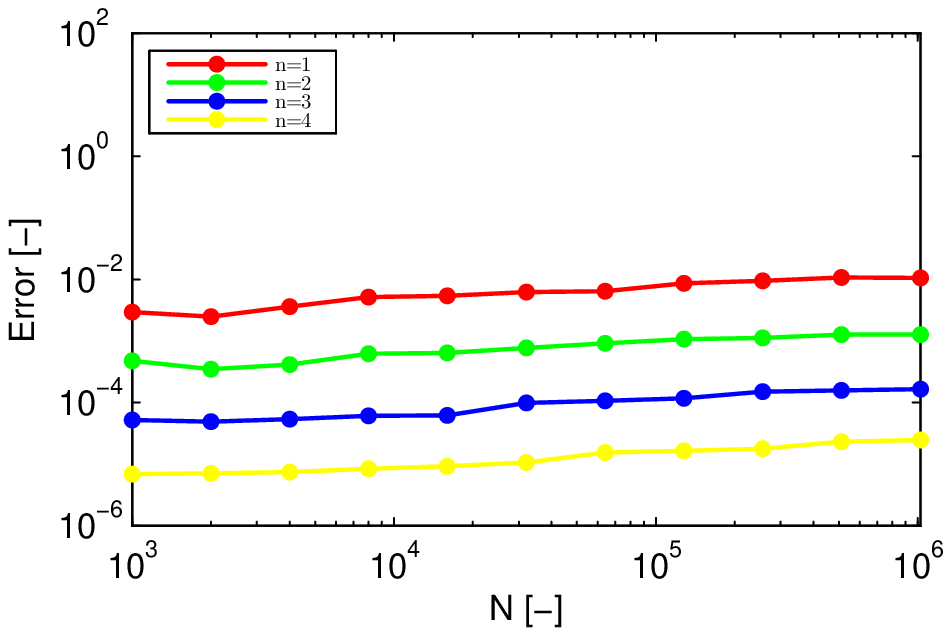}
	\end{center}
	\caption{(a)~Total CPU time and (b)~the relative error $\|\mathbf{x} -  \widehat{\mathbf{x}} \|_{2}/ \| \mathbf{x} \|_{2}$ as a function of the matrix size~$N$, for points randomly distributed in the cube $\left[-1,1\right]^3$.}
\label{fig:benchmark_ii_3D_random}
\end{figure}

Next, a more general example of $N$ points $\mathbf{r}$ randomly distributed in the cube $\left[-1,1\right]^3$ is investigated. The parameter $d$ varies in this case as $d =10^{-3} \, (N/1000)^{-1/3}$ to keep the condition number of $\mathbf{A}$ nearly constant for all $N$ (see \autoref{fig:benchmark_ii_conditionnumber}(b)). \autoref{fig:benchmark_ii_3D_random}(a) confirms that a linear scaling of the computation time is also obtained in this case (although only from $N \sim 10^{4}$ upwards), while \autoref{fig:benchmark_ii_3D_random}(b) demonstrates the convergence of the relative error as the number of Chebyshev nodes is increased. The computation times are slightly higher than in the case of points distributed on a sphere (\autoref{fig:benchmark_ii_3D_sphere}(a)). This is due to the fact that a uniform octree decomposition of a sphere leads to many empty clusters, implying that most of the clusters will have less than 27 neighbors, while it is likely that almost all clusters will possess the maximum number of neighbors in case the points are randomly distributed in a cube. The reader is referred to \autoref{app:} for a detailed breakdown of the computational cost for this example.

\begin{figure}[hbtp]
	\begin{center}
		\includegraphics*[width=0.475\linewidth,clip=true]{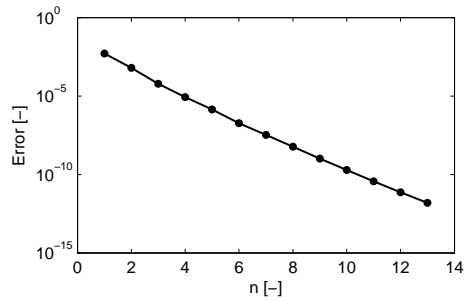}
	\end{center}
	\caption{The relative error $\|\mathbf{x} -  \widehat{\mathbf{x}} \|_{2}/ \| \mathbf{x} \|_{2}$ as a function of the number of Chebyshev nodes $n$, for $N=8000$ points randomly distributed in the cube $\left[-1,1\right]^3$.}
\label{fig:Sherlock_iFMM_paper_3D_randompoints_a_1e-3_scaled_N_8000_convergence}
\end{figure}

Depending on the application of interest, the errors shown in \autoref{fig:benchmark_ii_3D_sphere}(b) and~\ref{fig:benchmark_ii_3D_random}(b) might already be sufficiently small. The convergence of the IFMM is investigated in more detail in \autoref{fig:Sherlock_iFMM_paper_3D_randompoints_a_1e-3_scaled_N_8000_convergence}, which shows the relative error $\|\mathbf{x} -  \widehat{\mathbf{x}} \|_{2}/ \| \mathbf{x} \|_{2}$ as a function of the number of Chebyshev nodes for a fixed matrix size $N=8 \times 10^{3}$. An increase of $n$ leads to more accurate low-rank approximations of the far field interactions and consequently a smaller overall error in $\widehat{\mathbf{x}}$. This figure demonstrates that the method is capable of serving as a highly accurate direct solver. The number of Chebyshev nodes needs to be increased significantly if 10 or more digits of accuracy are desired, however, which leads to a strong increase of the computational cost (even if the linear complexity is preserved). In most cases (depending on the number of right-hand sides), it is more efficient to apply a low accuracy IFMM solver as a preconditioner in an iterative method rather than using a highly accurate direct IFMM solver, as will be illustrated in the next subsection.

\subsection{The IFMM as a preconditioner in an iterative scheme}\label{subsec:ifmm_precon}
In this subsection, the linear system $ \mathbf{A} \mathbf{x} = \mathbf{b}$ is solved iteratively using a non-restarted GMRES-algorithm~\cite{saad93a,saad86a} with a tolerance $\varepsilon = 10^{-10}$ for the relative residual $\| \mathbf{b} - \mathbf{A} \widehat{\mathbf{x}} \|_{2}/ \| \mathbf{b} \|_{2}$ and a maximum of $500$ iterations. The same kernel $K(r)$ as in \autoref{eq:K_ij} is used. The black-box FMM~\cite{fong09a} is employed for evaluating the matrix-vector products within GMRES, while a low accuracy IFMM solver is applied as right preconditioner in each iterative step. A comparison of right and left preconditioning is provided in \autoref{app:leftright}.

\begin{figure}[hbtp]
	\begin{center}
	(a)~\hspace{-1.5em}\includegraphics*[width=0.475\linewidth,clip=true]{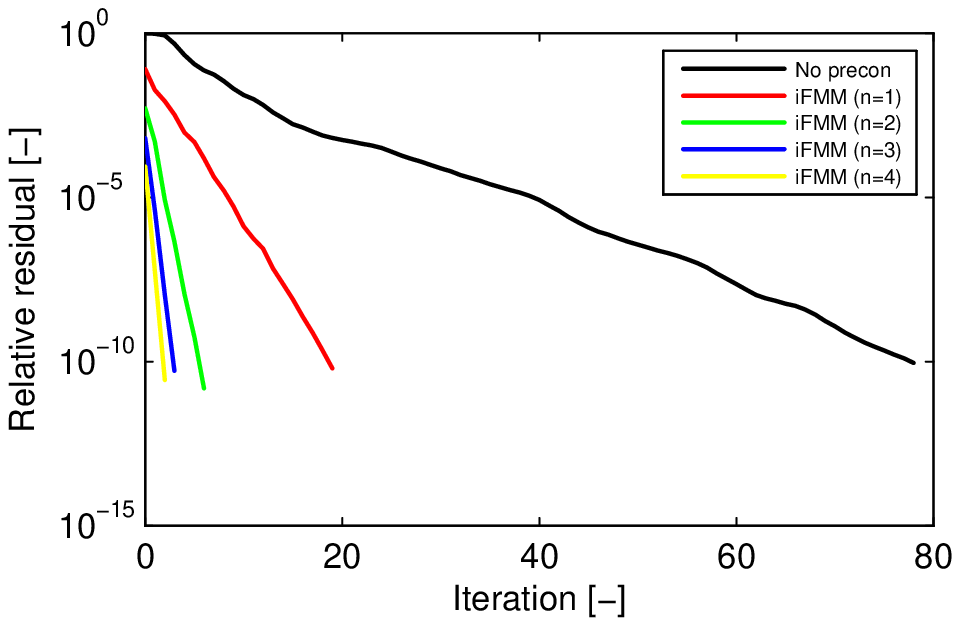}
	(b)~\hspace{-1.5em}\includegraphics*[width=0.475\linewidth,clip=true]{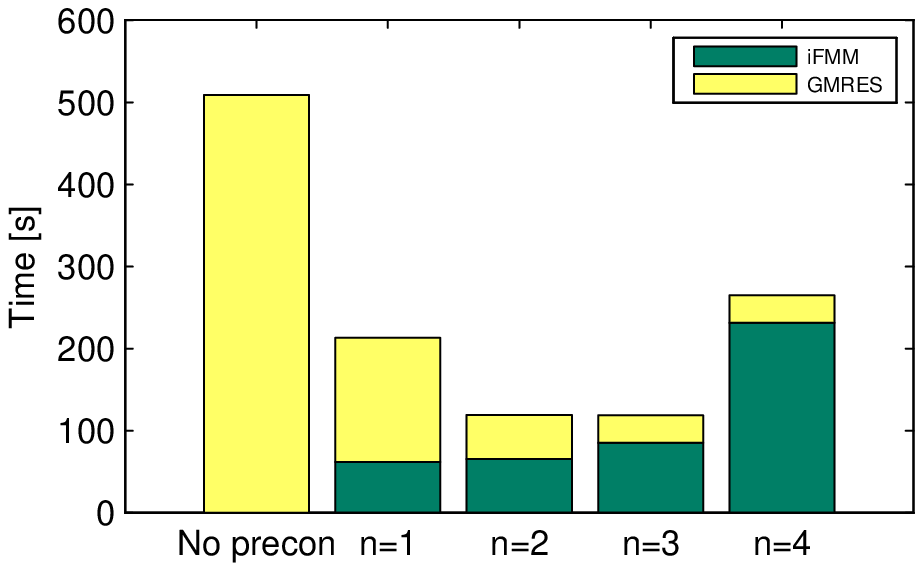}
	\end{center}
	\caption{(a)~Relative residual $\| \mathbf{b} -  \mathbf{A} \widehat{\mathbf{x}} \|_{2}/ \| \mathbf{b} \|_{2}$ as a function of the iteration step and (b)~total CPU time for $d=10^{-3}$.}
	\label{fig:GMRES_3D_randompoints_benchmark_ii_a_1e-3_plot}
\end{figure}

The first example involves $N =10^5$ points randomly distributed in the cube $\left[-1,1\right]^3$, with $d=10^{-3}$ (in \autoref{eq:K_ij}). \autoref{fig:GMRES_3D_randompoints_benchmark_ii_a_1e-3_plot}(a) shows the relative residual as a function of the number of GMRES-iterations for several preconditioning approaches. Without preconditioning, convergence is rather slow and a total of $78$ iterations are required to achieve convergence. The application of the IFMM as preconditioner results in a much faster convergence. Even with a rank-1 approximation of the far field interactions (i.e., $n=1$), a significant reduction of the number of iterations is obtained. An increase of the number of Chebyshev nodes leads to a more accurate preconditioner and hence a further reduction of the number of iterations; this comes at the price of a higher cost to construct the preconditioner, however. The computation times for the considered preconditioning strategies are presented in \autoref{fig:GMRES_3D_randompoints_benchmark_ii_a_1e-3_plot}(b), where a decomposition is made into the time required to construct the preconditioner (i.e., the elimination phase in the IFMM) and the actual iteration time (i.e., the substitution phase in the IFMM, together with all GMRES-related computations). For the example under concern, the use of three Chebyshev nodes per dimension leads to the smallest overall computation time. The results are summarized in \autoref{tbl:GMRES_3D_randompoints_benchmark_ii_a_1e-3_tbl}.

\begin{table}[htbp]
\footnotesize
\centering
\caption{Overview of the results for various preconditioning strategies for $d=10^{-3}$.}
\label{tbl:GMRES_3D_randompoints_benchmark_ii_a_1e-3_tbl}
\begin{tabular}{c| c c c c}
	\toprule
    	&	\# iterations	&	Total CPU time [s]	& 	$\| \mathbf{b} -  \mathbf{A} \widehat{\mathbf{x}} \|_{2}/ \| \mathbf{b} \|_{2}$	  &	$\|\mathbf{x} -  \widehat{\mathbf{x}} \|_{2}/ \| \mathbf{x} \|_{2}$\\
	\midrule
	No precon 	& 78& 509  & $9.1 \times 10^{-11}$ & $4.6 \times 10^{-10}$ \\
	IFMM ($n=1$)& 19& 213  & $6.1 \times 10^{-11}$ & $4.0 \times 10^{-11}$ \\
	IFMM ($n=2$)& 6	& 119  & $1.5 \times 10^{-11}$ & $9.6 \times 10^{-12}$ \\
	IFMM ($n=3$)& 3	& 118  & $5.2 \times 10^{-11}$ & $2.9 \times 10^{-11}$ \\
	IFMM ($n=4$)& 2	& 265  & $2.7 \times 10^{-11}$ & $2.1 \times 10^{-11}$ \\
	\bottomrule
\end{tabular}
\end{table}

\begin{figure}[hbtp]
	\begin{center}
	(a)~\hspace{-1.5em}\includegraphics*[width=0.475\linewidth,clip=true]{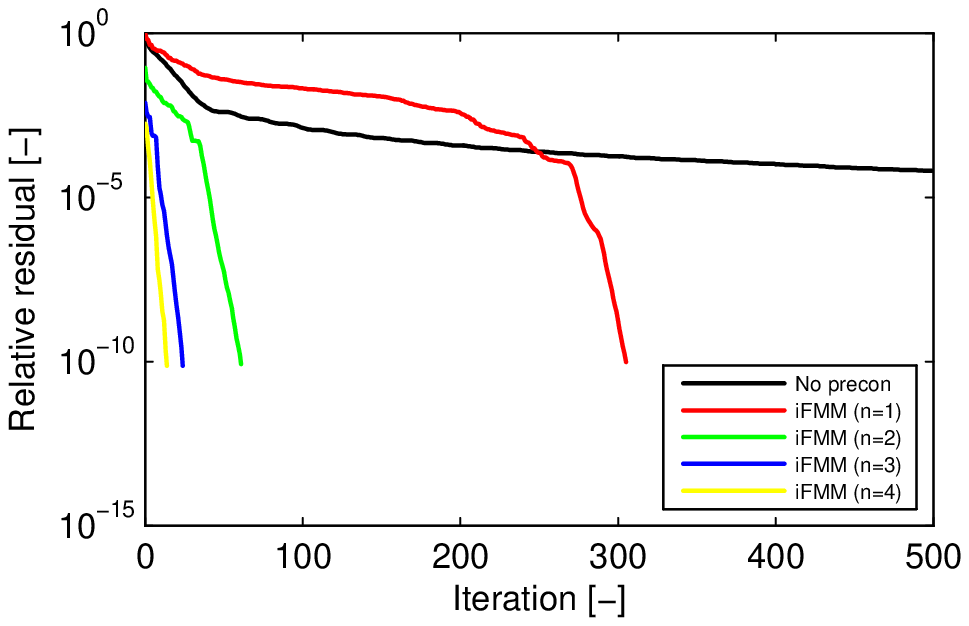}
	(b)~\hspace{-1.5em}\includegraphics*[width=0.475\linewidth,clip=true]{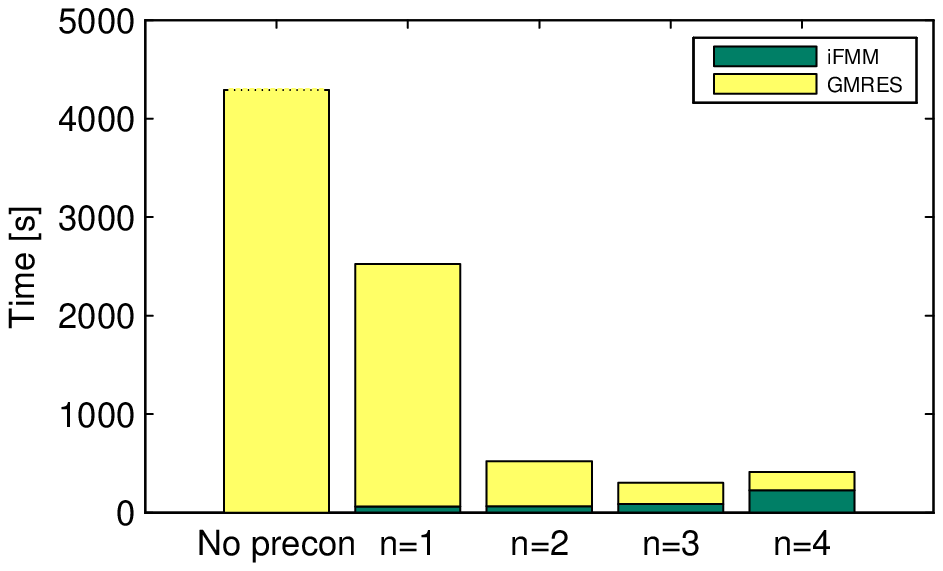}
	\end{center}
	\caption{(a)~Relative residual $\| \mathbf{b} -  \mathbf{A} \widehat{\mathbf{x}} \|_{2}/ \| \mathbf{b} \|_{2}$ as a function of the iteration step and (b)~total CPU time for $d=10^{-2}$.}
\label{fig:GMRES_3D_randompoints_benchmark_ii_a_1e-2_plot}
\end{figure}

In the second example, the parameter $d$ in \autoref{eq:K_ij} is increased from $10^{-3}$ to $10^{-2}$, resulting in a more ill-conditioned matrix $\mathbf{A}$. \autoref{fig:GMRES_3D_randompoints_benchmark_ii_a_1e-2_plot}(a) demonstrates that the unpreconditioned GMRES-algorithm is now unable to converge to the specified residual within $500$ iterations. Application of the IFMM as preconditioner ensures convergence, although a rank-1 approximation of the far field is not very efficient in this case, as the number of iterations remains rather large. As soon as two or more Chebyshev nodes are employed, a remarkable reduction of the number of iterations is obtained and convergence is much faster. \autoref{fig:GMRES_3D_randompoints_benchmark_ii_a_1e-2_plot}(b) and \autoref{tbl:GMRES_3D_randompoints_benchmark_ii_a_1e-2_tbl} illustrate that the use of three Chebyshev nodes per dimension results once more in the smallest overall computation time. It is also important to note that although the relative residual $\| \mathbf{b} -  \mathbf{A} \widehat{\mathbf{x}} \|_{2}/ \| \mathbf{b} \|_{2}$ reaches the specified tolerance $\varepsilon$, the actual relative error $\|\mathbf{x} -  \widehat{\mathbf{x}} \|_{2}/ \| \mathbf{x} \|_{2}$ is a few order of magnitudes larger. This is due to the large condition number of $\mathbf{A}$.

\begin{table}[htbp]
\footnotesize
\centering
\caption{Overview of the results for various preconditioning strategies for $d=10^{-2}$. The results are surrounded with parentheses if no convergence was achieved.}
\label{tbl:GMRES_3D_randompoints_benchmark_ii_a_1e-2_tbl}
\begin{tabular}{c| c c c c}
	\toprule
    	&	\# iterations	&	Total CPU time [s]	& 	$\| \mathbf{b} -  \mathbf{A} \widehat{\mathbf{x}} \|_{2}/ \| \mathbf{b} \|_{2}$	  &	$\|\mathbf{x} -  \widehat{\mathbf{x}} \|_{2}/ \| \mathbf{x} \|_{2}$\\
	\midrule
	No precon 	& (500)	& (4290) & ($6.3 \times 10^{-5}$) & ($4.2 \times 10^{-2}$) \\
	IFMM ($n=1$)& 305	& 2522   & $9.7 \times 10^{-11}$ & $3.8 \times 10^{-8}$ \\
	IFMM ($n=2$)& 61	& 521    & $8.3 \times 10^{-11}$ & $2.2 \times 10^{-7}$ \\
	IFMM ($n=3$)& 24 	& 301    & $7.3 \times 10^{-11}$ & $4.6 \times 10^{-7}$ \\
	IFMM ($n=4$)& 14	& 414    & $7.3 \times 10^{-11}$ & $7.9 \times 10^{-7}$ \\
	\bottomrule
\end{tabular}
\end{table}

\begin{figure}[hbtp]
	\begin{center}
	(a)~\hspace{-1.5em}\includegraphics*[width=0.475\linewidth,clip=true]{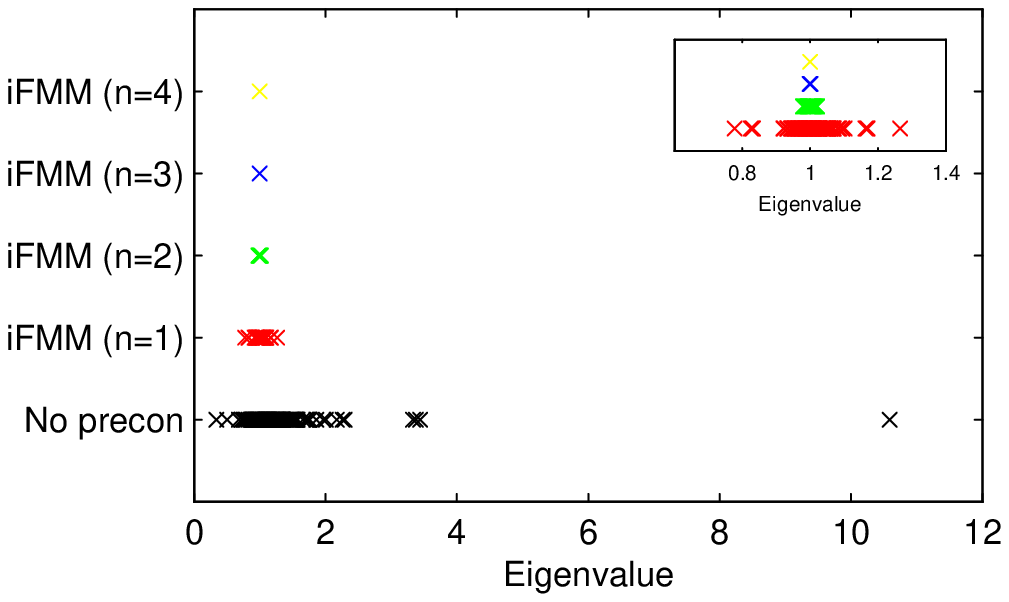}
	(b)~\hspace{-1.5em}\includegraphics*[width=0.475\linewidth,clip=true]{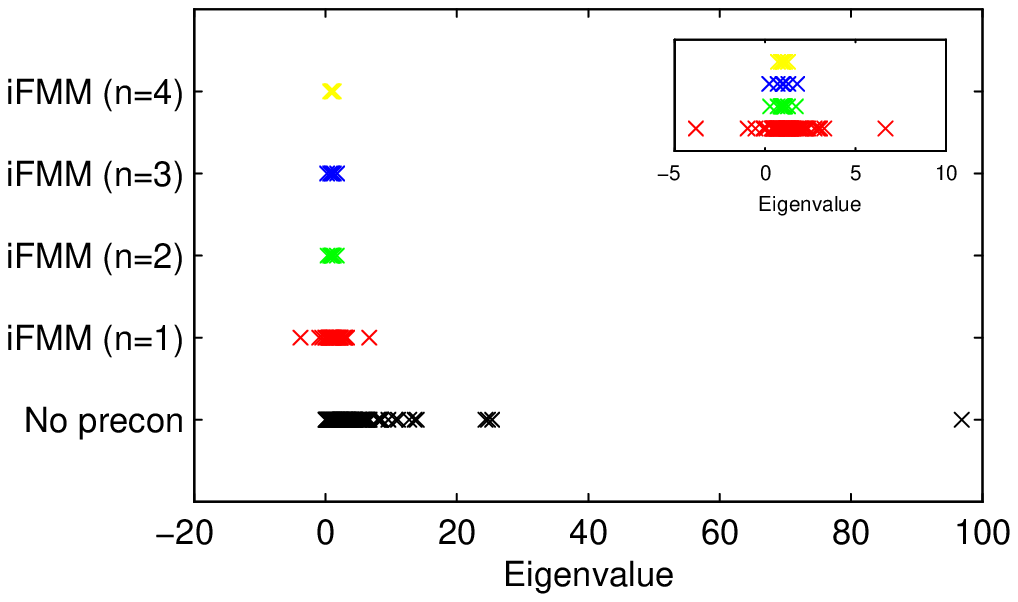}
	\end{center}
	\caption{Eigenvalue distribution of the original (black crosses) and the preconditioned (colored crosses) matrices of size $N =10^4$, for (a)~$d=10^{-3}$ and (b)~$d=10^{-2}$.}
\label{fig:3D_randompoints_benchmark_ii_iFMM_eigenvalues_N_10000}
\end{figure}

In order to provide more insight in the effect of the IFMM as preconditioner, the eigenvalue distributions of the original and preconditioned matrices $\mathbf{A}$ and $\mathbf{P}^{-1}\mathbf{A}$ are investigated, where $\mathbf{P}^{-1}$ represents the IFMM preconditioner. A smaller problem of size $N =10^4$ is considered here. It is furthermore emphasized that the `true' matrix $\mathbf{A}$ is used for this purpose, and not an FMM approximation of $\mathbf{A}$. \autoref{fig:3D_randompoints_benchmark_ii_iFMM_eigenvalues_N_10000}(a) and \ref{fig:3D_randompoints_benchmark_ii_iFMM_eigenvalues_N_10000}(b) show the eigenvalues of $\mathbf{A}$ and $\mathbf{P}^{-1}\mathbf{A}$ for $d=10^{-3}$ and $d=10^{-2}$, respectively. As $\mathbf{A}$ is symmetric, all eigenvalues are real. It is clear from these figures that the eigenvalues of $\mathbf{A}$ are not well clustered. This explains the slow convergence of GMRES if no preconditioning is employed. Application of the IFMM as preconditioner leads to a better clustering of the eigenvalues. An increase of the number of Chebyshev nodes results in a more compact clustering of the eigenvalues near $1$, and consequently a faster convergence of the iterative solver.

The results presented in this subsection clearly illustrate the effectiveness of the IFMM as preconditioner for an iterative solver. There is a trade-off between the preconditioner's accuracy (and thus the number of iterations) and its construction cost, which has to be evaluated on a case-by-case basis.

\section{Application: rigid bodies immersed in a Stokes flow}\label{sec:Application}
In order to illustrate the applicability of the IFMM to challenging science and engineering problems (as opposed to the synthetic examples discussed in \autoref{sec:NumericalExamples}), an application involving rigid bodies immersed in a Stokes flow is considered in this section. These rigid bodies represent, for example, colloidal particles suspended in a solvent~\cite{delo14a}. Each rigid body is discretized with a number of points $\mathbf{r}_i$ on the surface, and a force vector $\mathbf{F}_{i} \in \mathbb{R}^{3 \times 1}$ and velocity vector $\mathbf{U}_{i} \in \mathbb{R}^{3 \times 1}$ are associated with each point. The linear system of equations that needs to be solved reads as:
\begin{equation}
	\mathbf{\mathcal{{M}}} \mathbf{F} = \mathbf{U} 
	\label{eq:system_Stokes}
\end{equation}
where $\mathbf{F} \in \mathbb{R}^{3N \times 1} $ and $\mathbf{U} \in \mathbb{R}^{3N \times 1}$ collect the forces and velocities of all $N$ points, respectively. The mobility matrix $\mathbf{\mathcal{{M}}} \in \mathbb{R}^{3N \times 3N}$ consists of $N \times N$ blocks of size $3 \times 3$. A subblock $\mathbf{M}_{ij} \in \mathbb{R}^{3 \times 3}$ corresponding to the interaction between points $\mathbf{r}_i$ and $\mathbf{r}_j$ can be expressed as:
\begin{equation}
	\mathbf{M}_{ij} = f(\| \mathbf{r}_{ij} \|) \, \mathbf{I} + g(\| \mathbf{r}_{ij} \|) \, \hat{\mathbf{r}}_{ij} \otimes \hat{\mathbf{r}}_{ij} \label{eq:mobility_Mij}
\end{equation}
with $\mathbf{r}_{ij}=\mathbf{r}_i - \mathbf{r}_j$, while a hat in \autoref{eq:mobility_Mij} denotes a unit vector. $f(r)$ and $g(r)$ are scalar functions of distance. In this paper, the commonly used Rotne-Prager-Yamakawa (RPY) kernel~\cite{rotn69a} is employed for $f(r)$ and $g(r)$, leading to a symmetric positive-semidefinite mobility matrix $\mathbf{\mathcal{{M}}}$. The reader is referred to \cite{kall15a} for a detailed overview of the immersed boundary method that leads to \autoref{eq:system_Stokes}.

\begin{figure}[htbp]
	\begin{center}
	(a)~\hspace{-1.5em}\includegraphics*[width=0.475\linewidth,clip=true]{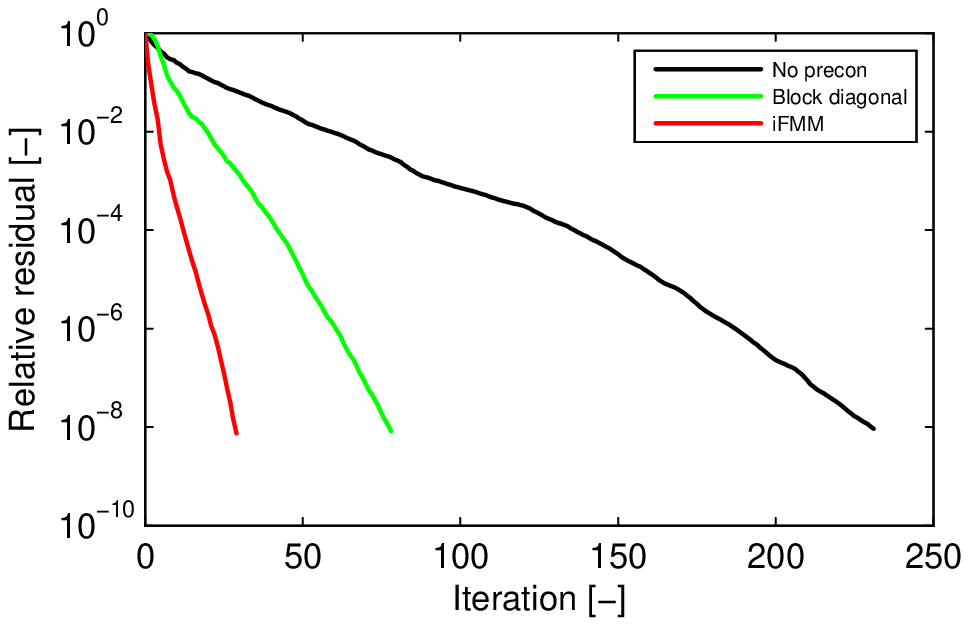}
	(b)~\hspace{-1.5em}\includegraphics*[width=0.475\linewidth,clip=true]{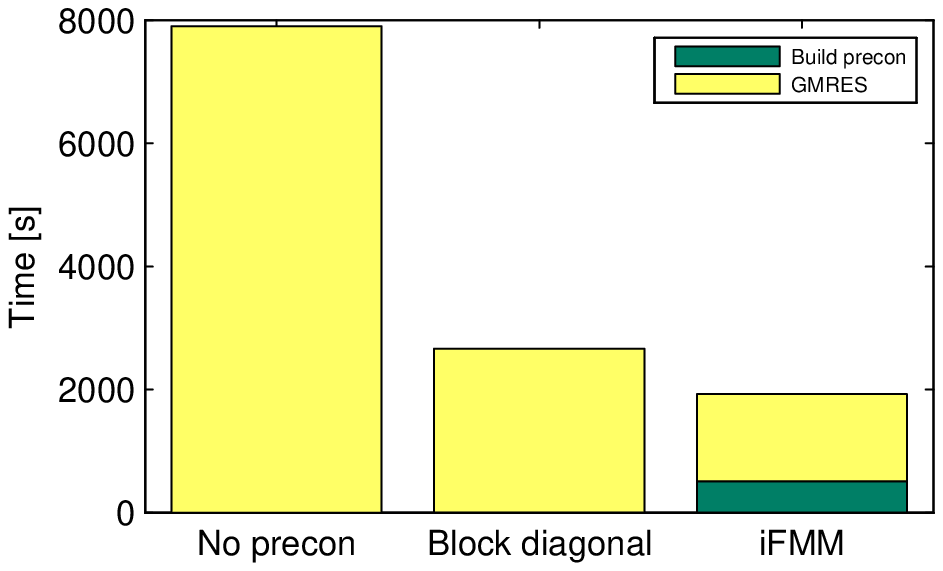}
	\end{center}
	\caption{(a)~Relative residual $\|\mathbf{U} -  \mathbf{\mathcal{M}} \hat{\mathbf{F}}\|_{2}/ \|\mathbf{U}\|_{2}$ as a function of the iteration step and (b)~total CPU time for a cubic lattice consisting of $4096$ spheres, each discretized with $42$ particles.}
\label{fig:residual_vs_iteration_4096_spheres}
\end{figure}

The dense system presented in \autoref{eq:system_Stokes} is solved with GMRES; the matrix-vector products are evaluated with an FMM-scheme~\cite{kall15a}. An arbitrary velocity vector~$\mathbf{U}$ is taken as right hand side. First, a regular cubic lattice consisting of $16 \times 16 \times 16$ spherical rigid bodies is considered. Each sphere is discretized with $42$ particles, resulting in a system of equations with $3 \times 42 \times 4096 = 516096$ unknowns. 
\autoref{fig:residual_vs_iteration_4096_spheres}(a) shows the relative residual as a function of the iteration step in the GMRES-algorithm, and this for three different preconditioning approaches. If no preconditioning is applied, convergence to the specified relative residual $\varepsilon = 10^{-8}$ is achieved after $231$ iterations. The use of a simple block-diagonal preconditioner (where the diagonal subblocks of size $126 \times 126$ correspond to the mobility matrices of the individual spheres) already leads to a significant reduction of the number of iterations and hence the overall computation time (see \autoref{fig:residual_vs_iteration_4096_spheres} and \autoref{tbl:GMRES_Stokes_lattice_tbl}). This is a very cheap preconditioning strategy (i.e., the time needed to construct the preconditioner is negligible) that performs very well in this particular case, as the subblocks in the block-diagonal matrix represent the mobilities of disjoint rigid bodies. Application of the IFMM as preconditioner (with $n=2$ Chebyshev nodes) results in a further reduction of the number of iterations. The overall speed-up remains rather limited, however, as the preconditioner's construction cost is considerably higher. It is nevertheless remarkable that the IFMM still results in a speed-up, given the fact that the block-diagonal preconditioner already performs very well.

\begin{table}[htbp]
\footnotesize
\centering
\caption{Overview of the results for various preconditioning strategies for a cubic lattice consisting of $4096$ spheres, each discretized with $42$ particles.}
\label{tbl:GMRES_Stokes_lattice_tbl}
\begin{tabular}{c| c c c }
	\toprule
    	&	\# iterations	&	Total CPU time [s]	& 	$\| \mathbf{U} -  \mathbf{\mathcal{M}} \hat{\mathbf{F}} \|_{2}/ \| \mathbf{U} \|_{2}$	  \\
	\midrule
	No precon 		& 231	& 7900 & $9.3 \times 10^{-9}$\\
	Block diagonal 	& 78	& 2662 & $8.3 \times 10^{-9}$\\
	IFMM ($n=2$)		& 29	& 1928 & $7.5 \times 10^{-9}$\\
	\bottomrule
\end{tabular}
\end{table}

\begin{figure}[hbtp]
	\begin{center}
	(a)~\hspace{-1.5em}\includegraphics*[width=0.475\linewidth,clip=true]{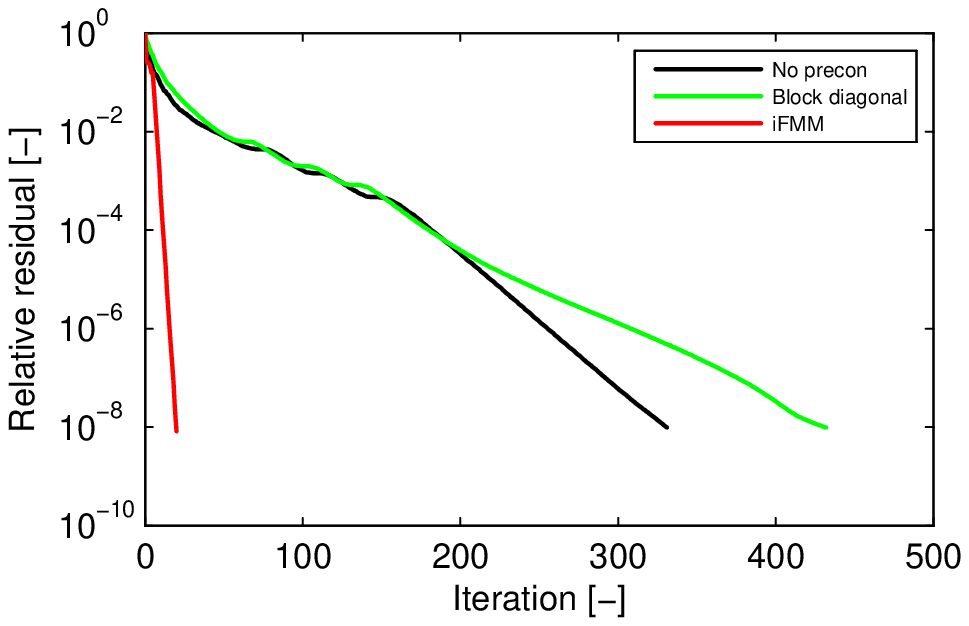}
	(b)~\hspace{-1.5em}\includegraphics*[width=0.475\linewidth,clip=true]{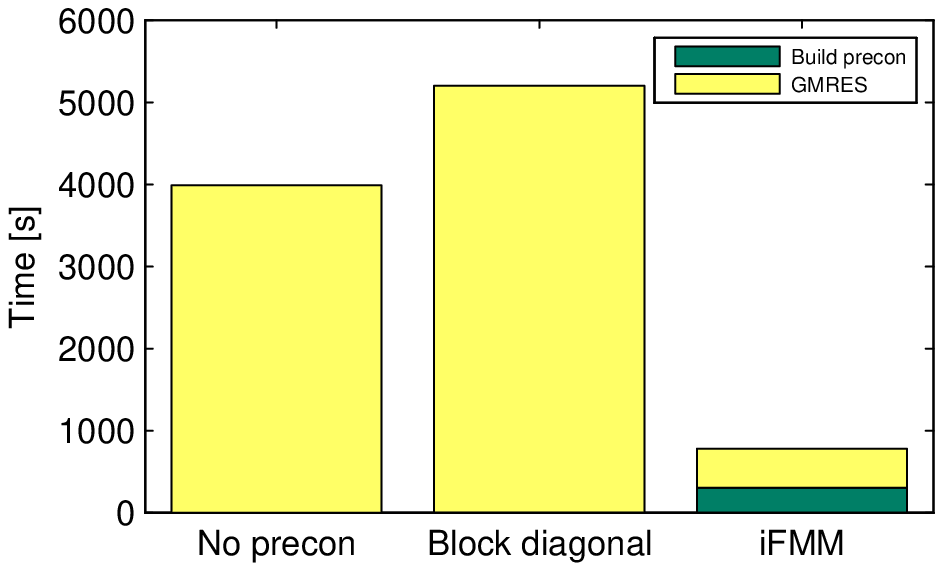}
	\end{center}
	\caption{(a)~Relative residual $\| \mathbf{U} -  \mathbf{\mathcal{M}} \hat{\mathbf{F}} \|_{2}/ \| \mathbf{U} \|_{2}$ as a function of the iteration step and (b)~total CPU time for six concentric spherical shells discretized with 54612 particles.}
\label{fig:residual_vs_iteration_concentricspheres}
\end{figure}

The second example involves six concentric spherical shells with a discretization ranging from $42$ particles for the smallest inner shell to $40962$ particles for the largest outer shell; the total number of particles equals $54612$. The mobility matrix~$\mathbf{\mathcal{{M}}}$ is more ill-conditioned compared to the previous example. Once again, GMRES is employed for solving \autoref{eq:system_Stokes}. 
\autoref{fig:residual_vs_iteration_concentricspheres}(a) and \autoref{tbl:GMRES_Stokes_concentricspheres_tbl} indicate that $331$ iterations are required if no preconditioner is incorporated. The use of a block-diagonal preconditioner (where diagonal subblocks of size $108 \times 108$ have been employed\footnote{Defining a block diagonal preconditioner in this case is rather arbitrary, as the diagonal subblocks do not correspond to the self-interactions of individual spheres anymore. The size of the subblocks was chosen as $108 \times 108$ to ensure that all diagonal subblocks are of equal size ($3 \times 54612/108=1517$), while keeping the trade-off between the size of the subblocks and the computational cost of assembling and applying the preconditioner in mind.}) does not result in a speed-up; it even deteriorates the convergence. While this was a successful preconditioning strategy in the previous example, such an approach is too simplistic to be effective in this more complicated test case, as the diagonal subblocks do not correspond any more to the mobility matrices of disjoint spheres. Application of the IFMM as preconditioner (with $n=3$ Chebyshev nodes), on the other hand, leads once more to a very fast convergence of the relative residual and consequently a considerable reduction of the computation time (see \autoref{fig:residual_vs_iteration_concentricspheres}(b)). Note that this is a relatively small problem compared to the previous example, as only $3 \times 54612 = 163836$ unknowns are considered. The relative speed-up of the IFMM is expected to be even more significant for larger problems.

\begin{table}[htbp]
\footnotesize
\centering
\caption{Overview of the results for various preconditioning strategies for six concentric spherical shells discretized with 54612 particles.}
\label{tbl:GMRES_Stokes_concentricspheres_tbl}
\begin{tabular}{c| c c c }
	\toprule
    	&	\# iterations	&	Total CPU time [s]	& 	$\| \mathbf{U} -  \mathbf{\mathcal{M}} \hat{\mathbf{F}} \|_{2}/ \| \mathbf{U} \|_{2}$	  \\
	\midrule
	No precon 		& 331	& 3988 & $9.9 \times 10^{-9}$\\
	Block diagonal 	& 432	& 5202 & $9.8 \times 10^{-9}$\\
	IFMM ($n=3$)		& 20	& 778  & $8.3 \times 10^{-9}$\\
	\bottomrule
\end{tabular}
\end{table}
%

\section{Conclusions and future work}\label{sec:Conclusions}
In this paper, the inverse fast multipole method has been presented as a fast direct solver for linear systems of equations involving dense $\mH^2$-matrices. The method is inexact, although an arbitrary high accuracy can be achieved by tuning the parameters appropriately. Because the error can be related to the decay of the singular values for low-rank fill-in blocks, we expect the error to decay roughly geometrically in many applications as a function of the rank. The overall computational cost scales like $\mO(r^2 N)$ where $r$ is the average rank at the leaf level. In terms of the tolerance criterion $\varepsilon$, the cost is $\mO(N \log(1/\varepsilon)^2)$.

The algorithm is based on two key ideas. First, the $\mH^2$-structure of the dense matrix is exploited to construct an extended sparse matrix. Second, fill-ins arising throughout the elimination of the latter are compressed and redirected if they correspond to well-separated clusters. As a result, the sparsity pattern is maintained (strictly speaking, there is a bound, independent of $N$, on the fill-in), leading to a fast solver that scales like $\mO(N)$ for 3D problems. Numerical examples have been presented to validate the linear scaling of the algorithm. It has furthermore been demonstrated that a low accuracy IFMM solver performs extremely well as a preconditioner for an iterative solver; using the IFMM as a preconditioner is often favorable, compared to its use as an accurate direct solver. Furthermore, an application related to the interaction of rigid bodies immersed in a Stokes fluid has been discussed to illustrate the method's applicability to challenging science and engineering problems.

An important aspect to consider in future work is the parallelization of the proposed methodology. The algorithm is expected to exhibit good parallel scalability, as each cluster only interacts with a limited number of other clusters (i.e., those in its neighbor and interaction list). Furthermore, alternative low-rank approximation schemes that might be more efficient than interpolation based on Chebyshev polynomials can be investigated. Finally, the order of elimination of the nodes can be optimized (at each level of the tree) to minimize the amount of fill-in and the number of (re)compressions, and could potentially also improve the numerical stability of the algorithm.

Although the examples presented in this paper only involve uniform octrees, the method can easily be applied to adaptive trees with a varying number of levels~\cite{pour15a,ying04a}.

\appendix
\section{Additional information - the IFMM as a direct solver}\label{app:}
\begin{figure}[htbp]
	\begin{center}
	\hspace{1.5em}~\includegraphics*[width=0.91\linewidth,clip=true,bb=55 94 438 116]{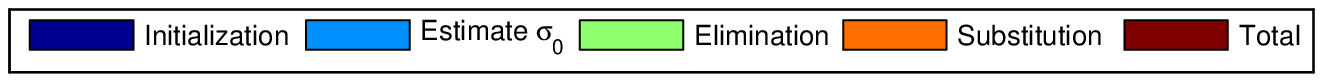}
	(a)~\hspace{-1.5em}\includegraphics*[width=0.475\linewidth,clip=true]{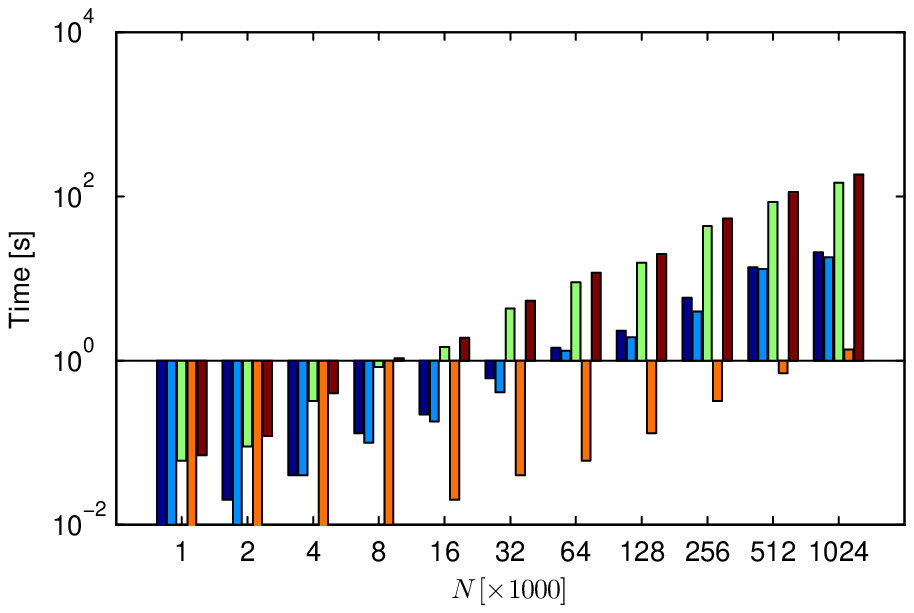}
	(b)~\hspace{-1.5em}\includegraphics*[width=0.475\linewidth,clip=true]{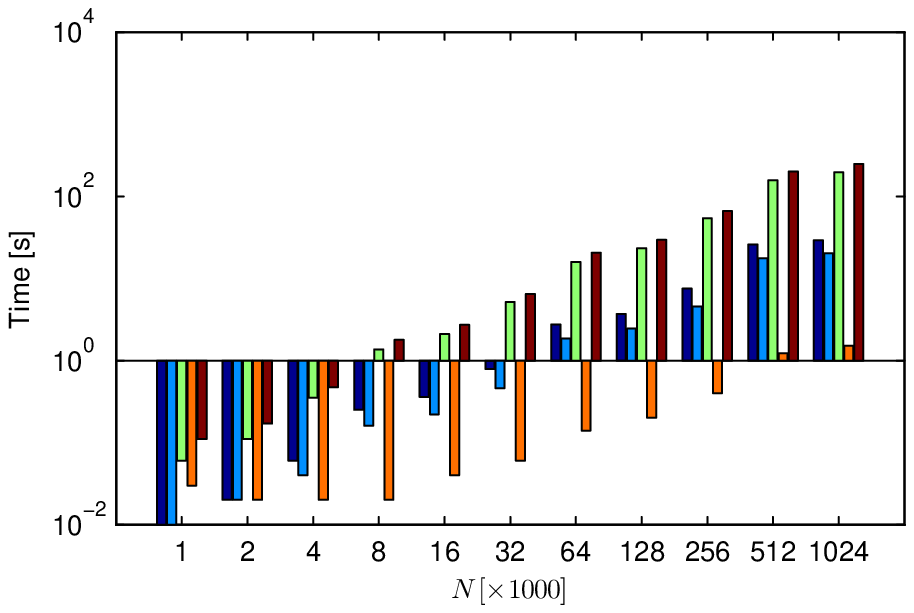}\\
	(c)~\hspace{-1.5em}\includegraphics*[width=0.475\linewidth,clip=true]{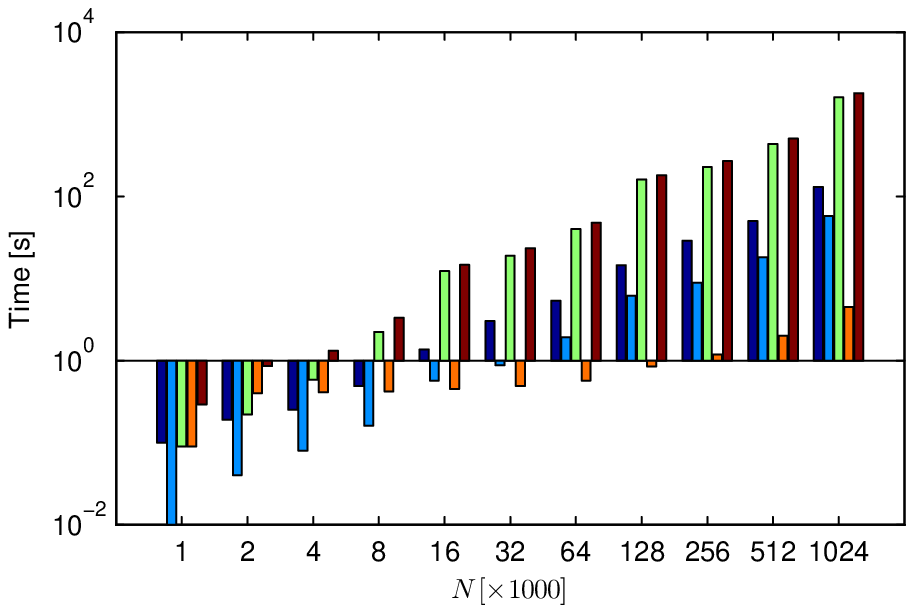}
	(d)~\hspace{-1.5em}\includegraphics*[width=0.475\linewidth,clip=true]{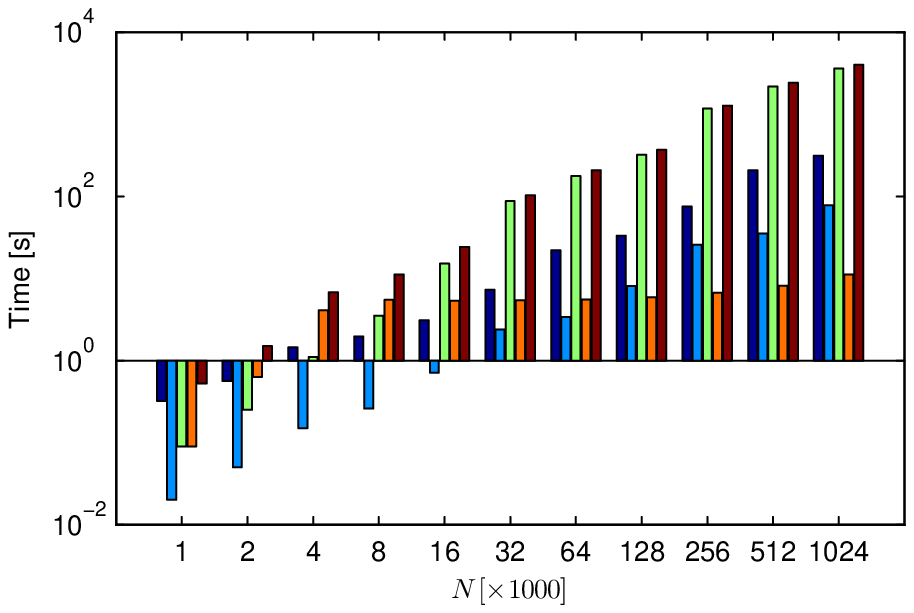}
	\end{center}
	\caption{Breakdown of the total CPU time for the IFMM as a direct solver in case of $N$ points randomly distributed in the cube $\left[-1,1\right]^3$. The number of Chebyshev nodes (in each dimension) used for the initial low-rank approximations corresponds to (a)~$n=1$, (b)~$n=2$, (c)~$n=3$, and (d)~$n=4$. From left to right, each bar represents: initialization, estimation of $\sigma_0(\mathbf{E})$, elimination, substitution, and total computation time (brown bars). The elimination step (light green bars) dominates the total time in most cases as expected.}
\label{fig:benchmark_ii_3D_random_paper_profiling_all}
\end{figure}

This appendix provides some additional data related to the benchmarks presented in \autoref{subsec:ifmm_direct}. More specifically, the computational cost for the IFMM as a direct solver, as depicted in \autoref{fig:benchmark_ii_3D_random}, is analyzed in detail in the following. \autoref{fig:benchmark_ii_3D_random_paper_profiling_all} shows a breakdown of the total CPU time into the time required for initializing the IFMM operators, estimating the largest singular value $\sigma_{0}(\mathbf{E})$ of the extended sparse matrix~$\mathbf{E}$, and performing the elimination and substitution. The number of Chebyshev nodes (in each dimension) is varied from $n=1$ to $n=4$; an increase of $n$ results in a higher overall computation time. It is clear from this figure that the elimination phase strongly dominates the total computation time, while the contribution of the other components, and especially of the substitution phase, is almost negligible.

\begin{figure}[htbp]
	\begin{center}
	\hspace{1.5em}~\includegraphics*[width=0.91\linewidth,clip=true,bb=82 97 409 116]{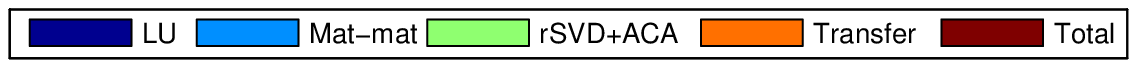}
	(a)~\hspace{-1.5em}\includegraphics*[width=0.475\linewidth,clip=true]{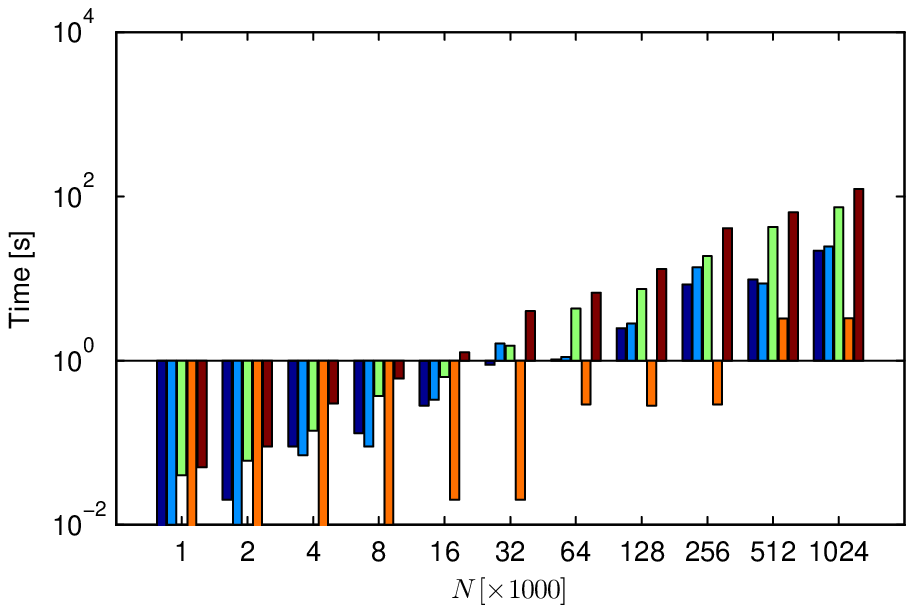}
	(b)~\hspace{-1.5em}\includegraphics*[width=0.475\linewidth,clip=true]{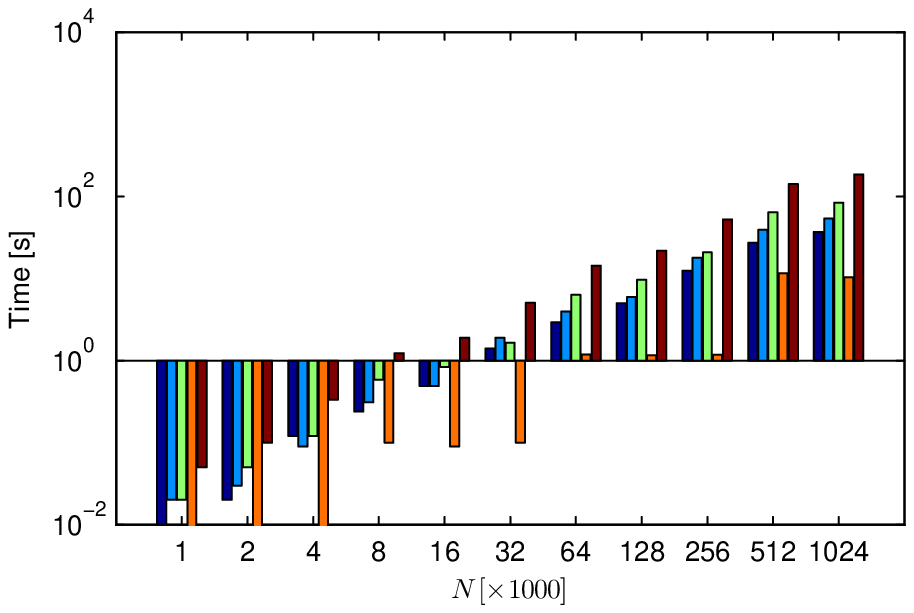}\\
	(c)~\hspace{-1.5em}\includegraphics*[width=0.475\linewidth,clip=true]{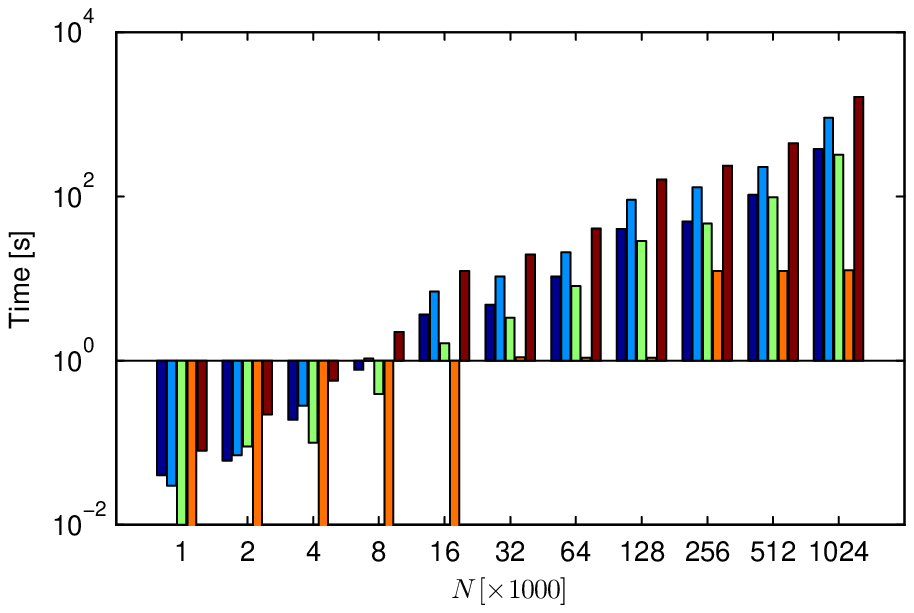}
	(d)~\hspace{-1.5em}\includegraphics*[width=0.475\linewidth,clip=true]{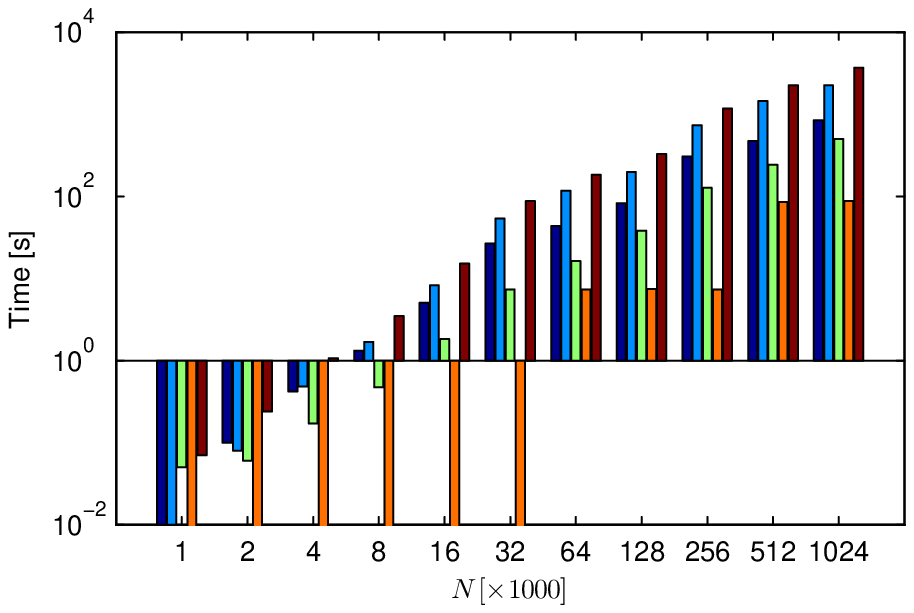}
	\end{center}
	\caption{Breakdown of the CPU time spent on different types of operations throughout the elimination phase of the IFMM algorithm. The number of Chebyshev nodes (in each dimension) used for the initial low-rank approximations corresponds to (a)~$n=1$, (b)~$n=2$, (c)~$n=3$, and (d)~$n=4$. From left to right, each bar represents: block LU factorization and triangular solves, matrix-matrix products (trailing matrix update in the block LU factorization), low-rank approximations, operator transfer (from child to parent), and total elimination time (brown bars). For small $n$, the low-rank approximations (light green bars) dominate, while for large $n$, the matrix-matrix products (light blue bars) dominate.}
\label{fig:benchmark_ii_3D_random_paper_profiling_elimination}
\end{figure}

A more detailed breakdown of the time spent on different types of operations throughout the elimination phase is provided in \autoref{fig:benchmark_ii_3D_random_paper_profiling_elimination}. Following operations are distinguished: performing LU-factorizations and -solves, executing matrix-matrix multiplications, obtaining low-rank approximations through rSVD and ACA, and transferring operators to the parent level (corresponding to \autoref{eq:U_merge}--\eqref{eq:K_merge}). The computation of low-rank approximations is dominant in case of $n=1$ and $n=2$, while matrix-matrix multiplications are the most time consuming operations for $n=3$ and $n=4$. In all cases, the contribution of transferring operators from one level to the parent level is significantly less important than the other components.

\begin{figure}[htbp]
	\begin{center}
		\includegraphics*[width=0.475\linewidth,clip=true]{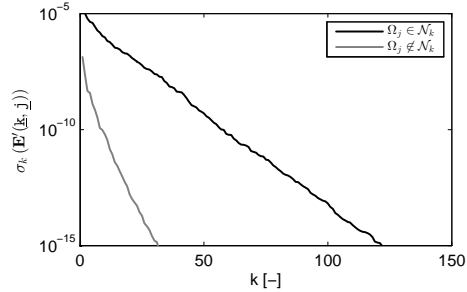}
	\end{center}
	\caption{Singular values $\sigma_k$ of fill-ins $\mathbf{E}^{\prime}( \underline{\texttt{k}}, \underline{\texttt{j}})$ arising between neighboring (black line) and well-separated (grey line) clusters $\Omega_{j}$ and $\Omega_{k}$.}
\label{fig:3D_randompoints_singularvalues}
\end{figure}

Finally, the different behavior of fill-ins arising between neighboring or well-separated clusters is illustrated. \autoref{fig:3D_randompoints_singularvalues} shows the singular values of some fill-ins $\mathbf{E}^{\prime}( \underline{\texttt{k}}, \underline{\texttt{j}})$ that are encountered throughout the elimination phase, both for neighboring ($\Omega_{j} \in \mathcal{N}_{k}$) and well-separated ($\Omega_{j} \not \in \mathcal{N}_{k}$) clusters $\Omega_{j}$ and $\Omega_{k}$. This figure confirms that the decay of singular values is significantly faster in the case of well-separated clusters, indicating that $\mathbf{E}^{\prime}( \underline{\texttt{k}}, \underline{\texttt{j}})$ can then indeed accurately be approximated using a low-rank representation.

%

\section{Additional information - the IFMM as a preconditioner}\label{app:leftright}

The results presented in \autoref{subsec:ifmm_precon} have been obtained by applying the IFMM as a right preconditioner in GMRES. This appendix compares these results to a left preconditioning approach.

The results of left and right preconditioning are compared in \autoref{tbl:GMRES_3D_randompoints_benchmark_ii_a_1e-3_tbl_leftright} and~\ref{tbl:GMRES_3D_randompoints_benchmark_ii_a_1e-2_tbl_leftright} for the two test cases considered in \autoref{subsec:ifmm_precon}. Note that the `relative residual' indicated in these tables corresponds to the actual residual $\| \mathbf{b} -  \mathbf{A} \widehat{\mathbf{x}} \|_{2}/ \| \mathbf{b} \|_{2}$ for right preconditioning, while it equals $\| \mathbf{P}^{-1}\mathbf{b} -  \mathbf{P}^{-1}\mathbf{A} \widehat{\mathbf{x}} \|_{2}/ \| \mathbf{P}^{-1} \mathbf{b} \|_{2}$ for left preconditioning (with $\mathbf{P}^{-1}$ representing the IFMM preconditioner).

For $d=10^{-3}$ (\autoref{tbl:GMRES_3D_randompoints_benchmark_ii_a_1e-3_tbl_leftright}), the effectiveness of the IFMM as left or right preconditioner is similar, both in terms of the number of iterations and the relative error $\|\mathbf{x} -  \widehat{\mathbf{x}} \|_{2}/ \| \mathbf{x} \|_{2}$. For $d=10^{-2}$ (\autoref{tbl:GMRES_3D_randompoints_benchmark_ii_a_1e-2_tbl_leftright}) and $n>1$, the left preconditioning approach needs more iterations to converge to the specified residual, but it leads to a significantly smaller relative error in the solution.

\begin{table}[htbp]
	\footnotesize
	\centering
	\caption{Overview of the results for left and right preconditioning for $d=10^{-3}$.}
	\label{tbl:GMRES_3D_randompoints_benchmark_ii_a_1e-3_tbl_leftright}
	\begin{tabular}{c| c c c}
		\toprule
	    	&	\# iterations	& 	relative residual  & $\|\mathbf{x} -  \widehat{\mathbf{x}} \|_{2}/ \| \mathbf{x} \|_{2}$\\
		\midrule
		$n=1$ (right) & 19& $6.1 \times 10^{-11}$ & $4.0 \times 10^{-11}$ \\
		$n=1$ (left)  & 19& $3.5 \times 10^{-11}$ & $4.5 \times 10^{-11}$ \\\midrule
		$n=2$ (right) & 6	& $1.5 \times 10^{-11}$ & $9.6 \times 10^{-12}$ \\
		$n=2$ (left)  & 6	& $9.1 \times 10^{-12}$ & $9.2 \times 10^{-12}$ \\\midrule
		$n=3$ (right) & 3	& $5.2 \times 10^{-11}$ & $2.9 \times 10^{-11}$ \\
		$n=3$ (left)  & 3	& $2.8 \times 10^{-11}$ & $2.8 \times 10^{-11}$ \\\midrule
		$n=4$ (right) & 2	& $2.7 \times 10^{-11}$ & $2.1 \times 10^{-11}$ \\
		$n=4$ (left)  & 2	& $2.0 \times 10^{-11}$ & $2.0 \times 10^{-11}$ \\
		\bottomrule
	\end{tabular}
	\end{table}

	\begin{table}[htbp]
	\footnotesize
	\centering
	\caption{Overview of the results for left and right preconditioning for $d=10^{-2}$.}
	\label{tbl:GMRES_3D_randompoints_benchmark_ii_a_1e-2_tbl_leftright}
	\begin{tabular}{c| c c c}
		\toprule
	    	&	\# iterations	& relative residual  &	$\|\mathbf{x} -  \widehat{\mathbf{x}} \|_{2}/ \| \mathbf{x} \|_{2}$\\
		\midrule
		$n=1$ (right)& 305& $9.7 \times 10^{-11}$ & $3.8 \times 10^{-8}$ \\
		$n=1$ (left) & 303& $5.4 \times 10^{-11}$ & $5.2 \times 10^{-8}$ \\\midrule
		$n=2$ (right)& 61 & $8.3 \times 10^{-11}$ & $2.2 \times 10^{-7}$ \\
		$n=2$ (left) & 62 & $9.3 \times 10^{-11}$ & $5.2 \times 10^{-9}$ \\\midrule
		$n=3$ (right)& 24 & $7.3 \times 10^{-11}$ & $4.6 \times 10^{-7}$ \\
		$n=3$ (left) & 27 & $3.8 \times 10^{-11}$ & $2.0 \times 10^{-9}$ \\\midrule
		$n=4$ (right)& 14 & $7.3 \times 10^{-11}$ & $7.9 \times 10^{-7}$ \\
		$n=4$ (left) & 18 & $3.8 \times 10^{-11}$ & $1.6 \times 10^{-9}$ \\
		\bottomrule
	\end{tabular}
	\end{table}

\section{Classification and relation between different methods}\label{classification}
We briefly present the different $\mH$-matrix formats and discuss connections between the algorithm presented in this paper and the fast solvers developed by Hackbusch et al.~\cite{bebe08a,borm03a,hack00a,hack02a,hack00b}.

$\mH$-matrices can be classified into four broad categories.

1) w$\mH$ with non-nested basis. Off-diagonal blocks that correspond to clusters that are well-separated are low-rank. The low-rank basis is non-nested, that is the low-rank basis of a node is unrelated to the basis of its children. This is the most general $\mH$-format. The matrix is generally speaking expanded in the form:
\[ \mathbf{A} = \mathbf{A}_s + \sum_{l=0}^{\ln N} \mathbf{U}_l \mathbf{\Sigma}_l \mathbf{V}_l^T \]
where $\mathbf{A}_s$ is a block sparse matrix.

2) w$\mH$ with nested basis. The low-rank basis of a node is a linear combination of the basis of its children nodes. In that case, the expansion can be written in recursive form:
\[ \mathbf{A} = \mathbf{A}_s + \mathbf{U} \mathbf{A}_{\mH} \mathbf{V}^T \]
where $\mathbf{A}_{\mH}$ is a w$\mH$ matrix of smaller size. FMM matrices~\cite{carrier1988fast,fong09a}, and the matrices discussed in this paper are of this form. This also corresponds to $\mH^2$-matrices~\cite{borm2007data,hackbusch2002data,hack02a}.

3) s$\mH$ with non-nested basis. All off-diagonal blocks are assumed to be low-rank. The matrix decomposition is:
\[ \mathbf{A} = \mathbf{A}_d + \sum_{l=0}^{\ln N} \mathbf{U}_l \mathbf{\Sigma}_l \mathbf{V}_l^T \]
where $\mathbf{A}_d$ is a block diagonal matrix. This is the same as the HODLR matrix format~\cite{ambi13b}.

4) s$\mH$ with nested basis:
\[ \mathbf{A} = \mathbf{A}_d + \mathbf{U} \mathbf{A}_{\mH} \mathbf{V}^T \]
where $\mathbf{A}_{\mH}$ is an s$\mH$ matrix. This is the HSS format~\cite{chan06a,shen07b}.

\bigskip

We now discuss the connection between our algorithm and Hackbusch et al.~\cite{bebe08a,borm03a,hack02a,hack00a,hack00b}. Both methods use a tree decomposition of the matrix in order to identify low-rank blocks. Although the Hackbusch algorithm was written differently, this algorithm can be derived using the extended sparse matrix introduced in \autoref{subsec:extended_sparsification}.

The key difference with this paper is that Hackbusch's algorithm performs the LU factorization using a depth-first in-order (symmetric) traversal. In contrast, our method uses a breadth-first or level-by-level traversal, starting from the leaf level and moving up to the root.

This can be illustrated on a simple example for concreteness. Consider the simplest $\mH$-matrix:
\[
\mathbf{A} = \begin{pmatrix}
\mathbf{A}_{00} & \mathbf{U}_0 \mathbf{V}_1^T \\
\mathbf{U}_1 \mathbf{V}_0^T & \mathbf{A}_{11}
\end{pmatrix}
\]
The extended sparse format is:
\[
\mathbf{E} = \begin{pmatrix}
\mathbf{A}_{00} & & \mathbf{U}_0 \\
& \mathbf{A}_{11} & & \mathbf{U}_1 \\
\mathbf{V}_0^T & & & -\mathbf{I} \\
& \mathbf{V}_1^T & -\mathbf{I}
\end{pmatrix}
\]
Hackbusch's algorithm is based on a block LU factorization. Such a factorization is obtained if, in $\mathbf{E}$, we eliminate first $\mathbf{A}_{00}$, then the auxiliary variables, then $\mathbf{A}_{11}$. Step-by-step, let's eliminate $\mathbf{A}_{00}$:
\[
\begin{pmatrix}
\mathbf{A}_{11} & & \mathbf{U}_1 \\
& -\mathbf{V}_0^T \mathbf{A}_{00}^{-1} \mathbf{U}_0  & -\mathbf{I} \\
\mathbf{V}_1^T & -\mathbf{I}
\end{pmatrix}
\]
If we eliminate pivots 2 and 3 (before $\mathbf{A}_{11}$), we simply get:
\[
\mathbf{A}_{11} - \mathbf{U}_1 \mathbf{V}_0^T \mathbf{A}_{00}^{-1} \mathbf{U}_0 \mathbf{V}_1^T
\]
This is the expected Schur complement from a block LU factorization.

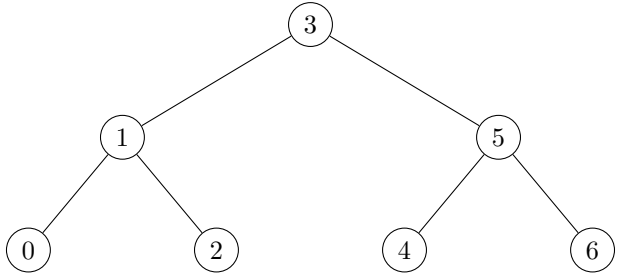
\begin{figure}[htbp]
\tikzset{
  tn/.style = {circle, draw=black, align=center, inner sep=3pt, text centered}
}

\begin{center}
\begin{tikzpicture}[>=stealth',level/.style={sibling distance = 5cm/#1,
  level distance = 1.5cm}] 
\node [tn] {3}
    child{  node [tn] {1}
              child{ node [tn] {0} }
              child{ node [tn] {2} }
         }
    child{  node [tn] {5}
              child{ node [tn] {4} }
              child{ node [tn] {6} }
         }         
; 
\end{tikzpicture}
\end{center}
\caption{Tree node elimination order following Hackbusch's algorithm. The traversal is a depth-first in-order tree traversal.}
\label{hack_elim}
\end{figure}

If we represent $\mathbf{E}$ using a tree, for a 3-level $\mH$-matrix, this elimination corresponds to the order shown in \autoref{hack_elim}.

In contrast, our algorithm follows a level-by-level elimination. Going back to $\mathbf{E}$, we first eliminate $\mathbf{A}_{00}$, then $\mathbf{A}_{11}$. This leads to a new linear system involving the auxiliary variables only:
\[
\begin{pmatrix}
-\mathbf{V}_0^T \mathbf{A}_{00}^{-1} \mathbf{U}_0 & -\mathbf{I} \\
-\mathbf{I} & -\mathbf{V}_1^T \mathbf{A}_{11}^{-1} \mathbf{U}_1
\end{pmatrix}
\]

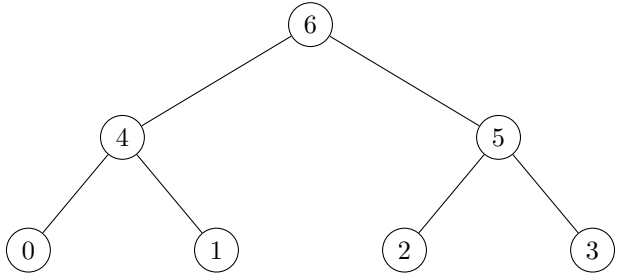
\begin{figure}[htbp]
\tikzset{
  tn/.style = {circle, draw=black, align=center, inner sep=3pt, text centered}
}

\begin{center}
\begin{tikzpicture}[>=stealth',level/.style={sibling distance = 5cm/#1,
  level distance = 1.5cm}] 
\node [tn] {6}
    child{  node [tn] {4}
              child{ node [tn] {0} }
              child{ node [tn] {1} }
         }
    child{  node [tn] {5}
              child{ node [tn] {2} }
              child{ node [tn] {3} }
         }         
; 
\end{tikzpicture}
\end{center}
\caption{Tree node elimination order following this paper's algorithm. The traversal is a level-by-level (breadth-first) traversal.}
\label{this_paper_elim}
\end{figure}

As a result, in terms of traversal ordering, this paper follows a level-by-level elimination. This is shown in \autoref{this_paper_elim}.

\bigskip \indent
{\bf{Acknowledgments}} The authors would like to thank Prof.\ Aleksandar Donev and Dr.\ Floren Balboa (Courant Institute of Mathematical Sciences, New York University) for their assistance in obtaining the results presented in \autoref{sec:Application}. 

The first author is a post-doctoral fellow of the Research Foundation Flanders (FWO) and a Francqui Foundation fellow of the Belgian American Educational Foundation (BAEF). The financial support is gratefully acknowledged. 

Funding from the ``Army High Performance Computing Research Center'' (AHPCRC), sponsored by the U.S. Army Research Laboratory under contract No.~W911NF-07-2-0027, at Stanford, supported in part this research.

\bibliographystyle{siam}
\bibliography{abbreviations,references}

\end{document}